\documentclass[12pt]{article}
\usepackage{amssymb,amsmath,latexsym,amsthm}
\usepackage[utf8]{inputenc}
\usepackage[T1]{fontenc}
\usepackage[english]{babel}
\usepackage{esint}
\usepackage{amssymb, amsthm, bbm, amsmath, latexsym, color}
\usepackage{enumerate, mathrsfs, hyphenat, graphicx}
\usepackage[font=sf, labelfont=sf]{caption}
\newcommand{\dom}{\mathcal{D}}
\usepackage{thmtools, thm-restate}
\usepackage{thm-patch, thm-kv, thm-autoref}
\usepackage[colorlinks=true, urlcolor=blue, linkcolor=blue, citecolor=black]{hyperref}
\theoremstyle{plain}
\declaretheorem[title=Theorem, parent=section]{theorem}
\declaretheorem[title=Lemma,sibling=theorem]{lemma}
\declaretheorem[title=Corollary,sibling=theorem]{corollary}
\declaretheorem[title=Proposition,sibling=theorem]{proposition}
\declaretheorem[title=Remark,sibling=theorem]{remark}
\theoremstyle{definition}
\declaretheorem[title=Definition,sibling=theorem]{definition}
\declaretheorem[title=Example,sibling=theorem]{example}

\numberwithin{equation}{section}
\newcommand{\eps}{\varepsilon}
\newcommand{\R}{  \mathbb{R}}

\addto\extrasenglish{

}

\usepackage[backgroundcolor=none, bordercolor=blue]{todonotes}
\makeatletter
\providecommand\@dotsep{5}
\renewcommand{\listoftodos}[1][\@todonotes@todolistname]{%
  \@starttoc{tdo}{#1}}
\makeatother

\def\bea{\begin{align*}}
\def\eea{\end{align*}}
\def\bee{\begin{equation}}
\def\eee{\end{equation}}

\def\LL{{\cal L}}
\def\DD{{\cal D}}
\def\AA{{\cal A}}

 \def\E {{\mathbb E}} 
 \def\bH {{\mathbb H}}

\def\P {{\mathbb P}}  \def\bR {{\mathbb R}}

\newcommand{\X}{\mathfrak{X}}
\def\EE{{\cal E}}
\def\FF{{\cal F}}
\def\sG {{\mathcal G}}
\def\sH {{\mathcal H}}
\def\sA {{\mathcal A}}
\def\sL{{\mathcal L}}
\def\U{{\mathcal U}}
\def\V{{\mathcal W}}
\def\1{{\mathbbm{1}}}
\def\nn{\nonumber}

\def\qed{{\hfill $\Box$ \bigskip}}
\def\eps{\varepsilon}
\def\wt{\widetilde}
\def\wh{\widehat}

\def\tx{{\wt{x}}}
\def\ty{{\wt{y}}}

\def\tu{{\wt{u}}}

\def\pf{\noindent{\bf Proof.} }

\allowdisplaybreaks

\begin{document}

\author{
{\bf Soobin Cho}\thanks{This work was  supported by the National Research Foundation of
Korea(NRF) grant funded by the Korea government(MSIP) (No. NRF-2015R1A4A1041675)}
\quad
{\bf Panki Kim}\thanks{This work was  supported by the National Research Foundation of
Korea(NRF) grant funded by the Korea government(MSIP) (No. NRF-2015R1A4A1041675)
}
\quad {\bf Renming Song\thanks{Research supported in part by a grant from
the Simons Foundation (\#429343, Renming Song)}} \quad and
\quad {\bf Zoran Vondra\v{c}ek\thanks{Research supported in part by the Croatian Science Foundation under the project 4197}}
}

\title{Factorization and estimates of Dirichlet heat kernels for non-local operators with critical killings} 
 \date{}

\maketitle
 \begin{abstract}
In this paper we discuss non-local operators with 
killing potentials, which may not be in the standard Kato class. 
We first discuss factorization of their Dirichlet heat kernels in metric measure spaces.
Then we establish explicit estimates of the Dirichlet heat kernels under critical killings 
 in $C^{1,1}$ open subsets of $\R^d$ or in $\R^d\setminus\{0\}$.
The decay rates of our explicit estimates come from the values of the multiplicative constants 
in the killing potentials.
Our method also provides an alternative and unified proof of the main results of 
\cite{CKS10a, CKS10b, CKS-AOP}.
\end{abstract}



\smallskip

\section{Introduction}\label{s:int}
Stability of  Dirichlet heat kernel estimates under certain Feynman-Kac transforms
was studied in the recent paper \cite{CKS15}.
To be more precise, let $X$ be a Hunt process on a Borel set $D\subset \R^d$ that admits a jointly continuous transition density $p_D(t,x,y)$ with respect to the Lebesgue measure. 
Let $\alpha\in (0,2)$ and $\gamma\in [0, \alpha\wedge d)$,  and define
$$
q_{\gamma}(t,x,y):=\left(1\wedge \frac{\delta_D(x)}{t^{1/\alpha}}\right)^{\gamma}
\left(1\wedge \frac{\delta_D(y)}{t^{1/\alpha}}\right)^{\gamma}\left(t^{-d/\alpha}\wedge \frac{t}{|x-y|^{d+\alpha}}\right),
$$
where $\delta_D(x)$ denotes the distance between $x\in D$ and $D^c$.
Assume that $p_D(t,x,y)$ is comparable to $q_{\gamma}(t,x,y)$ for $(t,x,y)\in (0,1]\times D\times D$. Examples of processes satisfying this assumption include killed symmetric stable processes in $C^{1,1}$ open sets $D$ 
(with $\gamma=\alpha/2$, cf. \cite{CKS10a}), 
and, when $\alpha\in(1,2)$, censored $\alpha$-stable processes in any $C^{1,1}$ 
 open sets $D$ (with $\gamma=\alpha-1$, cf. \cite{CKS10b}).
Consider the following Feynman-Kac transform:
$$
T_t f(x)=\E_x\left[\exp\left(-A_t\right)f(X_t)\right],
$$
where $A$ is a continuous additive functional of $X$ with  Revuz measure 
$\mu$.
If $\LL$ denotes the $L^2$-infinitesimal generator of $X$, then, informally, the semigroup $(T_t)$ has the $L^2$-infinitesimal generator 
$\AA f(x):=(\LL-\mu)f(x)$.
Under the assumption that $\mu$ belongs to some appropriate Kato class, 
one of the main results of \cite{CKS15} implies 
 that the semigroup $(T_t)$ admits a continuous density $q^D(t,x,y)$
which is comparable to $q_{\gamma}(t,x,y)$ for all $(t,x,y)\in (0,1]\times D\times D$. Hence a Kato class perturbation preserves the Dirichlet heat kernel estimates
and is in this sense subcritical. Related results on the stability of the 
 heat kernel estimates
 (without boundary condition) under Kato class perturbations were obtained earlier in 
 \cite{BHJ, S06, W08}.

Kato class perturbations of the Laplacian have been studied earlier and more thoroughly, e.g. 
\cite{AS82, BM90, KS07, S82}, 
with the same conclusion that Kato class perturbations preserve the 
(Dirichlet) heat kernel estimates. 
It is well known since \cite{BG84} that,
in the case of the Laplacian in the whole space, 
the inverse square potential $\kappa(x)=c|x|^{-2}$ 
is critical,
and, in the case of the Dirichlet Laplacian in a domain $D$, the potential $\kappa(x)=c\delta_D(x)^{-2}$ is critical.
Criticality of the potentials above can be explained as follows.
In the whole space case, consider the symmetric form
$$
\EE(u,u)=\int_{\R^d} \left(|\nabla u(x)|^2+\frac{cu(x)^2}{|x|^2}\right)dx, 
\quad u\in C_0^{\infty}(\R^d)\, .
$$
Due to Hardy's inequality, this form is  non-negative 
if and only if $c\ge -((d-2)/2)^2$ (here $d\ge 3$),
and the operator $\Delta-c|x|^{-2}$ is defined as 
the $L^2$-generator of the form above.
In the case of a bounded smooth domain $D$, consider the symmetric form
$$
\EE(u,u)=\int_D \left(|\nabla u(x)|^2+\frac{cu(x)^2}{\delta_D(x)^2}\right)dx, 
\quad u\in C_0^{\infty}(D)\, .
$$
Again due to Hardy's inequality, this form is  non-negative 
if and only if $c\ge -1/4$,
and the operator $\Delta-c\delta_D(x)^{-2}$ is defined as the
$L^2$-generator of the form above. 
In both cases above, when $c<0$, 
the potential $\kappa$ above can be interpreted as creation, and, 
when $c>0$, the potential $\kappa$ can be interpreted as killing.
Note that in both cases,
even when $c>0$, the potential $\kappa$ above does not belong to the Kato class.
There exists a large body of literature on the heat kernel estimates of 
critical perturbations of the (Dirichlet) Laplacian,
e.g.,~\cite{BFT04, FMT07, IKO17, MS03, MS04, MT07}. 
In this paper we use probabilistic methods to study 
sharp two-sided heat kernel estimates for  critical perturbations of 
the fractional Laplacian in a smooth domain $D$, 
as well as the fractional Laplacian in $\R^d$.
When the potential involves both killing and creation, there is no Markov process associated
with the corresponding Schr\"odinger type operator. Since our argument depends crucially on
properties of Markov processes, we will only deal with killing type potentials. 
To deal with
critical perturbations involving creation, different methods 
are needed, 
cf.~\cite{BDK, BGJP}.

Let us describe some of our results in more detail. Let $D\subset \R^d$, $d\ge 2$, be a $C^{1,1}$ open set and let $X=(X_t, \P_x)$ be the reflected $\alpha$-stable process in $\overline{D}$, $\alpha\in (0,2)$, cf.~\cite{BBC03}. When $\alpha\in (0,1]$, we denote by $X^D$ the process $X$ restricted to $D$ ($X$ does not hit $\partial D$ in this case), while for $\alpha\in (1,2)$, $X^D$ denotes the censored $\alpha$-stable process in $D$. 
In this paragraph, we will only describe some of our results in the case $\alpha\in (1,2)$. So in the rest of this paragraph, we will assume that $\alpha\in (1,2)$. Similar results are also valid
in the case $\alpha\in (0,1]$.
Consider the potential $\kappa(x)=c_1 \delta_D(x)^{-\alpha}$ where 
$c_1\in [0,\infty)$ (the family of potentials we study is in fact larger, 
see \eqref{e:C1} for the full picture). 
The criticality of this type of potentials can also be interpreted using the fractional Hardy 
inequality in, for instance,  \cite{CS-TMJ, BD}.
Let 
$$
T_t f(x):=\E_x\left[e^{-\int_0^t \kappa(X_s^D)ds} f(X_t^D)\right], \quad x\in D, t>0,
$$
be the Feynman-Kac semigroup of $X^D$ via the multiplicative functional $e^{-\int_0^t \kappa(X_s^D)ds}$. Alternatively, we can think of $(T_t)$ as the semigroup corresponding to the Schr\"{o}dinger operator $\LL-\kappa$, where $\LL$ is the generator of $X^D$. We show that the semigroup $(T_t)$ admits a continuous density $q^D(t,x,y)$ 
and that there exists $p\in [\alpha-1, \alpha)$ depending on the constant $c_1$ such that $q^D(t,x,y)$
is comparable to 
\begin{equation}\label{e:estimates-intro}
\left(1\wedge \frac{\delta_D(x)}{t^{1/\alpha}}\right)^{p}\left(1\wedge \frac{\delta_D(y)}{t^{1/\alpha}}\right)^{p}\left(t^{-d/\alpha}\wedge \frac{t}{|x-y|^{d+\alpha}}\right),
\end{equation}
for $(t,x,y)\in (0,1]\times D\times D$, see Theorem \ref{t:DHKE1}. Moreover, the mapping $[0,\infty)\ni c_1\mapsto p\in [\alpha-1, \alpha)$ is one-to-one and onto. If $c_1=0$ (no killing), then $p=\alpha-1$ and we recover the heat kernel estimates for the censored $\alpha$-stable  process obtained in \cite{CKS10b}. When $p=\alpha/2$, the semigroup $(T_t)$ corresponds to the $\alpha$-stable process killed upon exiting $D$ and we recover the heat kernel estimates from \cite{CKS10a}. The novelty of our approach is that by changing the constant $c_1$ in the potential $\kappa$ we can obtain a whole spectrum of boundary behaviors of the heat kernel. In particular, we construct semigroups whose heat kernels satisfy the assumptions from \cite{CKS15}. 
On the other hand, the interior estimates
$t^{-d/\alpha}\wedge t|x-y|^{-d-\alpha}$, which correspond to the reflected $\alpha$-stable process,  
remain unchanged.

Besides reflected $\alpha$-stable processes on $\overline{D}$,
we can also consider processes which are lower order perturbations of reflected $\alpha$-stable processes.
A typical example is the process $X^{\beta}$ on $\overline{D}$ with Dirichlet form of the type
$$
\EE^{\beta}(u,u)=\frac12 \int_D\int_D (u(x)-u(y))^2 \left(\frac{\sA(d,-\alpha)}{|x-y|^{d+\alpha}}+\frac{b}{|x-y|^{d+\beta}}\right)dxdy,
$$
where $D$ is a bounded $C^{1, 1}$ open set, $\beta<\alpha<2$, $\sA(d,-\alpha)=\alpha 2^{\alpha-1} \pi^{-d/2} \Gamma((d+\alpha)/2)\Gamma(1-\alpha/2)^{-1}$  and $b$ is a non-negative constant (see Subsection \ref{ss-nl} for the more general setup, 
in particular, $b$ need not be non-negative).
 Let $X^{\beta, D}$ denote the process $X^{\beta}$ killed upon exiting $D$. We prove that the Feynman-Kac semigroup of $X^{\beta, D}$ via the multiplicative functional  $e^{-\int_0^t \kappa(X_s^D)ds}$ has 
a continuous transition density comparable 
 to \eqref{e:estimates-intro}, where the critical potential 
  $\kappa(x)=c_1\delta_D(x)^{-\alpha}$ is the same as above.

Our final result concerns the isotropic $\alpha$-stable process $Z$ in $\R^d$, $d\ge 2$, and the singular potential $\kappa(x)=c_1|x|^{-\alpha}$, $c_1>0$ (again, we in fact consider more general potentials -- see \eqref{e:kappainf}). We show that the Feynman-Kac semigroup of $Z$ via the multiplicative functional $e^{-\int_0^t \kappa(Z_s)ds}$ admits a continuous density comparable to
$$
\left(1\wedge \frac{|x|}{t^{1/ \alpha}}\right)^p \left(1\wedge \frac{|y|}{t^{1/ \alpha}}\right)^p \left(t^{-d/\alpha}\wedge \frac{t}{|x-y|^{d+\alpha}}\right),
$$
with $p\in (0,\alpha)$ depending on $c_1$, see Theorem \ref{t:DHKE3}. 
For related results, cf.~\cite[Theorem 1.1]{BGJP}  and \cite[Theorem 1.1]{JW}.

\medskip
Organization of the paper: the paper is divided into two major parts and an appendix. The first part is Section 2 and the setup is quite general there. 
We consider a  Hunt process $X$ on a locally compact separable metric space $(\X, \rho)$. 
The process $X$ is not necessarily symmetric and may not be conservative.
We assume that the process $X$ is in strong duality (with respect to a Radon measure $m$ with full support) with another Hunt process $\wh{X}$. 
We further assume that both  
$X$ and $\wh{X}$ are 
Feller and strongly Feller, 
and that their semigroups admit a strictly positive and jointly continuous transition density $p(t,x,y)$. The main assumption is that $p(t,x,y)$ is, for small times, comparable to the function $\wt{q}(t,x,y)$ 
defined in \eqref{e:tildeq} in
terms of the volume of balls in $(\X, \rho)$ and a strictly increasing function $\Phi$ satisfying a weak scaling condition \eqref{e:scale_Phi}.
In Subsection \ref{s:pre} we argue that $X$ and $\wh{X}$ satisfy the scale invariant parabolic Harnack inequality with explicit scaling in terms of $\Phi$ and use this to obtain interior lower bound on the transition density $p_D(t,x,y)$ of the process $X^D$ - the process $X$ killed upon exiting an open subset $D$  of $\X$.
Subsection \ref{s:3p} contains the definition of the class $\mathbf{K}_T(D)$, $T>0$, of 
possibly critical smooth measures $\mu$ on $D$, cf.~Definition \ref{d:KT}. Using the positive additive functional $(A_t^{\mu})$ of $X^D$ with Revuz measure $\mu$, we define the Feynman-Kac semigroup of 
$X^D$ associated with $\mu$:
$$
T_t^{\mu, D}f(x)=\E_x\left[\exp(-A_t^{\mu})f(X_t^D)\right], \quad t\ge 0, x\in D\, . 
$$
The Hunt process corresponding to $(T_t^{\mu,D})$ is denoted by $Y$. 
We analogously define the dual semigroup $(\wh{T}_t^{\mu, D})$ and denote the corresponding Hunt process by $\wh{Y}$.
We argue that $(T_t^{\mu,D})$ has a transition density (with respect to the measure $m$) $q^D(t,x,y)$ and that there exists $C_0>0$ such that $q^D(t,x,y)\le C_0\wt{q}(t,x,y)$ for small $t$. The argument relies on the 3P inequality for $\wt{q}(t,x,y)$ proved in Lemma \ref{l:3P}.
In Subsection \ref{s:ie} we prove some interior estimates for the transition density $q^U(t,x,y)$, where $U$ is an open subset of  $D$. Two examples of critical potentials are given in Subsection \ref{s:ex}. Finally, in Subsection \ref{s:factorization}, we show that factorization of the transition density $q^D(t,x,y)$ in $\kappa$-fat open set $D$ holds true. 
The result is proved in Theorem \ref{t:f1} and states that for small
 time $t$, $q^D(t,x,y)$ is comparable to $\P_x(\zeta >t)\wh{\P}_y(\wh{\zeta}>t)\wt{q}(t,x,y)$, $x,y\in D$. Here $\zeta$ and $\wh{\zeta}$ are 
the lifetimes of $Y$ and $\wh{Y}$ respectively.
 We note that this is a quite general result and, besides critical perturbations, includes also subcritical perturbations (or no perturbation at all). 
 Criticality and subcriticality of the perturbation are hidden 
 in the tail behavior of the lifetimes, $\P_x(\zeta >t)$ and $\wh{\P}_y(\wh{\zeta}>t)$. 
 To prove Theorem \ref{t:f1}, we follow the ideas in the proof of 
  \cite[Theorem 1.3]{CKS-PLMS}, but we use Assumption {\bf U}, cf.~Subsection 
 \ref{s:factorization}, instead of  the boundary Harnack principle. 
Approximate factorization of heat kernels involving tails of lifetimes can be traced back to \cite{Var, BG}.
 If one can get explicit two-sided estimates on the survival probabilities, then one can combine them
with the approximate factorization to get explicit two-sided estimates on the heat kernel.
This is the strategy employed in \cite{BGR, CKS-PLMS}. We will also use this strategy in
Section \ref{s:3}, and
 as a by-product, 
 give  an alternative and unified proof of the main results of 
\cite{CKS10a, CKS10b, CKS-AOP}.

The second part is Section \ref{s:3}. In this section
we assume that  $\X$ is either the closure of a $C^{1,1}$ open subset $D$ of $\R^d$ 
or $\R^d$ itself, $d\ge 2$, and we assume that 
the underlying process 
$X$ is either a reflected $\alpha$-stable(-like) process  on $\overline {D}$ (or a non-local perturbation of it), or an $\alpha$-stable process in $\R^d$ 
(or a drift perturbation of it).
The critical potentials have been already described above and are essentially of the form either  $
c_1\delta_D(x)^{-\alpha}$ or $c_1|x|^{-\alpha}$.
The goal of this section is to estimate the tail of the lifetime $\P_x(\zeta >t)$ in terms of
 $\delta_D(x)$ and $|x|$ respectively. 
 Then, as was done in \cite{BGR, CKS-PLMS}, together 
 with the factorization obtained in Theorem \ref{t:f1}, this gives sharp two-sided estimates of the transition density of the Feynman-Kac semigroup. 
The main step in obtaining the estimates of the lifetime consists of finding appropriate superharmonic and subharmonic function for the process $Y$. 
This relies on quite detailed computations of the generator of $Y$ acting on some appropriate functions. Although similar methods have been already employed in some previous works, our calculations are quite involved and delicate.
\ Section \ref{s:3}  also  provides an alternative and unified proof of the main results of 
\cite{CKS10a, CKS10b, CKS-AOP}.

The paper ends with an appendix devoted to a result about a continuous additive  functional for killed non-symmetric process.

Notation: We will use the symbol ``$:=$'' to denote a definition, 
which is read as ``is defined to be.''
In this paper, 
for $a,b\in \R$, we denote $a\wedge b:=\min\{a,b\}$ and $a\vee b:=\max\{a,b\}$.  We also use the convention $0^{-1}=+\infty$. 
For two non-negative functions $f$ and $g$, the notation $f\asymp g$ means that there are strictly positive constants 
$c_1$ and $c_2$ such that $c_1g(x)\leq f (x)\leq c_2 g(x)$ in the common domain of the definition of $f$ and $g$.

Letters with subscripts $r_i$, $R_i$, $A_i$, $C_i$, $i=0,1,2,  \dots$, denote constants 
that will be fixed throughout the paper. 
Lower case letters $c$'s without subscripts denote strictly positive
constants  whose values
are unimportant and which  may change even within a line, while values of lower case letters with subscripts
$c_i, i=0,1,2,  \dots$, are fixed in each proof, 
and the labeling of these constants starts anew in each proof.
$c_i=c_i(a,b,c,\ldots)$, $i=0,1,2,  \dots$, denote  constants depending on $a, b, c, \ldots$.
The dependence on the dimension $d \ge 1$
may not be mentioned explicitly.

For  a function space $\bH(U)$
on an open set $U$ in $\X$, we let   
$\bH_c(U):=\{f\in\bH(U): f \mbox{ has  compact support}\},$
$\bH_0(U):=\{f\in\bH(U): f \mbox{ vanishes at infinity}\}$ and $\bH_b(U):=\{f\in\bH(U): f \mbox{ is bounded}\}$.

\section{Factorization of Dirichlet heat kernels  in metric measure spaces}\label{s:Fac}
\subsection{Setup}\label{s:setup}

We first spell out our assumptions on the state space:
$(\X,\rho)$ is a locally compact separable metric space such that  
all bounded closed sets are compact  and $m$ is a  Radon measure on $\X$ with full support. 
By $B(x,r)=\{y\in\X: \rho(x,y)<r\}$ we denote the open ball of radius $r$ centered at $x\in \X$. For any open set $V$ of $\X$ and $x\in \X$, we denote by $\delta_V (x)$
the distance between $x$ and $\X \setminus V$.

Let $R_0 \in (0, \infty]$ be the largest number
such that $\X \setminus B(x, 2 r) \neq \emptyset$ for all $x \in \X$ and all $r < R_0$. 
We call   $R_0$ the localization radius of $(\X,  \rho)$.

Let $V(x,r):= m(B(x,r)).$
We assume that  there exist constants  $d\ge d_0>0$ such that for every $M\ge 1$ there exists $\wt{C}_M\ge 1$ with the property that
\begin{equation}\label{volume condition}
	\wt C^{-1}_M\left(\frac{R}{r}\right)^{d_0}\le\frac{V(x,R)}{V(x,r)}\le \wt C_M\left(\frac{R}{r}\right)^{d}
\qquad \text{for all}\  x\in \X \text{ and }\ 0<r\le R< M R_0.
\end{equation}
\smallskip
Observe that 
\begin{equation}\label{volume condition3}
	V(x,n_0r) \ge 2V(x,r) \qquad 
\text{ for all}\; x \in \X \ \text{and}\ r \in (0, R_0/n_0),
\end{equation}
where  $n_0:=({2}{\wt C_1})^{1/d_0}$.

Now we spell out the assumptions on the processes we are going to work with.
We assume that $X=(X_t, \P_x)$ is a Hunt process  admitting a 
 (strong) dual Hunt process $\widehat{X}=(\wh{X}_t, \wh{\P}_x)$ with respect to the measure $m$. For the definition of (strong) duality, see \cite[Section VI.1]{BG68}. 
We further assume that the transition semigroups $(P_t)$ and $(\widehat{P}_t)$ of $X$ and $\widehat{X}$ are both 
 Feller and strongly Feller, and that all semipolar sets are polar. 
The condition that semipolar sets are polar is known as Hunt’s hypothesis (H). 
This guarantees  the duality between  the killed processes when the original processes are duals
(since $X$ never hits irregular points).
 See \cite[p.481]{BKK14} and the end of \cite[Section 13.6]{CW05}.

In the sequel, all objects related to the dual process $\widehat{X}$ will be denoted by a hat. We also assume that $X$ and $\widehat{X}$  admit a strictly positive and jointly continuous transition density $p(t,x,y)$ with respect to $m$
so that 
$$
P_tf(x)=\int_{\X}p(t,x,y) f(y)m(dy)\quad \text{and} \quad \widehat{P}_t f(x)=\int_{\X}p(t,y,x)f(y) m(dy).
$$
We will make some assumptions on the transition density $p(t, x, y)$. To do this, we first introduce some notation.

Let $\Phi:(0,\infty) \to (0,\infty)$ be a strictly increasing function satisfying the following weak scaling condition: there exist constants $\delta_l, \delta_u \in (0, \infty)$, $a_l\in(0,1]$, $a_u\in[1,\infty)$ 
 such that
\begin{align}\label{e:scale_Phi}
a_l\left(\frac{R}{r}\right)^{\delta_l} \le \frac{\Phi(R)}{\Phi(r)} \le a_u \left(\frac{R}{r}\right)^{\delta_u}, \quad  r\le R < R_0.
\end{align}
\begin{remark}\label{r:scale-Phi}
Since the function $\Phi$ is strictly increasing,
for every $\wt R\in (0, \infty)$,
there exist $\wt a_l\in(0,1]$ and $ \wt a_u\in[1,\infty)$  such that
\begin{align}\label{e:scale_Phi2}
\wt a_l\left(\frac{R}{r}\right)^{\delta_l} \le \frac{\Phi(R)}{\Phi(r)} \le \wt a_u \left(\frac{R}{r}\right)^{\delta_u}, \quad  0<r\le R \le \wt R.
\end{align}
 Indeed, if $R_0 \le \wt R <\infty$, then for $r \vee  R_0 \le R \le \wt R$, we have $$\frac{\Phi(R)}{\Phi(r)} \le \frac{\Phi(\wt R)}{\Phi(R_0/2)} \frac{\Phi(R_0/2)}{\Phi(r \wedge (R_0/2))} \le a_u\frac{\Phi(\wt R)}{\Phi(R_0/2)}\left(\frac{R_0/2}{r \wedge (R_0/2)}\right)^{\delta_u}\le a_u\frac{\Phi(\wt R)}{\Phi(R_0/2)}\left(\frac{R}{r}\right)^{\delta_u}  $$ 
and
$$
\frac{\Phi(R)}{\Phi(r)} \ge \frac{\Phi(R_0/2)}{\Phi(r \vee (R_0/2))} \ge a_l\frac{\Phi(R_0/2)}{\Phi(\wt R)}\left(\frac{R_0/2}{r}\right)^{\delta_l} \ge a_l\frac{\Phi(R_0/2)}{\Phi(\wt R)}\left(\frac{R_0}{2 \wt R}\right)^{\delta_l}\left(\frac{R}{r}\right)^{\delta_l}.
$$ 
We will use \eqref{e:scale_Phi2} instead of \eqref{e:scale_Phi} whenever necessary.
From \eqref{e:scale_Phi} we can also get the scaling condition for the inverse of $\Phi$:
\begin{align}\label{e:scale_inv_Phi}
a_u^{-1/\delta_u}\left(\frac{R}{r}\right)^{1/\delta_u} \le \frac{\Phi^{-1}(R)}{\Phi^{-1}(r)} \le a_l^{-1/\delta_l}\left(\frac{R}{r}\right)^{1/\delta_l}, 
\quad 0<r \le R < \Phi(R_0).
\end{align}
\end{remark}
 Define 
\begin{align}\label{e:tildeq}
&\wt q(t,x,y) := \frac{1}{V(x,\Phi^{-1}(t))} \wedge \frac{t}{V(x, \rho(x,y)) \Phi(\rho(x,y))}, \quad t>0, \ x,y\in \X.
\end{align}

\begin{remark}\label{r:qcom}
It is easy to see that 
\begin{align*}
\wt q(t,x,y)\asymp \wt q(t,y,x)\asymp 
\frac{1}{V(x,\Phi^{-1}(t))} \wedge \frac{t}{V(y, \rho(x,y)) \Phi(\rho(x,y))}.
\end{align*}
See \cite[Remark 1.12]{CKW16a}.
Moreover, by integrating $\wt{q}(t,x,y)$ over the set $\{y: \rho(x,y) \le \Phi^{-1}(t))\}$, 
one easily gets that for all $t>0$ and $x\in \X$,
\begin{equation}\label{e:integral-q-lb}
\int_{\X}\wt{q}(t,x,y)\, m(dy)\ge 1\, .
\end{equation}
\end{remark}

We will assume that there exists a constant $C_0\ge 1$ such that
\begin{align}\label{e:est_p_D}
C_0^{-1} \wt q(t,x,y) \le p(t,x,y) \le C_0 \wt 
q(t,x,y), \qquad (t,x,y) \in (0,{\wt T})\times \X\times \X
\end{align} 
for some ${\wt T} \in (0, \infty]$.
Then \eqref{e:integral-q-lb} and the lower bound in \eqref{e:est_p_D} yield that
\begin{align}\label{e:int_q}
1 \le \int_{\X} \wt q(t,x,y) \ m(dy) \le C_0
\end{align}
for all $(t,x) \in (0,{\wt T})\times \X$. 
 The processes  
 $X$ (and $\wh X$) may not be  conservative so the lifetimes may be finite. We add  an extra point  $\partial$ (which is called the cemetery point) to $\X$ and assume 
our processes stay at the cemetery point after their lifetimes.
If ${\wt T} =\infty$, we assume $R_0=m(\X)=\infty$. Note that, if 
 ${\wt T} =\infty$ and 
both $X$ and $\wh X$
admit no killing inside $\X$, 
then it follows that $R_0=m(\X)=\infty$, and  $X$ and $\wh X$ are conservative
(see the proof of \cite[Proposition 2.5]{KKW17}, which still works under the non-symmetric setting). 
All functions $h$ on $\X$ will be automatically extended to $\X\cup\{\partial\}$
by setting $h(\partial)=0$.

\begin{remark}\label{r:wtT}
When ${\wt T} \in (0, \infty)$, the value of ${\wt T}$ is not important. That is, 
when ${\wt T} \in (0, \infty)$, for every $T>0$ there exists a constant $\overline C_0=\overline C_0(T)\ge 1$ such that
\begin{align}\label{e:est_p_1}
\overline C_0^{-1} \wt q(t,x,y) \le p(t,x,y) \le \overline C_0 \wt 
q(t,x,y), \qquad (t,x,y) \in (0,{T})\times \X\times \X.
\end{align}
This is a consequence of the  semigroup property of $p(t, x, y)$, 
\eqref{volume condition},  \eqref{e:scale_inv_Phi}  and \eqref{e:est_p_D}.
Indeed, assume $T \ge \wt T$, let $n:=\lfloor 2T/\wt T\rfloor \ge 2$ and fix it.
It follows from \eqref{volume condition},  
 \eqref{e:scale_inv_Phi}  and \eqref{e:est_p_D} that
$$
\wt q(t/n,z,w) \asymp \wt q(t/n^2,z,w) \asymp p(t/n^2,z,w), \quad (t, z,w)\in
(0,T)\times \X \times \X
$$
and
$$
p(t/n,z,w) \asymp \wt q(t/n,z,w)  \asymp \wt q(t,z,w), \quad (t, z,w)\in
(0,T)\times \X \times \X.
$$
Thus by the semigroup property of $p(t,x,y)$,
we have for each $(t, x,y)\in (0,T)\times \X\times \X $,
\begin{align*}
 &p(t, x, y) =\int_{\X} \cdots \int_{\X} p(t/n,x,z_1) \cdots  p(t/n,z_{n-1} ,y)m(dz_1) \cdots m(dz_{n-1}) \\
& \asymp
 \int_{\X} \cdots \int_{\X} \wt q(t/n,x,z_1) \cdots  \wt q(t/n,z_{n-1} ,y)m(dz_1) \cdots m(dz_{n-1})\\
 &\asymp
  \int_{\X} \cdots \int_{\X} p(t/n^2,x,z_1) \cdots  p(t/n^2,z_{n-1} ,y)m(dz_1) \cdots m(dz_{n-1})\\
 &= p(t/n,x,y)   \asymp \wt q(t,x,y).
\end{align*}
\end{remark}

Let $C_0(\X)$ stand for the Banach space of bounded continuous functions on $\X$ vanishing at infinity.
Let $(\sL, \dom(\sL))$ and $(\widehat \sL, \dom(\widehat \sL))$ be the generators of $(P_t)$ and $(\widehat{P}_t)$ in $C_0(\X)$ respectively. 
We assume the following Urysohn-type condition.

\smallskip
\noindent
  {\bf Assumption A}:  
There is a linear subspace $\dom$ of $\dom(\sL) \cap \dom(\widehat{\sL})$ satisfying the following condition:
For any compact $K$ and open $U$ with $K\subset U\subset \X$, 
there is a nonempty collection 
$\dom(K, U)$ of functions $f\in \dom$ satisfying the conditions (i) $f(x)=1$ for $x\in K$; (ii) $f(x)=0$ for $x\in \X\setminus U$; (iii) $0\le f(x)\le 1$ for $x\in \X$, and (iv) the boundary of the set $\{x \, : \, f(x) > 0\}$ has
zero $m$ measure. 
\smallskip

Assumption {\bf A} implies that there exists a kernel $J(x,dy)=J(x, y)m(dy)$ 
(satisfying $J(x, \{x\})=0$ for all $x\in\X$) 
such that $X$
satisfies the following L\'evy system formula (see  \cite[p.482]{BKK14}):
for every stopping time $T$,
\begin{equation}\label{e:levy_system}
\E_x\sum_{s\in (0, T]}f(X_{s-}, X_s)= \E_x\int^T_0\int_{\X}f(X_s, z)J(X_s, dz)ds.
\end{equation}
Here $f:\X\times\X\to[0, \infty]$, $f(x, x)=0$ for all $x\in \X$.
The kernel $J(x,dy)=J(x, y)m(dy)$ is called the jumping kernel of $X$.

Since $J$ satisfies
\begin{equation}\label{e:levykernel}
\int_{\X}f(y)J(x,dy)=\lim_{t\downarrow0}\frac{\E_xf(X_t)}{t}
\end{equation} for all bounded continuous function $f$ on $\X$ and $x\in \X\setminus
\mbox{supp} (f)$,
 we have from \eqref{e:est_p_D} that 
\begin{align}\label{e:c}
C_0^{-1}\frac{1}{V(x, \rho(x,y)) \Phi(\rho(x,y))} \le J(x,y) \le C_0\frac{1}{V(x, \rho(x,y)) \Phi(\rho(x,y))}. 
\end{align} 
Similarly, $\wh X$ has a jumping kernel $\wh J(x,dy)=\wh J(x, y)m(dy)$ with $\wh J(x, y)=J(y, x)$.

There are plenty of examples of processes satisfying the assumptions of this subsection.
Reflected stable-like processes in a closed $d$-set $D\subset \R^d$ satisfy the assumptions
 of this subsection, see \cite{BKK14, CKS15}.
Unimodal L\'evy processes in $\R^d$ with L\'evy exponents
satisfying weak upper and lower scaling conditions at infinity,
in particular, isotropic stable processes, satisfy the assumptions of this subsection, see, for example, \cite{BGR14, CKS-PLMS}.
Another typical example is given at the end of this section.

\subsection{Preliminaries}\label{s:pre}

In this section, we will  argue  that
 $X$ and $\wh X$ satisfy the scale-invariant
parabolic Harnack inequality  with  explicit scaling in terms of $\Phi$. 
Note that $X$ may not be symmetric and may not be conservative.

Let $\tau^X_U:=\inf\{t>0:\, X_t\notin U\}$ be  the first exit time from $U$ for $X$.
If $D$ is an open subset of $\X$, the killed process $X^D$ is defined by $X_t^D=X_t$ 
if $t<\tau^X_D$ and $X_t^D=\partial$ if $t\ge \tau^X_D$,
where 
$\partial$ is the cemetery point added to $\X$.

 Similarly, we define the killed process $\wh X^D$. 
It is well known that $X^D$ and $\wh X^D$ are strong duals of each other with respect to $m_D$, the restriction of $m$ to $D$ 
 (see \cite[p.481]{BKK14} and the end of \cite[Section 13.6]{CW05}). 
For $t>0$, $x,y \in D$,
define 
\begin{align}
\label{e:pD}
 p_D(t,x,y)=p(t,x,y)-\E_x\big[p(t-\tau^X_D,X_{\tau^X_D},y) : \tau^X_D<t<
 \zeta^X\big],
\end{align}
where $ \zeta^X$ is the lifetime of $X$.
By the strong Markov property, $p_D(t,x,y)$ is the transition density of  $X^D$ and, by the continuity of $p(t,x,y)$, \eqref{e:est_p_1}, 
 the Feller and the strong Feller properties of $X$ and $\wh X$, 
it is easy to see that $p_D(t,x,y)$ is jointly continuous (see \cite[pp.34--35]{CZ95} and \cite[Lemma 2.2 and Proposition 2.3]{KSV}).

The following lemma is basically \cite[Lemma 3.8]{BBCK}, 
except that we require neither symmetry nor conservativeness.
\begin{lemma}\label{l:BBCK}
Suppose that there exist positive constants $r, t$ and $p$ such that
\begin{equation}\label{e:BBCK3.22}
\P_x\big(X_s\notin B(x, r), s< \zeta^X\big)\le p, \qquad x\in \X, s\in [0, t].
\end{equation}
Then
$$
\P_x\big(\sup_{0\le s\le t}\rho(X_s, X_0)>2r, t< \zeta^X\big)\le 2p, \qquad x\in \X.
$$
\end{lemma}

\pf Let  $S:=\inf\{s>0: \rho(X_s, X_0)>2r\}$ 
Then using the strong Markov property of 
 $X$ 
and
\eqref{e:BBCK3.22},
\begin{align*}
&\P_x\big(\sup_{0\le s\le t}\rho(X_s, X_0)>2r, t< \zeta^X\big)=\P_x\big(S\le t< \zeta^X\big)\\
&\le \P_x\big(\rho(X_t, X_0)>r, t< \zeta^X\big)+\P_x\big(S\le t< \zeta^X, \rho(X_t, X_0)\le r\big)\\
&\le p+ \P_x\big(S\le t< \zeta^X, \rho(X_t, X_0)\le r\big)\\
&  \le p+\E_x\big[{\bf 1}_{S\le t< \zeta^X}\P_{X_S}(\rho(X_{t-S}, X_0)>r, t-S< \zeta^X)\big]\le 2p.
\end{align*}
\qed

Combining this lemma with \eqref{e:est_p_1} and 
\eqref{e:pD}, we can repeat
the argument of the proof of \cite[Proposition 2.3]{CKK09} word for word to get the following result. Note that conservativeness is not needed.

\begin{proposition}\label{P:3.5-550}
For every 
$a>0$,
there exist  constants
$c>0$ and $\eps \in (0,1/2)$ such that for all $x_0\in \X$ and $r  \in (0, a R_0)$,
\begin{equation}\label{e:aa11}
 p_{B(x_0,r)}(t, x,y)\ge    \, \frac {c}{V(x, \Phi^{-1}(t))}
 \qquad \hbox{for   } x, y\in B(x_0, \eps  \Phi^{-1}(t)) \hbox{ and }
 t\in (0,   \Phi(\eps r)].
\end{equation}
\end{proposition}
\pf Note that by Remarks \ref{r:qcom} and \ref{r:wtT}, there exists $c_1 \ge 1$ such that
\begin{equation}\label{e:push}
c_1^{-1} \wt q(t,y,x) \le p(t,x,y) \le c_1 \wt q(t,x,y), \qquad (t,y,x) \in (0,\Phi(aR_0))\times \X\times \X.
\end{equation}
Let $\eps \in (0,1/2)$ be a small constant which will be chosen later. Observe that for every $x, y\in B(x_0, \eps  \Phi^{-1}(t))$ and $t\in (0,   \Phi(\eps r)]$,
\begin{equation*}
\rho(X_{\tau^X_{B(x_0,r)}},y) \ge \rho(X_{\tau^X_{B(x_0,r)}},x_0) - \rho(x_0, y) \ge r-\eps \Phi^{-1}(t) \ge (\eps^{-1}-\eps)\Phi^{-1}(t).
\end{equation*}
Thus, by \eqref{e:pD}, \eqref{e:push},  \eqref{e:tildeq},  \eqref{volume condition} and \eqref{e:scale_Phi2}, we have that for every $x, y\in B(x_0, \eps  \Phi^{-1}(t))$ and $t\in (0,   \Phi(\eps r)]$,
\begin{align*}
p_{B(x_0,r)}&(t, x, y) = p(t, x, y) - \E_x\big[p(t-\tau^X_{B(x_0,r)},X_{\tau^X_{B(x_0,r)}},y) : \tau^X_{B(x_0,r)}<t<
 \zeta^X\big]\\
& \ge c_1^{-1} \wt q(t, y, x) - c_1 \E_x \left[\frac{t-\tau^X_{B(x_0,r)}}{V(y, \rho(X_{\tau^X_{B(x_0,r)}},y)) \Phi(\rho(X_{\tau^X_{B(x_0,r)}},y))}  : \tau^X_{B(x_0,r)}<t<
 \zeta^X \right] \\
& \ge c_1^{-1} \frac{1}{V(y, \Phi^{-1}(t))} - c_1   \frac{t}{V(y,(\eps^{-1}-\eps)\Phi^{-1}(t)) \Phi((\eps^{-1}-\eps)\Phi^{-1}(t))} \\
& \ge (c_1^{-1}- c_1 \wt C_a \wt a_l^{-1} (\eps^{-1}-\eps)^{-(d_0 + \delta_l)} ) \frac{1}{V(y, \Phi^{-1}(t))}.
\end{align*}
Therefore, by choosing $\eps$ sufficiently small so that $c_1^{-1}- c_1 \wt C_a \wt a_l^{-1} (\eps^{-1}-\eps)^{-(d_0 + \delta_l)} \ge 2^{-1}c_1^{-1}$ and using  \eqref{volume condition}, we conclude the result as
$$
p_{B(x_0,r)}(t, x, y) \ge 2^{-1}c_1^{-1}\frac{1}{V(y, \Phi^{-1}(t))} \ge c_2 \frac{1}{V(x, \Phi^{-1}(t))}.
$$
\qed

Let $\Xi_s:=(V_s, X_s)$ be
the time-space process of $X$, where
$V_s=V_0- s$.
The law of the time-space process $s\mapsto \Xi_s$ starting from
$(t, x)$ will be denoted as $\mathbb{P}^{(t, x)}$.

\begin{definition}
A non-negative Borel   function
$h(t,x)$ on $\bR \times \X$ is said to be {\it parabolic}
(or {\it caloric})
on $(a,b]\times B(x_0,r)$ with respect to $X$
if for every relatively compact open subset $U$ of $(a,b]\times B(x_0,r)$,
$h(t, x)=\mathbb{E}^{(t,x)} [h (\Xi_{\tau^\Xi_{U}})]=\mathbb{E}^{(t,x)}
[h (\Xi_{\tau^\Xi_{U}}):\tau^\Xi_{U} < \zeta^X ]$
for every $(t, x)\in U\cap ([0,\infty)\times
  \X)$,
where
$\tau^\Xi_{U}:=\inf\{s> 0: \, \Xi_s\notin U\}$.

\end{definition}

\begin{theorem}\label{t:PHI} For every $a>0$,
there exist
$c>0$ and $c_1, c_2 \in (0,1)$ depending on $d$,  $\wt T$ and $a$
such that for all $x_0\in
\R^d$, $t_0\ge 0$, $R \in (0, aR_0)$ and every non-negative function $u$ on $[0,
\infty)\times \R^d$ that is
  parabolic
 on $(t_0,t_0+4c_1
\Phi(R)]\times B(x_0,
R)$ with respect to $X$ or $\wh X$,
$$
\sup_{(t_1,y_1)\in Q_-}u(t_1,y_1)\le c \, \inf_{(t_2,y_2)\in
Q_+}u(t_2,y_2),
$$
where $Q_-=(t_0+c_1 \Phi(R),t_0+2c_1  \Phi(R)]\times B(x_0,
c_2R)$
and $Q_+=[t_0+3c_1 \Phi(R),t_0+ 4c_1 \Phi(R)]\times B(x_0,
c_2R)$.
\end{theorem}
\pf
Note that by \eqref{e:c}, 
there exists $c>0$ such that
for all  $x \not=y\in \X$ with $r\le  \rho(x,y)/2 <R_0$,
$$
 J(x, y) \le \frac{c}{V(x,r)}
\int_{B(x,r)} J(z,y)m(dz) \text{ and }
J(x,y) \le \frac{c}{V(y,r)}
\int_{B(x,r)} J(y, z)m(dz). 
 $$
Thus, using this, Proposition \ref{P:3.5-550} and \eqref{e:est_p_D}, we see that   
 the proof of Theorem \ref{t:PHI}  is almost identical to the proof for the symmetric case in \cite[Theorem 4.3]{CKW16b}.
 We emphasize that the conservativeness is not used in the proofs of \cite[Lemmas 3.7, 4.1, 4.2 and Theorem 4.3]{CKW16b}.
We omit the details.
\qed.

Theorem  \ref{t:PHI} clearly implies the elliptic Harnack inequality. Using Theorem  \ref{t:PHI}, 
we have the following result. In the remainder of this section, $D$ will always stand for an open
subset of $\X$.

\begin{proposition}\label{step1}
For all  $a, b>0$,  
  there exists
 $c=c(a, b)>0$ such that for every open set $D\subset \X$, 
$$
p_D(t,x,y) \,\ge\,\frac{c}{V(x, \Phi^{-1}(t))}
$$
for all $(t, x, y)\in (0,  a R_0)\times D\times D$ with
$\delta_D (x) \wedge \delta_D (y)
 \ge b \Phi^{-1}(t) \geq 4 \rho(x,y)$.
 \end{proposition}
\pf
See the proof of \cite[Proposition 3.4]{KM18}.
\qed

The proof of the next result is also standard.
\begin{proposition}\label{step3}
For all $a, b>0$,
there exists  $c=c(a,b)>0$ such that for every open set $D\subset\X$, 
$p_D(t, x, y)\,\ge \,c\, t\, 
 J(x,y)$
 for all $(t, x, y)\in (0,a R_0) \times D\times D$ with $\delta_D(x)\wedge \delta_D (y) \ge b \Phi^{-1}(t)$ and
$b \Phi^{-1}(t) \leq 4 \rho(x,y)$.
\end{proposition}

\pf
See the proofs of \cite[Lemma 3.4 and Proposition 3.5]{CKS-PLMS}.
\qed

Combining the two propositions above we get

\begin{proposition}\label{p:26}
For all  $a, b>0$, there exists $c=c(a, b)>0$ 
such that for every open set $D\subset\X$,
$
p_D(t,x,y)\ge c\wt q(t,x,y)
$
for all $t\in (0,a R_0), x,y\in D$ with $\delta_D(x) \wedge \delta_D(y)\ge b\Phi^{-1}(t)$.
\end{proposition}

\subsection{3P inequality and Feynman-Kac perturbations}\label{s:3p}

We first prove  the following 3P inequality.
\begin{lemma}\label{l:3P}
For every $a \in (0, \infty)$, there exists $c>0$ such that for all $0<s<t < a R_0$,
\begin{align}\label{e:3P}
\frac{\wt q(s,x,z)\wt q(t-s,z,y)}{\wt q(t,x,y)} \le c  (\wt q(s,x,z)+\wt q(t-s,z,y)), \quad x,y,z \in \X.
\end{align}
\end{lemma}
\pf
Note that, by the triangle inequality, either $\rho(x,z) \ge 2^{-1} \rho(x,y)$ or $\rho(z,y)\ge 2^{-1} \rho(x,y)$. 
Since the argument is the same, we 
only give the proof for the case $\rho(z,y)\ge 2^{-1} \rho(x,y)$.
Thus we assume that
$\rho(z,y)\ge 2^{-1} \rho(x,y)$ and set
\begin{equation}\label{e:4.6}
\bar q( t,x, y, r):= {t V(x, \Phi^{-1}(t)) + \Phi (r) \, V(y, r)}, 
\quad t, r >0 \hbox{ and } x, y\in \X.
\end{equation}
Note that by Remark \ref{r:qcom} and the fact $1 \wedge (1/r)  \asymp 1/(1+r) $ for $r>0$, we have that for all
  $t>0$  and  $x, y \in \X$,
\begin{equation}\label{e:4.7}
\wt q(t, x, y) \asymp \frac{t}{\bar q(t, x, x, \rho(x, y))} 
\end{equation}
and
\begin{equation}\label{e:4.71}
{\bar q(t, x, x, \rho(x, y))} \asymp {\bar q(t, y, y, \rho(x, y))} \asymp{\bar q(t, x,y, \rho(x, y))} 
\asymp{\bar q(t, y,x, \rho( x,y))}.
\end{equation}

We claim that 
\begin{align}
\label{e:cl2.1}
{\bar q(t, x, x, \rho(x, y))} \le c (\bar q( t-s,y, y,  \rho(z,y))+\bar q( s,x, x, \rho(x,z))  ).
\end{align}
The assertion of the lemma follows easily from \eqref{e:cl2.1}.
Indeed,  \eqref{e:cl2.1} implies that 
\begin{align*}
&\frac{\wt q(s,x,z)\wt q(t-s,z,y)}{\wt q(t,x,y)}
\asymp\frac{s(t-s)}{t}
 \frac{{\bar q(t, x, x, \rho(x, y))}}{
\bar q( t-s,y, y,  \rho(z,y))\bar q( s,x, x, \rho(x,z)) }
 \\  \le &c (s \wedge (t-s)) \left(\frac{1}{\bar q( s,x, x, \rho(x,z))}+ \frac1{\bar q( t-s,y, y,  \rho(z,y))}
 \right)\\
 \le &c\left( \frac{s}{\bar q(s,x,x,\rho(x,z))}+\frac{t-s}{\bar q(t-s,y,y,\rho(z,y))}\right)
  \asymp \wt q(s,x, x, \rho(x,z))+\wt q(t-s, y, y,\rho(z,y)). 
\end{align*}

We now prove the claim  \eqref{e:cl2.1} by considering two cases separately.

\noindent
(1) {  $\rho(x,y) \le \Phi^{-1}(t)$:}
In this case, ${\bar q(t, x, x, \rho(x, y))} \asymp t V(y, \Phi^{-1}(t)) \asymp  t V(x, \Phi^{-1}(t))$.
By \eqref{volume condition} and \eqref{e:scale_inv_Phi},
if $t/2 \ge s$, 
$$ t V(y, \Phi^{-1}(t)) \le c (t-s)V(y, \Phi^{-1}(t-s)) \le c \left[\bar q( t-s,y, y,  \rho(z,y))+\bar q( s,x, x, \rho(x,z)  ) \right].
$$
Similarly, if $t/2 < s$,
$$
t V(x, \Phi^{-1}(t)) \le c sV(x, \Phi^{-1}(s)) \le c (\bar q( t-s,y, y,  \rho(z,y))+\bar q( s,x, x, \rho(x,z))  ).
$$
Thus, by \eqref{e:4.71}, \eqref{e:cl2.1} holds true.

\noindent
(2) {  $\rho(x,y) > \Phi^{-1}(t)$:}
The assumption $\rho(z,y)\ge 2^{-1} \rho(x,y)$ and  \eqref{volume condition} imply that
\begin{align}
\label{e:comp21}
\frac{V(x,2\rho(y,z))}
{V(y,\rho(y,z))} \le \frac{V(y,2\rho(x,y)+2\rho(y,z))}
{V(y,\rho(y,z))}\le \frac{V(y,6\rho(y,z))}
{V(y,\rho(y,z))} \le c.
\end{align}
Combining \eqref{e:comp21} and \eqref{volume condition} we get that for
$x,y,z \in D$ and $0<s<t < a R_0$,
\begin{align*}
{\bar q(t, x, x, \rho(x, y))}&\asymp \Phi (\rho(x,y)) \, V(x,\rho(x,y))\\
&\le
 \Phi (\rho(x,z)+\rho(z,y)) \, V(x,\rho(x,z)+\rho(z,y))
 \\ & \le 
 \Phi \big(2 (\rho(x,z)\vee \rho(z,y))\big) \, V\big(x,2 (\rho(x,z)\vee \rho(z,y))\big)\\
&\le 
 \Phi (2 \rho(z,y)) \, V(x,2 \rho(z,y))+
 \Phi (2\rho(x,z) ) \, V(x,2 \rho(x,z))\\
&\le  c\left[
\Phi (\rho(z,y)) \, V(y,\rho(z,y))+
 \Phi (\rho(x,z) ) \, V(x,\rho(x,z))\right]
\\
 &\le c \left[\bar q( t-s,y, y,  \rho(z,y))+\bar q( s,x, x, \rho(x,z))\right].
\end{align*}
We have proved \eqref{e:cl2.1}.
\qed

Recall that, for an open set $D\subset\X$,  a measure $\mu$ on $D$ is said to be a smooth measure
of $X^D$ with respect to the reference measure $m_D$
if there is a positive continuous additive functional
$A$ of $X^D$ such that for any
bounded non-negative Borel function $f$ on $D$,
\begin{equation}\label{eqn:Revuz1}
\int_D f(x) \mu (dx) = \lim_{t\downarrow 0}
\E_{m} \left[ \frac1t \int_0^t f(X^D_s) dA_s \right],
\end{equation}
cf. \cite{Revuz}.
The additive functional $A$ is called the positive continuous
additive functional of $X^D$ with Revuz measure $\mu$
with respect to the reference measure $m_D$.

It is known (see \cite{FG96}) that for any $x\in D$, $\alpha\ge 0$ and
bounded non-negative Borel function $f$ on $D$,
$$
\E_x \int^{\infty}_0e^{-\alpha t}f(X^D_t) dA_t =
\int^{\infty}_0e^{-\alpha t}\int_D
p_D(t,x,y)
f(y) \mu (dy)dt,
$$
and we have for any $x\in D$,  $t>0$ and
non-negative Borel function $f$ on $D$,
\begin{equation}\label{eqn:revuz2}
\E_x \int^t_0f(X^D_s) dA_s=\int^t_0\int_D
p_D(s, x, y)
f(y) \mu (dy)ds.
\end{equation}

We first introduce our class of possibly critical perturbations.
For an open set $D\subset\X$, 
a smooth Radon measure $\mu$ of $X^D$, $t>0$ and  $a \ge 0$, we define
\begin{align}\label{d:N_local}
    N^{D, \mu}_a(t) := \sup_{x\in \X} 
 \int_0^t \int_{z\in D: \delta_D(z) >a \Phi^{-1}(t)}  \wt q(s,x,z)\mu(dz)ds.
\end{align}
\begin{definition}\label{d:KT}
Let $\mu$ be a smooth measure for  both $X^D$ and $\wh X^D$ with respect to the reference measure $m_D$ and let 
$T\in (0, \infty]$. 
The measure $\mu$ is said to be in the  class 
$\mathbf{K}_T(D)$ if

\noindent
(1) $\displaystyle \sup_{t <T}N^{D, \mu}_a(t)<\infty$ for all $a \in (0, 1]$;

\noindent
(2)
 $\displaystyle\lim_{t \to 0}N^{U, \mu}_0(t)=0$ for every relatively compact open set $U$ of $D$.
\end{definition}

 For $\mu\in \mathbf{K}_T(D)$, using condition (2) in the definition above, 
one can show that, for any relatively compact open subset $U$ of $D$, $A_{t \wedge \tau^X_U}$ is a positive continuous
additive functional of $X^U$ with Revuz measure $\mu_U$,
where $\mu_U$  is the measure $\mu$ restricted to $U$. 
See Appendix for the proof.

\begin{remark}
Note that by the semigroup property, it is easy to check that 
$$
N^{D, \mu}_a(t) \le N^{D, \mu}_a(s) +  \overline  C_0(T)^2  
N^{D, \mu}_a(t-s), \quad 0<s<t \le T,$$
where $\overline C_0(T)$ is the constant in \eqref{e:est_p_1}.
Thus, if $\mu$ is in the  class 
$\mathbf{K}_1(D)$,
then $\sup_{t < T}N^{D, \mu}_a(t)<\infty$ for all $a>0$ and $T \in (0, \infty)$. 
\end{remark}

For $\mu \in \mathbf{K}_1(D)$, we denote by
$A^{\mu}_t $  the positive continuous additive functional  of $X^D$ 
with Revuz measure $\mu$ and denote by 
$\wh A^{\mu}_t $  the 
positive continuous additive functional of
 $\wh X^D$ 
with Revuz measure $\mu$. 
For any non-negative Borel function $f$ on $D$, we define
\begin{align*}
T^{\mu,D}_t f(x) = \E_x\left[ \exp (-A^{\mu}_t) f(X^D_t)\right], \quad \wh T^{\mu,D}_t f(x) = \wh \E_x\left[ \exp (- \wh A^{\mu}_t) f(\wh X^D_t)\right],  
\quad t\ge 0 , x\in D.
\end{align*}
The semigroup $(T^{\mu,D}_t: t\ge 0)$ (respectively $(\wh T^{\mu,D}_t: t\ge 0)$) 
is called the Feynman-Kac semigroup of $X^D$ (respectively $\wh X^D$) 
associated with $\mu$. 
By \cite[Theorem 6.10(2)]{Y96}, 
$T^{\mu, D}_t$ and $\wh T^{\mu, D}_t$ are duals of each other
with respect to the measure $m_D$ so that
\begin{align}
\label{e:Tdual}
\int_D T^{\mu,D}_t f(x) g(x)m(dx)=\int_D  f(x) \wh T^{\mu,D}_tg(x)m(dx).
\end{align}

Let $Y$ ($\wh Y$, respectively) be a Hunt process on $D$ corresponding to the  transition semigroup $(T_t^{\mu,D})$ ($(\wh T_t^{\mu,D})$, respectively). 
For an open subset $U \subset D$, we denote by $Y^{U}$  ($\wh Y^{U}$, respectively)
the process $Y$  ($\wh Y$, respectively) killed upon exiting  $U$.

Suppose that $U \subset D$ is a relatively compact open subset of $D$.
Since for any relatively compact open set $U$, $A^{\mu, U}_t:=A^{\mu}_{t \wedge \tau^X_U}$
is a positive continuous additive functional of $X^U$ with Revuz measure $\mu_U$,
the transition semigroup of $Y^{U}$ is  $(T_t^{\mu_U,U})$.
For simplicity, in the sequel we denote this semigroup as $(T_t^{\mu, U})$.
Moreover, for any $t\ge 0 , x\in U$, 
\begin{align*}
T^{\mu,U}_t f(x) =\E_x\left[  f(Y^U_t)\right]= \E_x\left[ \exp \big(-A^{\mu}_{t\wedge 
 \tau^X_U}\big) f(X^U_t)\right]
\end{align*}
and
\begin{align*}
\E_x \int^t_0f(X^U_s) dA^{\mu, U}_s=\int^t_0\int_U
p_U(s, x, y)
f(y) \mu (dy)ds.
\end{align*}

It follows from Definition \ref{d:KT}(2) that, for all relatively compact open subset $U$ of $D$,
$$
\lim_{t\to 0}\sup_{x\in U}\int^t_0\int_Up_U(s, x, y)\mu(dy)ds=0,
$$
i.e., $\mu_U$ is in  the standard Kato class of $X^U$.
Thus, according to the discussion in \cite[Section 1.2]{CKS15}, 
we have for any non-negative bounded  Borel function $f$ on $U$,
\begin{align*}
T^{\mu,U}_t f(x) = \E_x\left[  f(X^U_t)\right]+ \E_x\left[ f(X^U_t)\sum_{n=1}^\infty (-1)^n
\int^t_0\int^{s_n}_0\cdots \int^{s_2}_0 
dA^{\mu, U}_{s_1}
 \cdots dA^{\mu, U}_{s_{n-1}}
dA^{\mu, U}_{s_n}\right].
\end{align*}
Define $p_U^0(t,x,y) := p_U(t, x, y)
$ and, for $k \ge 1$, 
\begin{align}
p_U^{k}(t,x,y) =&- \int_0^t \int_U p_U(s, x, z)
p_U^{k-1}(t-s,z,y) \mu(dz) ds
.\label{d:pk}
\end{align}
Repeating the discussion in \cite[Section 1.2]{CKS15}, one can conclude that
\begin{align*}
T^{\mu,U}_t f(x) = \int_U q^U(t,x,y) f(y) m(dy), \quad (t,x) \in (0,\infty) \times U,
\end{align*}
where 
$
q^U(t,x,y) := \sum_{k=0}^\infty p_U^{k}(t,x,y).$

By Lemma \ref{l:3P},
we have that for any $\mu$ in $\mathbf{K}_1(D)$,
any relatively compact open set $U$ of $D$ and any $(t,x,y)\in(0,1]\times U\times U$,
\begin{align}
&\int_0^t \int_{U}  \wt q(t-s,x,z) \wt q(s,z,y) \mu(dz) ds \nn\\
&\le c \wt q(t,x,y) \sup_{u\in \X} \int_0^t \int_{U}   \wt q(s,u,z)\mu(dz)ds =c \wt q(t,x,y) N^{U, \mu}_0 (t).\label{e:int_qqmu}
\end{align}
Using \eqref{e:int_qqmu} and the semigroup property, it is standard to show that 
$p_U^{k}(t,x,y)$ is continuous in $(t,y)$ for each fixed $x$, continuous in $(t,x)$ for each fixed $y$, and
$\sum_{k=0}^\infty p_U^{k}(t,x,y)$ converges absolutely and uniformly so that   $q^U(t,x,y) $ is continuous in $(t,y)$ for each fixed $x$, 
and also continuous in $(t,x)$ for each fixed $y$  (for example, see \cite{CKS15}).

Define 
$
q^D(t,x,y) :=\lim_{n \to \infty} q^{D_n} (t,x,y),
$
where 
$D_n \subset D$ are bounded increasing open sets such that $ \overline{D_n} \subset D_{n+1}$ and $\cup_{n=1}^\infty D_n=D$. 
Then, using the monotone convergence theorem and 
$$q^{D_n}(t,x,y)\le p_{D_n}(t,x,y)\le p(t,x,y)\le \overline C_0(T) \wt q(t,x,y), \quad t <T,$$
we see that 
$q^D(t,x,y)$ is the transition density of $Y$ and 
$q^{D} (t,x,y) \le 
\overline C_0(T)\wt q(t,x,y)$ for $t <T$.
Therefore, we obtain the following 
\begin{proposition}\label{t:exsitq}
Suppose that $D$ is an open set in $\X$ and $\mu \in \mathbf{K}_1(D)$. Then
the Hunt process $Y$ on $D$ corresponding to the  transition semigroup 
$(T_t^{\mu,D})$ has a transition density $q^D(t,x,y)$ 
with respect to $m$ 
such that for each $T \in (0, \infty)$, 
$q^{D} (t,x,y) \le 
\overline C_0(T)\wt q(t,x,y)$ for $t <T$. 
Furthermore, if $D$ is relatively compact, then
$q^D(t,x,y)$ is continuous in $(t,y)$ for each fixed $x$, 
and continuous in $(t,x)$ for each fixed $y$.
If $\mu\in \mathbf{K}_{\infty}(D)$ 
and $\wt T=\infty$, then the estimate $q^{D} (t,x,y) \le 
 c\wt q(t,x,y)$ holds for every $t>0$.
\end{proposition}

\subsection{Interior estimates}\label{s:ie}
In this subsection, we prove some interior estimates for the transition density $q^U(t,x,y)$, where $U$ is an open subset of  $D$.
Recall that we assume $R_0=m(\X)=\infty$ when  ${\wt T} =\infty$.

\begin{theorem}\label{t:sharpest_qD}
Suppose that $\mu \in \mathbf{K}_1(D)$. 
Then for every $T \in (0, \infty)$ and $a\in (0, 1]$,
 there exists a constant  
 $ c:=c(a, \Phi, C_0,M,\sup_{t \le T}N^{D, \mu}_{ 2^{-1}a }(t)) >0$ 
 such that for every open $U \subset D$,
\begin{align}\label{e:lower_bound}
q^U(t,x,y)  \ge c \wt q(t,x,y) 
\end{align}
for all $t\in (0,T), x,y\in U$ with $\delta_U(x) \wedge \delta_U(y)\ge a\Phi^{-1}(t)$. Moreover, if 
$\mu \in \mathbf{K}_{\infty}(D)$ and
 ${\wt T} =\infty$,
 then \eqref{e:lower_bound} holds for all $t>0$.
\end{theorem}

\pf 
Fix $t\in (0,T)$, $x,y\in U$ with $\delta_U(x) \wedge \delta_U(y)\ge a\Phi^{-1}(t)$ and choose a bounded open set  $V \subset
\{z \in U:\delta_U(z) > 2^{-1}a\Phi^{-1}(t)\}$ containing $x,y$ such that $\delta_V(x) \wedge \delta_V(y)\ge 2^{-2}a\Phi^{-1}(t)$   
(for example, one can take $V =
\{z \in U:\delta_U(z) > 2^{-1}a\Phi^{-1}(t)\} \cap B(x, 2(\rho(x,y)+\Phi^{-1}(t))$). Note that  
$q^U(t,x,w) \ge q^V(t,x,w)$ for all $w \in V$ and 
$ w \mapsto q^V(t,x,w)$  is continuous.

For $w \in V$, let 
\begin{align*}
\wt{p}^1_V(t,x,w) &: = \int_0^t\left(\int_V 
p_V(t-s, x, z) p_V(s, z, w)
\mu(dz) \right) ds .
\end{align*}
Then for any bounded Borel function $f$ on $V$,  by the Markov property of $X^V$,  we have
\begin{align}\label{e:AV}
\E_x \left[ A^{\mu_V}_t f(X^V_t)\right] =  \E_x \left[ \int_0^t \E_{X_s^V}[f(X_{t-s}^V)] dA^{\mu_V}_s \right] = 
\int_V \wt{p}^1_V(t,x,w) f(w) m(dw).
\end{align}

Observe that, since 
$\delta_D(z) \ge \delta_U(z)   > 2^{-1}a\Phi^{-1}(t) 
$ for $z \in V$,
by  \eqref{e:pD}, \eqref{e:est_p_D},  Lemma \ref{l:3P} and Proposition \ref{p:26}, we have that for $w \in B(y, 2^{-3} a\Phi^{-1}(t))$,
\begin{align*}
&\wt{p}^1_V(t, x, w)\, \le\, 
\overline C_0^2 \int_0^t \int_V \wt q(t-s,x,z)\wt q(s,z,w) \mu(dz) ds\\
 &\le \overline C_0^2 \int_0^t \int_{z \in D: \delta_D(z) > 2^{-1}  a \Phi^{-1}(t)}  \wt q(t-s,x,z)\wt q(s,z,w) \mu(dz) ds\\
& \le \overline C_0^2c\left(\sup_{s \le T}N^{D, \mu}_{ 2^{-1}a }(s)  \right)\wt q(t, x, w)
 \le \overline C_0^2c\left(\sup_{s \le T}N^{D, \mu}_{ 2^{-1}a }(s)  \right) C_*^{-1}
p_V(t, x, w)\\
&=:\,(k/2) p_V(t, x, w).
\end{align*}
 Hence, for  $w \in B(y, 2^{-3} a\Phi^{-1}(t))$, we have $p_V(t, x, w)
- k^{-1} \wt{p}^1_V(t, x, w)\ge 2^{-1}p_V(t, x, w)$, which implies that for any $r<2^{-3}a\Phi^{-1}(t)$,
\begin{equation}\label{e:lbd1}
\frac{1}{2} \E_x[{\bf 1}_{B(y,r)}(X_t^V)] \le \E_x\left[ \big(1- A_t^{\mu_V}/k\big){\bf 1}_{B(y,r)}(X_t^V) \right].
\end{equation}
Using the elementary fact that
$
1-A^{\mu_V}_t/k\le \exp\left(-A^{\mu_V}_t/k\right),
$
we get that for any $r < 2^{-3} a\Phi^{-1}(t)$,
$$
\frac1{V(y, r)}\E_x\left[\big(1-A^{\mu_V}_t/k \big)
{\bf 1}_{B( y ,r)}(X_t^V) \right] \le \frac1{V(y, r)}
\E_x\left[\exp (-A^{\mu_V}_t/k)
{\bf 1}_{B(y, r)}(X_t^V) \right].
$$
Thus, by \eqref{e:lbd1}, \eqref{e:AV} and H\"{o}lder's inequality, we have
\begin{eqnarray*}
&&\frac12\frac1{V(y, r)}\E_x\left[
{\bf 1}_{B(y, r)}(X_t^V) \right]\,\le \, \frac1{V(y, r)}\E_x\left[\exp (-A^{\mu_V}_t/k)
{\bf 1}_{B(  y, 
r)}(X_t^V) \right]\\
&&\le
\left(\frac1{V(y, r)}\E_x\left[\exp (-A^{\mu_V}_t)
{\bf 1}_{B(y, r)}(X_t^V) \right] \right)^{1/k}
\left(\frac1{V(y, r)} \E_x\left[
{\bf 1}_{B(y, r)}(X_t^V) \right]\right)^{1-1/k}.
\end{eqnarray*}
Therefore,
$$
\frac1{2^k}\frac1{V(y, r)} \E_x\left[
{\bf 1}_{B(y, r)}(X_t^V) \right]
\le \frac1{V(y, r)}
 \E_x\left[\exp (-A^{\mu_V}_t )
{\bf 1}_{B(y, r)}(X_t^V)
\right].
$$
Since $w \to q^V(t, x, w)$ is continuous by Proposition \ref{t:exsitq}, 
we conclude by sending $r\downarrow 0$ and applying Proposition \ref{p:26} again that for every
$t\in (0,T], x,y\in U$ with $\delta_U(x) \wedge \delta_U(y)\ge a\Phi^{-1}(t)$,
$$
 q^U(t, x, y) \geq  q^V(t, x, y) \geq 2^{-k}
p_V(t, x, y)
\ge c 2^{-k}\wt q(t, x, y) .
$$ 
\qed

Let $\tau_{U}:=\inf\{s> 0: \, Y_s\notin U\}$ and $\wh \tau_{U}:=\inf\{s> 0: \, \wh Y_s\notin U\}$.
Since the proofs for the dual processes are  same, throughout the paper we give the proofs for $Y$ only.

\begin{corollary}\label{c:18}
(1)
Suppose that $\mu \in \mathbf{K}_1(D)$.
For any positive constants $R, T$ and $a$, there exists $c_1=c_1(a,T)>0$ such that for all $t \in (0, T)$ and $B(x,\Phi^{-1}(t)) \subset D$,
\begin{equation}\label{e:c181}
\inf_{z \in B(x,a\Phi^{-1}(t)/2)}\P_z(\tau_{B(x,a\Phi^{-1}(t))}> t) \wedge \inf_{z \in B(x,a\Phi^{-1}(t)/2)}\wh \P_z(\wh \tau_{B(x,a\Phi^{-1}(t))}> t) \ge c_1
\end{equation}
and
\begin{equation}\label{e:c1820}
 \E_x[\tau_{B(x,r)}] \wedge \wh \E_x[\wh \tau_{B(x,r)}] \ge c_1 \Phi(r).
\end{equation}
Moreover,  there exist $r_1, c_2 >0$ such that  for all $r \in (0,r_1]$ and $B(x,r) \subset D$, 
\begin{equation}\label{e:c182}
 \E_x[\tau_{B(x,r)}] \vee \wh \E_x[\wh \tau_{B(x,r)}] \le c_2 \Phi(r).
\end{equation}
(2) 
If $\mu \in  \mathbf{K}_{\infty}(D)$ and
 ${\wt T} =\infty$ (and $R_0=\infty$), 
then \eqref{e:c181}--\eqref{e:c182} hold for all $r,t>0$.
\end{corollary}
\pf (1)
For any $z \in B(x,a\Phi^{-1}(t)/2)$, 
we have by Theorem \ref{t:sharpest_qD} that
\begin{align*}
&\P_z(\tau_{B(x,a\Phi^{-1}(t))}> t) \ge \P_z(\tau_{B(z,a\Phi^{-1}(t)/6)}> t) \\
& \qquad  =\int_{B(z,a\Phi^{-1}(t)/6)} q^{B(z,a\Phi^{-1}(t)/6)}(t,z,y)m(dy) \\
& \qquad \ge c \int_{B(z,a\Phi^{-1}(t)/12)} 
 \frac1{V(z, \Phi^{-1}(t))}
\wedge \frac{t}{V(z,\rho(z,y))\Phi(\rho(z,y))} m(dy)  \ge c_0.
\end{align*}

 \eqref{e:c1820} is clear from \eqref{e:c181}. In fact, $$\E_x[\tau_{B(x,r)}] \ge \Phi(r) \P_x(\tau_{B(x,r)}>\Phi(r)) \ge c_1 \Phi(r).$$
 
We now prove   \eqref{e:c182}.
By the semigroup property and
Proposition \ref{t:exsitq},
  we have that for 
   $t>2^ks$ and $s\le1$, 
\begin{align*}
&q^{B(x,r)} (t,x,y)
=\int_{B(x,r)} q^{B(x,r)} (t-2^ks,x,z)q^{B(x,r)} (2^ks,z,y) m(dz)\\
& \le \sup_{z,y \in B(x,r)} q^{B(x,r)} (2^ks,z,y) 
  = \sup_{z,y \in B(x,r)}  \int_{B(x,r)} q^{B(x,r)} (2^{k-1}s,z,w) q^{B(x,r)} (2^{k-1}s,w,y) dw\\
&  \le \left(\sup_{z,y \in B(x,r)} q^{B(x,r)} (2^{k-1}s,z,y) \right)^2 V(x, r)\\ 
& \le \cdots \le 
 \left(\sup_{z,y \in B(x,r)} q^{B(x,r)} (s,z,y) \right)^{2^k} V(x, r)^{2^{k-1}+\cdots+2+1} \\
 &\le
 \left(\sup_{z \in B(x,r)}\frac{c} {V(z, \Phi^{-1}(s))}\right)^{2^k} V(x, r)^{2^k-1}.
\end{align*}
Using this, 
we have that for all $A>0$ and $r\le\Phi^{-1}(1/A)$,
\begin{align*}
&\E_{x}[\tau_{B({x},r)}] 
      \le
 A \Phi(r)+ \sum_{k=0}^\infty \int_{ A 2^{k} \Phi(r)}^{ A 2^{k+1}  \Phi(r)} \int_{B({x},r)}q^{B({x},r)} (t,{x},y) m(dy) dt\\
  & \le
  A \Phi(r)\left(1+ \sum_{k=0}^\infty 2^k \left(\sup_{z,y \in B({x},r), t \ge A 2^{k} \Phi(r)} q^{B({x},r)} (t,z,y) \right)V({x}, r) \right)\\
& \le A \Phi(r)\left(1+ \sum_{k=0}^\infty 2^k  \left(\sup_{z \in B(x,r)}\frac{cV(x, r)} {V(z, \Phi^{-1}(A\Phi(r)))}\right)^{2^k} \right).
\end{align*}
Using \eqref{volume condition} 
 and \eqref{e:scale_inv_Phi}, 
for all $z \in B(x,r)$,  $A>1$  and $r\le \Phi^{-1}(1/A)$, 
\begin{align*}
&\frac{cV(x, r)} {V(z, \Phi^{-1}(A\Phi(r)))}
=c\frac{V(x, \Phi^{-1}(A\Phi(r)))} {V(z, \Phi^{-1}(A\Phi(r)))}\frac{V(x, r)} {V(x, \Phi^{-1}(A\Phi(r)))}\\
&\le c\frac{V(z, r+\Phi^{-1}(A\Phi(r)))} {V(z, \Phi^{-1}(A\Phi(r)))}c_2\left(\frac{r} { \Phi^{-1}(A\Phi(r))}\right)^{d_0}\\
&\le c\frac{V(z, 2\Phi^{-1}(A\Phi(r)))} {V(z, \Phi^{-1}(A\Phi(r)))}c_2\left(\frac{\Phi^{-1}(\Phi(r))} { \Phi^{-1}(A\Phi(r))}\right)^{d_0}
\le c_3 2^d c_2c_4 A^{-d_0/\delta_u}.
\end{align*}
Choose $A >1$ large so that $\eps:=c_3 2^d c_2c_4 A^{-d_0/\delta_u}<1$. Then for $r\le r_1:=\Phi^{-1}(1/A)$
\begin{align*}
&\E_x[\tau_{B(x,r)}]  \le A \Phi(r)\left(1+ \sum_{k=0}^\infty 2^k  \eps^{2^k} \right) \le c  \Phi(r).
\end{align*}

(2) The proof is similar to the proof of (1). We omit the details. 
\qed

\subsection{Examples of critical potentials}\label{s:ex}

In this subsection, we give two examples of critical potentials.

\begin{example}\label{e:1}
Suppose $\mu(dz) =  q(z)m(dz) $ 
with $ 0\le q(z) \asymp 1/\Phi( \delta_D(z) \wedge 1)$. 
Since $q$ is bounded on every relatively compact open set $U\subset D$, 
$N^{U, \mu}_0(t) \le C_0 t\|q\|_{L^\infty(U)} \to 0$ as $t \to 0$.
Moreover, 
for $x\in D$, $a \in (0, 1]$ and $t  < 1$,
\begin{align}\label{e:ex1}
 &\int_0^t \int_{z\in D: \delta_D(z) >a\Phi^{-1}(t)}  \wt q(s,x,z)q(z)m(dz)ds\nn\\
  \le & ct
  +c \int_0^t \int_{z\in D: 1>\delta_D(z) >a\Phi^{-1}(t)} \Phi( \delta_D(z))^{-1} \wt q(s,x,z)m(dz)ds\nn\\
 \le &ct+c
 \frac1{\Phi(a\Phi^{-1}(t) )}\int_0^t \int_{D}  \wt q(s,x,z)m(dz)ds
 \le ct+
c \frac{t}{\Phi(a\Phi^{-1}(t) )} <c<\infty.
\end{align}
Thus $\sup_{t < 1}N^{D, \mu}_a(t)<\infty$ for all $a \in (0, 1]$ and so $\mu$ is in the class 
$ \mathbf{K}_1(D)$.
\end{example}

\begin{example}\label{e:2}
Suppose $\wt T=\infty$  and $\mu(dz) =  q(z)m(dz) $ 
with $ 0\le q(z) \asymp 1/\Phi( \delta_D(z))$. 
Then $R_0=\infty$ and  for all $a \in (0, 1]$ and  $t < \infty$,
\begin{align}\label{e:ex2}
 N^{D, \mu}_a(t) \le c
 \frac1{\Phi(a\Phi^{-1}(t) )}\sup_{x \in \X}\int_0^t \int_{D}  \wt q(s,x,z)m(dz)ds\le 
 \frac{ct}{\Phi(a\Phi^{-1}(t) )} <c<\infty.
\end{align}
Thus $\mu$ is in the class 
$\mathbf{K}_{\infty}(D)$.
\end{example}

\subsection{Factorization of Dirichlet heat kernel in $\kappa$-fat open set}\label{s:factorization}

Recall that $\dom(K, U)$
 is the subset of $\dom$ in Assumption {\bf A}.
Let 
$$
A(z_0,p,q) := \{ x\in \X : p<\rho(x,z_0)<q \}\quad \text{ and } \quad \overline A(z_0,p,q) := \{ x\in \X : p\le \rho(x,z_0) \le q \}.
$$
Note that, due to our assumption that  all bounded closed sets are compact, 
$\overline A(z_0,p,q)$ is compact. Thus by
Assumption \textbf{A}, 
for any  $1/2<b<a<1$, the set  $\mathcal{D}(\overline{A}(z_0, br, ar), A(z_0, r/2, r))$ is non-empty.
We now add the final assumption saying that 
there exist proper bump functions in each non-empty set $\mathcal{D}(\overline{A}(z_0, br, ar), A(z_0, r/2, r))$
 providing scale-invariant control on the action of the generator.

\medskip
\noindent    {\bf Assumption U}: 
There exists $r_0 \in (0, \infty]$ such that for any  $1/2<b<a<1$, there exists $c=c(a,b)$ such that for every $z_0\in \X$ and $r<r_0$,
\begin{equation}\label{uryshon cond}
\inf_{f\in \mathcal{D}(\overline{A}(z_0, br, ar), A(z_0, r/2, r)) } \sup_{x \in \X} \max(\sL f(x),  \widehat{\sL}f(x)) \le  \frac{c}{\Phi(r)}.
\end{equation}

\smallskip

 This assumption 
 is used in connection with Dynkin's formula in Lemma 2.20 to get 
a scale-invariant estimate of the exit probability.

\begin{definition}
Let $0<\kappa \le 1/2$. We say that an open set $D$ is $\kappa$-fat if there is 
$R_1\in (0,\infty]$ such that for all 
$x \in \overline{D}$ 
and all $r\in(0,R_1)$, 
there is a ball $B(A_r(x),\kappa r) \subset D \cap B(x,r)$. The pair $(R_1,\kappa)$ is called the characteristics of the $\kappa$-fat open set $D$.
\end{definition}

In the remainder of this subsection,
$T>0$ is a fixed constant and $D$ is a fixed $\kappa$-fat open set with characteristics $(R_1,\kappa)$. Without loss of generality, we can 
assume that $R_1\le R_0 \wedge r_0 \wedge r_1$, where $r_1$  is the constant in Corollary \ref{c:18}(1).
For $(t,x)\in(0,T]\times D$, set $r_t=
\Phi^{-1}(t)R_1/(3\Phi^{-1}(T))
 \le R_1/3.$ An open neighborhood $\U(x,t)$ of $x\in D$ and an open ball $\V(x,t)\subset D\setminus \U(x,t)$ are defined as follows:

\begin{figure}[!t]
	\centering
	\includegraphics[width=0.455\columnwidth]{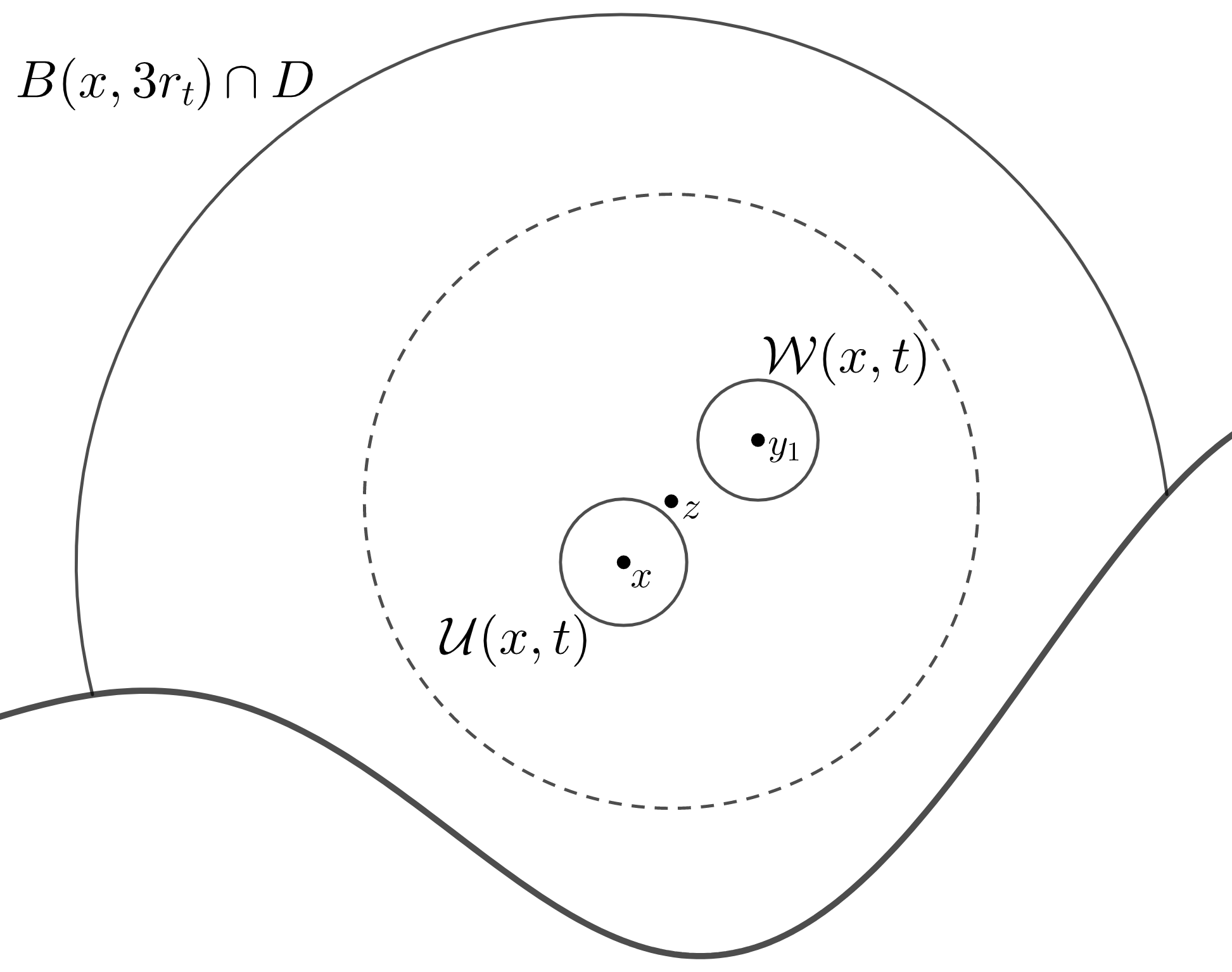} \hspace{5mm}
	\includegraphics[width=0.485\columnwidth]{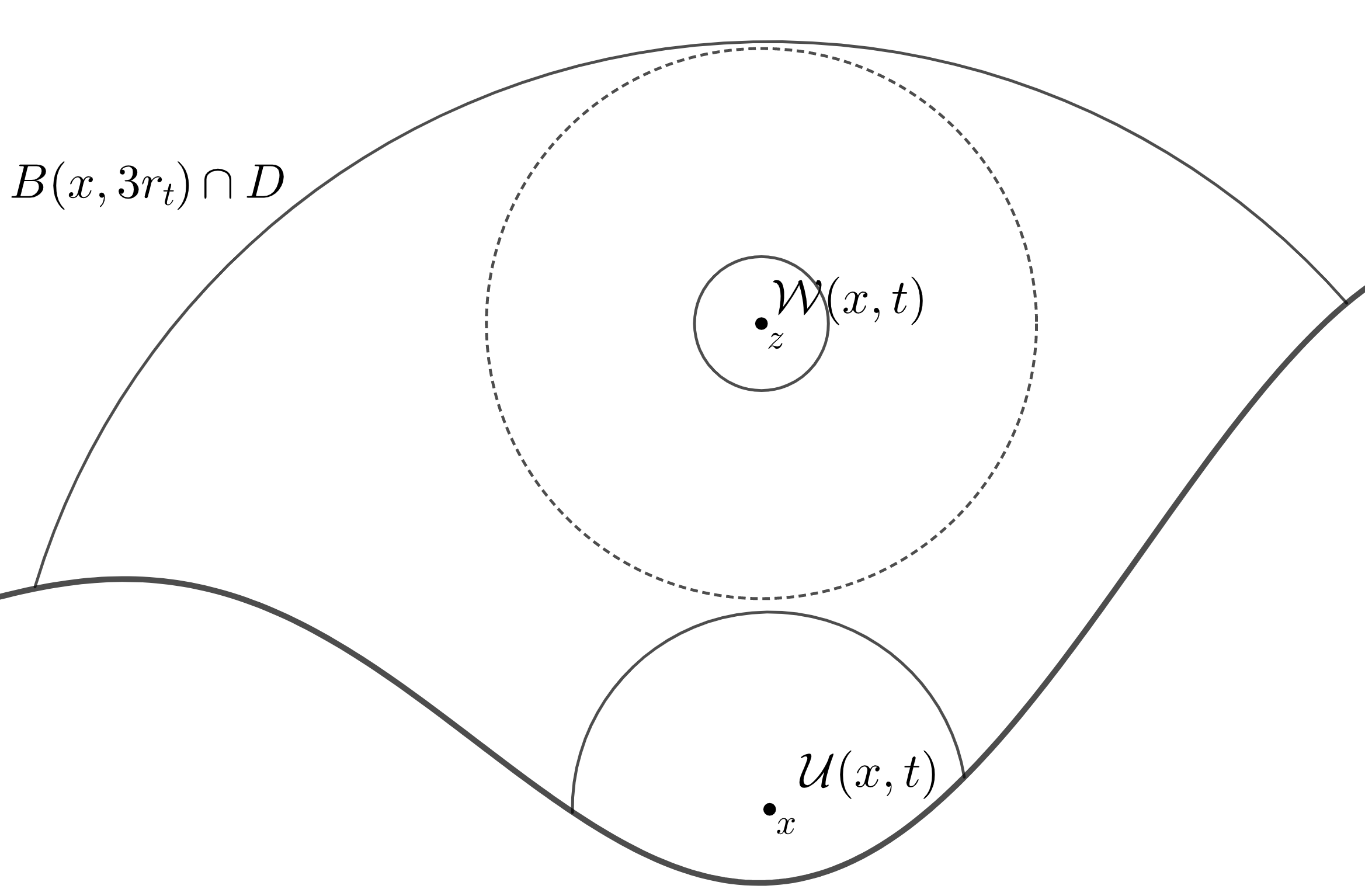}
	Figure 1. $\rho(x,z) \le 3\kappa r_t/2$. \hspace{36mm} Figure 2.  $\rho(x,z) > 3\kappa r_t/2$.
\end{figure}

\smallskip

By the definition of $\kappa$-fat open set, we can find $z=z_{x,t} \in D$ such that $B(z,3\kappa r_t) \subset B(x,3r_t) \cap D$.

(i) If  $\rho (x, z) \le 3 \kappa r_t/2$, we choose 
 $y_1\in \X$ 
such that $\kappa r_t/n_0 \le \rho(x,y_1)\le \kappa r_t$, where 
 $n_0>1$  
is the constant in \eqref{volume condition3}.
Then we define $\U(x,t)=B(x,\kappa r_t/(4n_0))$ and $\V(x,t)=B(y_1,\kappa r_t/(4n_0)).$ We can easily check that $ \U(x,t) \cup\V(x,t) \subset B(x, 3\kappa r_t/2) \subset B(z,3\kappa r_t) \subset D$ and $\U(x,t) \cap \V(x,t) = \emptyset$.

(ii) If $\rho (x, z) > 3 \kappa r_t/2$, we 
define $\U(x,t)=B(x,\kappa r_t) \cap D$ and $\V(x,t)=B(z,\kappa r_t/(4n_0))$.

Note that in either case, we have,
\begin{equation}\label{e:UV}
\kappa r_t/(2n_0) \le \rho(u,v) \le 
 4r_t  
\qquad \text{for all }\; u\in \U(x,t) \text{ and } v\in\V(x,t).
\end{equation}

 See Figures 1 and 2 for some illustration of the sets $\U(x,t)$ and $\V(x,t)$.

It follows from \cite[Theorem I.3.4]{Y96} 
that the L\'evy system of $Y$ is the same as that  of $X$, hence the following L\'evy system formula is valid: for any $f:D\times D \to[0, \infty]$ vanishing on the diagonal and every stopping time 
$S$,
\begin{equation}\label{e:levy_systemY}
\E_x\sum_{t\in (0, S]}f(Y_{t-}, Y_t)= \E_x\int^S_0\int_{D}f(Y_t, z)J(Y_t, z)m(dz)dt.
\end{equation}

Recall that $\tau_{U}=\inf\{s> 0: \, Y_s\notin U\}$ and 
$\wh \tau_{U}=\inf\{s> 0: \, \wh Y_s\notin U\}$.
Note that $\P_x(Y_{\tau_{\U(x,t)}} \in D)=\P_x(\tau_{\U(x,t)} <\zeta)$, 
where $\zeta$ is the lifetime of $Y$.

\begin{lemma}\label{l:22}
Suppose that $\mu \in \mathbf{K}_1(D)$.
  There exists 
  $c=c(T)>0$ such that for all $(t,x)\in(0,T]\times D$ and $z=z_{x,t} \in D$ with $B(z,3\kappa r_t) \subset B(x,3r_t) \cap D$ and  $\rho (x, z) > 3 \kappa r_t/2$, we have 
\begin{equation}
\P_x\big(Y_{\tau_{\U(x,t)}} \in D\big) \le c \,\P_y\big(Y_{\tau_{\U(x,t)}} \in D\big) \frac{\P_x\big(Y_{\tau_{\U(x,t)}} \in \V(x,t)\big)}{\P_y\big(Y_{\tau_{\U(x,t)}} \in \V(x,t)\big)}
\end{equation}
and
\begin{equation}
\wh \P_x\big(\wh Y_{\wh \tau_{\U(x,t)}} \in D\big) \le c \,\wh \P_y\big(\wh Y_{ \wh \tau_{\U(x,t)}} \in D\big) \frac{\wh \P_x\big( \wh Y_{ \wh \tau_{\U(x,t)}} \in \V(x,t)\big)}{\wh \P_y\big( \wh Y_{ \wh\tau_{\U(x,t)}} \in \V(x,t)\big)},
\end{equation}
for every $y\in B(x, \kappa r/2) \cap D$.
\end{lemma}
\smallskip

\pf
Fix $(t,x)\in(0,T]\times D$ and assume that $B(z,3\kappa r_t) \subset B(x,3r_t) \cap D$ and  $\rho (x, z) > 3 \kappa r_t/2$.
Recall that $\U(x,t)=B(x,\kappa r_t) \cap D$ and $\V(x,t)=B(z,\kappa r_t/(4n_0))$.
Define $D_1:=B(x,9\kappa r_t/8) \cap D$ and $D_2:=B(x,9\kappa r_t/8)^c \cap D$. Take any
$$
f\in \mathcal{D}\big(\overline{A}(z, \kappa r_t, 9\kappa r_t/8), A(z, 5 \kappa r_t/8, 5\kappa r_t/4)\big).
$$
Then, 
by Dynkin's formula for $X$ 
(see \cite[(2.11)]{BKK14} and the proof of \cite[(4.6)]{BKK14}), 
 we have for all $y\in B(x, \kappa r_t/2) \cap D$,
\begin{align*}
& \P_y\big(Y_{\tau_{\U(x,t)}} \in D_1\big) = \E_y\left[\;\exp\big(-A^{\mu}_{\tau^X_{\U(x,t)}}\big) : X_{\tau^X_{\U(x,t)}} \in D_1\right]  \\
&\quad \le \E_y\left[f\big(X_{\tau^X_{\U(x,t)}}\big)\; \exp\big(-A^{\mu}_{\tau^X_{\U(x,t)}}\big)\right]-f(y)\\
&\quad= \E_y\left[\int_0^{\tau^X_{\U(x,t)}} \sL f(X_s)\exp(-A^{\mu}_s)ds \right] + \E_y\left[\int_0^{\tau^X_{\U(x,t)}} f(X_s) d\exp(-A^{\mu}_s) \right] \\
&\quad \le \left(\sup_{z \in \X} \sL f(z) \right) \E_y\left[\int_0^{\tau^X_{\U(x,t)}} \exp(-A^{\mu}_s)ds \right] =
\left(\sup_{z \in \X} \sL f(z) \right) \E_y[\tau_{\U(x,t)}]. 
\end{align*}

By Assumption $\mathbf{U}$, taking infimum over $f$ on both sides gives
$$
\P_y(Y_{\tau_{\U(x,t)}} \in D_1) \le \frac{c_1}{\Phi(r_t)}\E_y[\tau_{\U(x,t)}], \quad y\in B(x, \kappa r_t/2) \cap D,
$$
for some constant $c_1>0$.

On the other hand, by \eqref{e:UV} and \eqref{e:levy_systemY}, 
we have that for all $y\in B(x, \kappa r_t/2) \cap D$,
\begin{align*}
&\P_y(Y_{\tau_{\U(x,t)}} \in D_2) 
 \ge \P_y(Y_{\tau_{\U(x,t)}} \in \V(x,t))\\
&= \E_y\left[\int_0^{\tau_{\U(x,t)}} \int_{\V(x,t)} J(Y_s,w)m(dw)ds\right] \\
&\asymp  \E_y[\tau_{\U(x,t)}] \int_{\V(x,t)} \frac{1}{V(y, \rho(y,w))\Phi(\rho(y,w))}m(dw) \asymp \frac{1}{\Phi(r_t)}\E_y[\tau_{\U(x,t)}].
\end{align*}
Thus, using that $J(y,w) \asymp J(v,w )$ for $(w, y, v) \in  D_2 \times (B(x, \kappa r_t/2) \cap D) \times \U(x,t)$, we conclude that for all $y\in B(x, \kappa r_t/2) \cap D$,
\begin{equation*}
\P_y(Y_{\tau_{\U(x,t)}} \in D) \asymp \P_y(Y_{\tau_{\U(x,t)}} \in D_2) \asymp   \frac{1}{\Phi(r_t)} \E_y[\tau_{\U(x,t)}] \asymp  \E_y[\tau_{\U(x,t)}] \int_{D_2}J(y, w)m(dw).
\end{equation*}

Finally,  we have 
\begin{align*}
\frac{\P_x(Y_{\tau_{\U(x,t)}} \in D)}{\P_x(Y_{\tau_{\U(x,t)}} \in \V(x,t))} &\le c_2 \Phi(r_t) \int_{D_2}J(x,w)m(dw)\\
&\le c_2^2 \frac {\P_y(Y_{\tau_{\U(x,t)}} \in D)}{\P_y(Y_{\tau_{\U(x,t)}} \in \V(x,t))},  \quad y\in B(x, \kappa r_t/2) \cap D.
\end{align*}
\qed

Recall that  $\zeta$ is the lifetime of $Y$. We also denote by $\wh \zeta$ the lifetime of $\wh Y$.
\begin{lemma}\label{l:23}
Suppose that $\mu \in \mathbf{K}_1(D)$.
For all $M, T\ge 1$, 
we have that, for all $t \in (0,T)$ and $x \in D$,
\begin{align*}
&\P_x(\zeta>t) \asymp \P_x(\zeta>t/M) \asymp \P_x(\tau_{\U(x,t)}>t) \asymp \P_x(Y_{\tau_{\U(x,t)}} \in D)\\
& \asymp
\P_x(Y_{\tau_{\U(x,t)}} \in \V(x,t)) \asymp  t^{-1} \E_x[\tau_{\U(x,t)}]
\end{align*}
and 
\begin{align*}
&\wh \P_x(\wh \zeta>t) \asymp \wh \P_x(\wh \zeta>t/M) \asymp \wh \P_x(\wh \tau_{\U(x,t)}>t) \asymp \P_x(\wh  Y_{\tau_{\U(x,t)}} \in D) \\
&\asymp \P_x(\wh  Y_{\tau_{\U(x,t)}} \in \V(x,t))  \asymp t^{-1} \wh \E_x[\wh  \tau_{\U(x,t)}],\end{align*}
where $\U(x,t)$  and $\V(x,t)$ are the open sets defined 
in the beginning of this subsection and 
the comparison constants depend only on $d_0,d,\delta_l, \delta_u,T,M,R_1$ and $\kappa.$
\end{lemma}
\smallskip
\pf
Fix $t \in (0,T]$, $x \in D$ and set $r:=r_t=\Phi^{-1}(t)R_1/3\Phi^{-1}(T).$ 

\noindent
Case (i):  $\rho (x, z) \le 3 \kappa r/2$.
By \eqref{e:c181}, we have
$$
 1\ge \P_x(\zeta>t/M) \ge \P_x(\zeta>t) \ge \P_x(\tau_{\U(x,t)}>t)=\P_x(\tau_{B(x,\kappa r/(4n_0))}>t)\ge c>0.
$$

On the other hand, by \eqref{e:UV},  \eqref{e:levy_systemY} and 
\eqref{e:c1820},
\begin{align*}
1&\ge \P_x(Y_{\tau_{\U(x,t)}} \in D) \ge \P_x(Y_{\tau_{\U(x,t)}} \in \V(x,t))= \E_x\left[\int_0^{\tau_{\U(x,t)}} \int_{\V(x,t)}J(Y_s,v)m(dv)ds\right]\\
&\ge c_1
\frac{m(\V(x,t))\E_x[\tau_{\U(x,t)}]}{V(x, 3r)\Phi(3r)} \ge c_2  \frac{\E_x[\tau_{\U(x,t)}]}{\Phi( 3r)}\ge c_3  t^{-1}\E_x[\tau_{B(x,\kappa r/(4n_0)}]\ge c_4.
\end{align*}
Therefore, we arrive at the assertion of the lemma in this case.\\

\noindent
Case (ii): $\rho (x,z) > 3 \kappa r/2$. Note that
\begin{align*}
\P_x(\zeta>t/M) & \le \P_x(\tau_{\U(x,t)}>t/M)+\P_x(Y_{\tau_{\U(x,t)}} \in D) \le  M  t^{-1}\E_x[\tau_{\U(x,t)}] +\P_x(Y_{\tau_{\U(x,t)}} \in D).
\end{align*}
Fix a $y\in D$ such that $B(y, \kappa^2 r/2) \subset B(x, \kappa r/2) \cap D$.
Then by Lemma \ref{l:22} we have
$$
\P_x(Y_{\tau_{\U(x,t)}} \in D) \le c_5 \P_y(Y_{\tau_{\U(x,t)}} \in D) \frac{\P_x(Y_{\tau_{\U(x,t)}} \in \V(x,t))}{\P_y(Y_{\tau_{\U(x,t)}} \in \V(x,t))}. 
$$
By \eqref{e:UV}, \eqref{e:levy_systemY} and 
\eqref{e:c1820},
\begin{align*}
&\P_y(Y_{\tau_{\U(x,t)}} \in \V(x,t)) \ge \P_y(Y_{\tau_{B(y,\kappa^2 r/2)}} \in \V(x,t))\\
&= \E_y\left[\int_0^{\tau_{B(y,\kappa^2 r/2)}} \int_{\V(x,t)}J(Y_s,v)m(dv)ds\right] \ge c_6 \frac{V(x, r)\E_y[\tau_{B(y,\kappa^2 r/2)}]}{V(x, 3r)\Phi(3r)} \ge c_7
\end{align*}
and
\begin{align*}
\P_x(Y_{\tau_{\U(x,t)}} \in \V(x,t))\asymp \frac{V(x,r)\E_x[\tau_{\U(x,t)}]}{V(x, r)\Phi(r)}=c_8t^{-1}\E_x[\tau_{\U(x,t)}].
\end{align*}
\smallskip
It follows that 
$$
\P_x(\zeta>t/M) \le c_9  t^{-1}\E_x[\tau_{\U(x,t)}]\asymp \P_x(Y_{\tau_{\U(x,t)}} \in \V(x,t)) \asymp \P_x(Y_{\tau_{\U(x,t)}} \in D).
$$

Note that $B(x, (3-2\kappa)r) \cap D \supset \U(x,t) \cup \V(x,t)$ for every $(t,x) \in (0,T] \times D$.
Thus by \eqref{e:c181},
\begin{align*}
&\P_x(\zeta>t) \ge \P_x(\tau_{B(x, 3r) \cap D}>t)\ge 
\E_x\left[\inf_{w\in \V(x,t)} \P_w(\tau_{B(x, 3r) \cap D}>t) : Y_{\tau_{\U(x,t)}} \in \V(x,t)\right] \\
&\ge \E_x\left[\inf_{w\in \V(x,t)} \P_w(\tau_{B(w,\kappa r)}>t) : Y_{\tau_{\U(x,t)}} \in \V(x,t)\right]\ge c_{10}\P_x(Y_{\tau_{\U(x,t)}} \in \V(x,t)).
\end{align*}
The proof is now complete. \qed

\begin{theorem}\label{t:f1} 
Let $D$ be a $\kappa$-fat set with characteristics $(R_1, \kappa)$.
Suppose that $\mu \in \mathbf{K}_1(D)$.
Then for
all $T>0$,
there exists $c \ge 1$ such that for all $(t,x,y) \in (0,T) \times D \times D$,
\smallskip
\begin{equation}
 c^{-1} \P_x(\zeta>t) \wh  \P_y( \wh \zeta>t) \wt q(t,x,y) \le q^D(t,x,y)\le c \P_x(\zeta>t) \wh  \P_y(\wh \zeta>t) \wt q(t,x,y).
\end{equation}

\end{theorem}

\pf
Fix $t \in (0,T]$ and set $r:=\Phi^{-1}(t)R_1/(3\Phi^{-1}(T))$.
 
(1) We first prove the upper bound.
By  Lemma \ref{l:23},
for any $x,y \in D$ with $\rho(x,y) \le 4r$, we have
\begin{align*}
&q^D(t/2,x,y)  = \int_D q^D(t/4,x,w) q^D(t/4,w, y) m(dw) \\
& \le C_0\int_D q^D(t/4,x,w)  \wt q(t/4,w, y) m(dw) \\
& \le c_1  \P_x(\zeta>t/4) V(y,\Phi^{-1}(t))^{-1}  \le c_2   \P_x(\zeta>t) p(t/2,x,y).
\end{align*}

Now, we assume $\rho(x,y) > 4r$. Let $U_1:=\U(x,t)$ be the set defined before, $U_3:=\{w \in D : \rho(x,w) > \rho(x,y)/2 \}$, and $U_2:=D \setminus (U_1 \cup U_3)$. Since $x \in U_1$, $y \in U_3$ and $U_1 \cap U_3 = \emptyset$, by the strong Markov property, we have
\begin{align*}
q^D(t/2,x,y) & = \E_x[q^D(t/2-\tau_{U_1} , Y_{\tau_{U_1}} , y) : \tau_{U_1} < t/2, Y_{\tau_{U_1}} \in U_2] \\
& \quad +\E_x[q^D(t/2-\tau_{U_1} , Y_{\tau_{U_1}} , y) : \tau_{U_1} < t/2, Y_{\tau_{U_1}} \in U_3] =: I + II.
\end{align*}

First, note that for every $u \in U_2$, $\rho(u,y) \ge \rho(x,y)-\rho(x,u) \ge \rho(x,y)/2$, which implies that
$$
V(y, \rho(x,y)) \le V(u, \rho(x,y)+\rho(u,y)) \le V(u, 3 \rho(u,y)).$$
 Therefore, using \eqref{volume condition}, for all $(s,u) \in (0,t/2] \times U_2$, 
\begin{align*}
& q^D(s,u,y) \le c_3 \wt q(s,u,y) \le  c_4 \frac{t}{V(y,\rho(x,y)) \Phi(\rho(x,y))} \le c_5 p(t/2,x,y).
\end{align*}
Now it follows from Lemma \ref{l:23} that
$$
I \le c_5 p(t/2,x,y) \P_x(Y_{\tau_{U_1}} \in D) \asymp \P_x(\zeta>t) p(t/2,x,y).
$$

On the other hand, for all $u \in U_1$ and $w \in U_3$, we have $\rho(u,x) \le \kappa r < 4^{-1}\kappa \rho(x,y)$ and $\rho(u,w) \ge \rho(x,w)-\rho(x,u)
\ge  \rho(x,y)/2-\kappa r \ge \rho(x,y)/4$, which implies that
$$V(x, \rho(x,y)) \le V(u, \rho(x,y)+\rho(u,x))  \le V(u, (1+\kappa/4) \rho(x,y))  \le V(u, 8 \rho(u, w)).$$
 Thus, by \eqref{volume condition}, \eqref{e:scale_Phi2}, Lemma \ref{l:23} and the L\'evy system formula, and using the assumption $\rho(x,y)>4r$, 
\begin{align*}
& II = \int_0^{t/2} \int_{U_1} \int_{U_3} q^{U_1} (s,x,u) J(u,w) q^D(t/2-s,w,y) m(dw) m(du) ds \\
& \quad \le c_6 \frac{1}{V(x,\rho(x,y))\Phi(\rho(x,y))} \int_0^{t/2} \P_x(\tau_{U_1}>s)
\wh\P_y( \wh \zeta>t/2-s) ds \\
& \quad \le c_6 \frac{1}{V(x,\rho(x,y))\Phi(\rho(x,y))} \int_0^{\infty} \P_x(\tau_{U_1}>s)  ds \\
& \quad = 2c_6 \frac{t/2}{V(x,\rho(x,y))\Phi(\rho(x,y))}  t^{-1}\E_x[\tau_{U_1}] \asymp \P_x(\zeta>t) 
p(t/2,x,y).
\end{align*}
Thus for all $x,y \in D$,
$
q^D(t/2,x,y) \le c_7 \P_x(\zeta>t) p(t/2,x,y),
$
and, similarly, for all $x,y \in D$,
$
q^D(t/2,x,y) \le c_7 \wh \P_y( \wh \zeta>t) p(t/2,x,y).
$

Finally, by the semigroup property, we conclude that
\begin{align*}
& q^D(t,x,y) = \int_D q^D(t/2,x,w) q^D(t/2,w,y)m(dw) \\
& \quad \le c_7^2 \P_x(\zeta>t) \wh \P_y( \wh \zeta>t) \int_{\X} p(t/2,x,w)p(t/2,w,y) m(dw) \\
& \quad \le c_7^2 C_0 \P_x(\zeta>t) \wh  \P_y( \wh \zeta>t) \wt q(t,x,y).
\end{align*}

(2)
For the lower bound, we use the notation $\V$ as before. By the semigroup property,
\begin{align*}
& q^D(t,x,y) = \int_D \int_D q^D(t/3,x,u) q^D(t/3,u,w) q^D(t/3,w,y) m(dw) m(du) \\
& \qquad \ge \int_{\V(x,t/3)} \int_{\V(y,t/3)} q^D(t/3,x,u) q^D(t/3,u,w) q^D(t/3,w,y) m(dw) m(du).
\end{align*}
First, observe that for all $(u,w) \in \V(x,t/3) \times \V(y,t/3)$,
$$
\delta_D(u), \delta_D(w) \ge \kappa (a_l/3)^{1/\delta_l} r/4, \quad (\rho(x,y) -   6(3a_u)^{1/\delta_u}r)_+ \le \rho(u,w) \le \rho(x,y) +  6  (3a_u)^{1/\delta_u}r. 
$$ 
Here is an explanation of the last inequality above, the others being similar. By the triangle inequality 
and symmetry, it suffices to show that $\rho(u,x)\le 3(3a_u)^{1/\delta_u}r$. 
Since $\V(x,t/3) \subset B(x, 3r_{t/3})$, this will be so provided that $r_{t/3}\le (3a_u)^{1/\delta_u}r$. But this immediately follows from \eqref{e:scale_inv_Phi} by estimating $\Phi^{-1}(t/3)/\Phi^{-1}(t)$.
By considering cases 
$\rho(x,y) >  12(3a_u)^{1/\delta_u}r $ and $\rho(x,y) \le   12(3a_u)^{1/\delta_u}r $ separately,
we get  from Theorem \ref{t:sharpest_qD} and \eqref{e:tildeq}    that for all $(u,w) \in \V(x,t/3) \times \V(y,t/3)$, 
$$
q^D(t/3,u,w) \asymp \wt q(t/3,u,w) \asymp \left(\frac{1}{V(u, \Phi^{-1}(t))} \wedge \frac{t}{V(u, \rho(w,u))\Phi(\rho(w,u))} \right) \asymp \wt q(t,x,y).
$$ 

Next, let $c_9:= \kappa (a_l/3)^{1/\delta_l}/8$.
  By Theorem \ref{t:sharpest_qD} and  \eqref{e:tildeq},    for all $(s, u) \in (t/6, t/3) \times \V(x,t/3)$ and $w \in B(u,c_9  r 
 )$,  we have 
$q^D(s,u,w) \asymp \wt q(s,u,w) \asymp V(u, r)^{-1}.
$ 
Moreover, by \eqref{e:UV}  and \eqref{e:c},   for all $u \in \V(x,t/3)$ and 
$(v, w) \in \U(x,t/3) \times B(u,c_9  r 
)$,
$$J(v,w) \asymp \frac1{V(v, r) \Phi(r)}\asymp  \frac1{V(x, r) \Phi(r)}.$$
Thus, by the L\'evy system formula and Lemma \ref{l:23}, for every $u \in \V(x,t/3)$, we have
\begin{align*}
& q^D(t/3,x,u) \ge 
\E_x[q^D(t/3-\tau_{\U(x,t/3)} , Y_{\tau_{\U(x,t/3)}} , u) : \tau_{\U(x,t/3)} < t/3, Y_{\tau_{\U(x,t/3)}} \in B(u,c_9r)]\\
& \ge  \int_0^{t/6} \int_{\U(x,t/3)} \int_{B(u,c_9r)} q^{\U(x,t/3)}(s,x,v) J(v,w) q^D(t/3-s,w,u) m(dw) m(dv) ds \\
&  \ge c_{10} \frac{1}{V(x,r)\Phi(r)} \int_{t/6}^{t/3} \int_{B(u,c_9r)} \P_x(\tau_{\U(x,t/3)}>s) V(u, r)^{-1}m(dw)ds \\
& \ge c_{11} \frac{t}{V(x,r)\Phi(r)} \P_x(\tau_{\U(x,t/3)}>t/3) \asymp \frac{1}{V(x,r)} \P_x(\zeta>t).
\end{align*}
Similarly for $w \in  \V(y,t/3)$, 
$ q^D(t/3,w,y) \ge c_{12} \frac{1}{V(y, r)} \wh \P_y(\wh \zeta>t).$
Therefore, we conclude that 
\begin{align*}
q^D(t,x,y) &\ge c_{13} \wt q(t,x,y) \int_{\V(x,t/3)} q^D(t/3,x,u) m(du) \int_{\V(y,t/3)} q^D(t/3,w,y) m(dw) \\
&\ge c_{14} \P_x(\zeta>t) \wh \P_y(\wh \zeta>t)\wt q(t,x,y).
\end{align*}
\qed

Using Theorem \ref{t:sharpest_qD} and Corollary \ref{c:18}(2), the following global estimates can be proved by the same argument. We omit the proof.  
\begin{theorem} \label{t:f2}
Let $D$ be a $\kappa$-fat set with characteristics $(\infty, \kappa)$.
Suppose that $\mu \in \mathbf{K}_1(D)$ and
$R_0=m(\X)=\wt T=r_0=\infty$, 
where $r_0$ is the constant in Assumption {\bf U}. 
Then
 there exists $c_1(\kappa)>0$ such that for 
all
  $(t,x,y) \in (0, \infty) \times D \times D$,
\smallskip
\begin{equation}\label{e:f2}
 c_1^{-1} \P_x(\zeta>t) \wh  \P_y( \wh\zeta>t) \wt q(t,x,y) \le q^D(t,x,y)\le c_1 \P_x(\zeta>t) \wh  \P_y(\wh \zeta>t) \wt q(t,x,y).
\end{equation}

\end{theorem}

\bigskip

\begin{example}\label{e:3_1}

Suppose that $(\X, \rho, m)$ is an 
unbounded Ahlfors regular $n$-space for some 
$n\in (0, \infty)$,  that is, 
for all $x \in \X$ and $r \in (0, 1]$, $m(B(x,r)) \asymp r^n.$ 
Assume that $\rho$ is uniformly equivalent to the shortest-path metric in $\X$.
Suppose that
there is a diffusion process $\xi$ with a symmetric, continuous transition density
$p^\xi(t,x, y)$ satisfying the following sub-Gaussian bounds 
\begin{align}\label{e:subgaussian}
&\frac{c_1}{t^{n/d_w}}\exp\left(-c_2\left(\frac{\rho(x, y)^{d_w}}{t}\right)^{1/(d_w-1)}\right)
\le p^\xi(t, x, y) \nonumber\\
&\,\,\,\, \le \frac{c_3}{t^{n/d_w}}\exp\left(-c_4\left(\frac{\rho(x, y)^{d_w}}{t}\right)^{1/(d_w-1)}\right),
\end{align}
for all $x, y\in \X$ and $t\in (0, \infty)$. Here $d_w\ge 2$ is the walk dimension of the space $\X$.
Examples of $\xi$ include 
Brownian motions on unbounded Riemannian manifolds, Brownian motions on Sierpinski gaskets, Sierpinski carpets or more general fractals.
Let $\alpha\in (0, d_w)$ and let 
$T$ be an $(\alpha/d_w)$-stable subordinator independent of $\xi$.
We define a process $X$ by $X_t=\xi_{T_t}$. 
Then $X$ is a symmetric Feller process.
It is easy to check that $X$ has a transition density $p(t, x, y)$ satisfying 
\begin{equation}\label{e:desnitySD}
p(t, x, y) \asymp \left(t^{-\frac{n}{\alpha}}\wedge \frac{t}{\rho(x, y)^{n+\alpha}}\right),
\end{equation}
for all $x, y\in \X$ and $t\in (0, \infty)$.
It follows from \cite[Appendix A]{BKK14} that Assumptions {\bf A} and {\bf U} above  are also satisfied with $\Phi(r)=r^\alpha$.
Therefore, by Theorem \ref{t:f1} and \eqref{e:desnitySD}, if $D$ is a $\kappa$-fat open set in  $\X$ and 
$\mu \in \mathbf{K}_1(D)$, 
then 
 for all $(t,x,y) \in (0,r_0) \times D \times D$,
$$ q^D(t,x,y)\asymp \P_x(\zeta>t)   \P_y( \zeta>t) \left(t^{-\frac{n}{\alpha}}\wedge \frac{t}{\rho(x, y)^{n+\alpha}}\right).
$$
\end{example}

\section{Dirichlet heat kernel estimates of regional fractional Laplacian with critical killing}\label{s:3}

In this section we assume that $d\ge 2$, $\X$ is either the closure of a 
$C^{1, 1}$ open subset $D$ of $\R^d$ or $\R^d$ itself, and 
the underlying process is either a reflected $\alpha$-stable
process in $\overline{D}$ (or  a non-local perturbation of it), 
or an $\alpha$-stable process in $\R^d$ (or a drift perturbation of it).
We investigate Dirichlet heat kernel estimates
under critical killing. We first recall the definition of reflected $\alpha$-stable processes.

Let $0<\alpha< 2$ and $\sA(d,-\alpha)=\alpha 2^{\alpha-1} \pi^{-d/2} \Gamma((d+\alpha)/2)\Gamma(1-\alpha/2)^{-1}$. Here 
$\Gamma$ is the gamma function defined by 
$\Gamma(\lambda):=\int^{\infty}_0 t^{\lambda-1} e^{-t}dt$, $\lambda > 0$. 
For a $C^{1,1}$ open subset $D$ of $\R^d$, let $(\overline{\EE}, \overline{\FF})$ be the
Dirichlet space on $L^2(D, dx)$ defined by
\begin{eqnarray*}
\overline{\FF}&:=& \left\{ u\in L^2(D);\,\int_D\int_D \frac{(u(x)-u(y))^2}
{|x-y|^{d+\alpha}}\, dxdy<\infty \right\}, \\
\overline{\EE} (u, v)&:=&  \frac{1}{2} {\cal A} (d, \, -\alpha)
\int_D \int_D\frac{(u(x)-u(y))(v(x)-v(y))}
{|x-y|^{d+\alpha}}\, dxdy, \quad u, v \in \overline{\FF}.
\end{eqnarray*}
It is  well known that 
$W^{\alpha/2, 2}(D)=\overline{\FF}$ and the Sobolev norm
$ \|\cdot \|_{\alpha/2, 2;D}$ is equivalent to $ \sqrt{\overline{\EE}_1}$ where
$\overline{\EE}_1:=\overline{\EE} +(\,\cdot\, ,\, \cdot\, )_{L^2(D)}$.
As noted in \cite{BBC03}, $(\,\overline{\EE},\, \overline{\FF}\,)$ is a regular
 Dirichlet form on
$\overline D$ and
its associated Hunt process $X$ lives on $\overline D$.
 We call the process $X$  a reflected $\alpha$-stable
process in $\overline{D}$. When $D$ is the whole $\R^d$, $X$ is simply an $\alpha$-stable process.

It folows from \cite{CK03} that
$X$  admits a strictly positive and jointly continuous transition density $p(t,x,y)$ with respect to the Lebesgue measure $dx$ and that
\begin{align}\label{e:est_p_Dr}
p(t,x,y) \asymp \left(t^{-d/\alpha} \wedge \frac{t}{|x-y|^{d+\alpha}}\right) \qquad (t,x,y) \in (0,1)\times \overline{D}\times \overline{D}.
\end{align}

When $\alpha \in (1,2)$, the killed  process $X^D$ is the censored stable process in  $D$.
When $\alpha\in (0, 1]$, it follows from \cite[Section 2]{BBC03} that, starting from inside $D$, the process $X$ 
neither hits nor approaches $\partial D$ at any finite time. 
Thus, 
the killed process $X^D$ is simply $X$ restricted to $D$ (without killing). 

We will see that, for all $\alpha\in (0, 2)$, the killed isotropic $\alpha$-stable process $Z^D$ 
can be obtained from $X^D$
through a Feynman-Kac perturbation of the form \eqref{e:rsFK} with $\kappa$ satisfying \eqref{e:C1}.

It follows from \cite{CKS10b} that, when  $\alpha\in (1,2)$, the transition density 
$p_D^X(t,x,y)$
of $X^D$ has the following estimates:
\begin{align}\label{e:est_p_Dc1}
p_D^X(t,x,y)
\asymp \left(1 \wedge \frac{\delta_D(x)}{ t^{1/\alpha}}\right)^{\alpha-1} \left(1 \wedge \frac{\delta_D(y)}{ t^{1/\alpha}}\right)^{\alpha-1} \left(t^{-d/\alpha} \wedge \frac{t}{|x-y|^{d+\alpha}}\right)
\end{align}
for $(t,x,y) \in (0,1)\times D\times D$. 

It follows from \cite{CKS10a} that the transition density 
$p_D^Z(t,x,y)$ of 
$Z^D$ has the following estimates:
\begin{align}\label{e:est_p_Dc}
 p_D^Z(t,x,y) 
\asymp \left(1 \wedge \frac{\delta_D(x)}{t^{1/\alpha}}\right)^{\alpha/2} \left(1 \wedge \frac{\delta_D(y)}{t^{1/\alpha}}\right)^{\alpha/2} \left(t^{-d/\alpha} \wedge \frac{t}{|x-y|^{d+\alpha}}\right) 
\end{align}
for $(t,x,y) \in (0,1)\times D\times D$. 

In Subsection \ref{sc:6.1}, we will establish 
explicit Dirichlet heat kernel estimates under critical killing, which also provides 
 an alternative and unified proof of \eqref{e:est_p_Dc1} and \eqref{e:est_p_Dc}.

In Subsection \ref{ss-nl}, we consider non-local perturbations of $(\,\overline{\EE},\, \overline{\FF}\,)$ when $D$ is a bounded $C^{1,1}$ open set. Subsection \ref{ss-eorigin} 
covers the case $D=\R^d\setminus \{ 0\}$ and drift perturbations.

\subsection{$C^{1,1}$ open set}\label{sc:6.1}

 In this subsection, we assume that  $D$ 
 is a $C^{1, 1}$ open set in $\R^d$ with
characteristics ($R_2$, $\Lambda$), and that $X$ is a reflected $\alpha$-stable process
in $\overline{D}$.  
Without loss of generality, we will always assume that $\Lambda\ge 1$.
It is easy to check that the process $X$ satisfies the assumptions in Subsection \ref{s:setup}
and Assumption {\bf U}.

Let $\R^d_+:=\{y=(y_1, \dots, y_d) \in \R^d: y_d >0\}$. 
For $d\ge 2$ and $p\in (-1, \alpha)$, 
we define 
 $w_p(y) = (y_d)^p$ for 
$y \in \R^d_+$
and $w_p(y)=0$ otherwise. According to \cite[(5.4)]{BBC03}, we have
\begin{align}\label{e:upperhalfspace}
\sA(d,-\alpha)\; \lim_{\epsilon \downarrow 0} \int_{\R^d_+, |y-z|>\epsilon} \frac{w_p(y)-w_p(z)}{|y-z|^{d+\alpha}} dy = C(d,\alpha,p)z_d^{p-\alpha}, \quad z \in \R^d_+,
\end{align}
where $C(d,\alpha,p):= \sA(d,-\alpha) \frac{\omega_{d-1}}{2} \beta(\frac{\alpha+1}{2}, \frac{d-1}{2}) \gamma(\alpha,p)$,
$\beta(\cdot, \cdot)$ is the beta 
function,
$\omega_{d-1}$ is the $(d-2)$-dimensional Lebesgue measure of the unit sphere in $\R^{d-1}$ and 
$$\gamma(\alpha,p)=\int_0^1 \frac{(t^p-1)(1-t^{\alpha-p-1})}{(1-t)^{1+\alpha}} dt.$$ 
Observe that
$$
\frac{d\gamma(\alpha,p)}{dp} = \int_0^1 \frac{(t^{\alpha-p-1}-t^p)|\log t|}{(1-t)^{1+\alpha}}dt
$$
is positive for $p>(\alpha-1)/2$ and thus $p\mapsto \gamma(\alpha, p)$ is 
strictly
increasing on 
$((\alpha-1)/2, \alpha)$. 
Moreover, we have 
\begin{align}
\label{e:palpha}
C(d,\alpha,\alpha-1)=C(d,\alpha,0)=0 \quad \text{and} \quad \lim_{p \uparrow \alpha} C(d,\alpha,p) = \infty. 
\end{align}

Let $\sH_\alpha$ be the collection of  non-negative functions $\kappa$ on $D$ 
with the property  that
there exist constants $C_1,C_2 \ge 0 $ and $\eta \in [0, \alpha)$ such that $\kappa (x) \le C_2$ for all $x\in D$ with $\delta_D(x) \ge 1$ and 
\begin{align}\label{e:C1}
|\kappa(x)-C_1 \delta_D(x)^{-\alpha}| \le C_2 \delta_D(x)^{-\eta},
\end{align}
for all $x \in D$ with $\delta_D(x)<1$. If $\alpha \le 1$, then we further assume that $C_1>0$.
It follows from \eqref{e:palpha} that we can find a unique 
$p\in [\alpha-1, \alpha)\cap (0, \alpha)$ such that $C_1 = C(d,\alpha,p)$. 
For any $p\in [\alpha-1, \alpha)\cap (0, \alpha)$, define
\begin{align}
\label{e:palphap}
\sH_\alpha (p):=\{\kappa \in \sH_\alpha:  \text{ the constant } C_1 \text{ in \eqref{e:C1} is } C(d,\alpha,p)\}.
\end{align}
Note that $\sH_\alpha=\cup_{p\in [\alpha-1, \alpha)\cap (0, \alpha)} \sH_\alpha (p)$. 
We fix a $\kappa \in \sH_\alpha (p)$ and let $Y$  be a Hunt process on $D$ corresponding to the Feynman-Kac semigroup of $X^D$  through the 
multiplicative 
functional $e^{-\int_0^t \kappa (X^D_s)ds}$.
That is, 
\begin{align}\label{e:rsFK}
\E_x\left[  f(Y_t)\right] = \E_x\left[ e^{-\int_0^t \kappa (X^D_s)ds} f(X^D_t)\right], \quad t\ge 0 , x\in D.
\end{align}
Since, by Example \ref{e:1}, $\kappa(x)dx \in \mathbf{K}_1(D)$, 
it follows from Theorem \ref{t:f1} that
 $Y$ has a transition density $q^D(t,x,y)$ with the following estimate:
\begin{align}\label{e:qDf}
q^D(t,x,y) \asymp \P_x(\zeta>t) \P_y(\zeta>t)\left[t^{-d/\alpha} \wedge \frac{t}{|x-y|^{d+\alpha}}\right],
\end{align} 
for $(t,x,y) \in (0,1] \times D \times D$.
To get  explicit estimate of $\P_x(\zeta>t)$, 
we will estimate 
$\P_x(Y_{\tau_{\U(x,t)}} \in D)$
and use Lemma \ref{l:23}.

For $f \in C_c^2(D)$, define 
\begin{equation}\label{e:defL}
L_\alpha f(x) :=\sA(d,-\alpha)\; \lim_{\epsilon \downarrow 0} \int_{D, |y-x|>\epsilon}
\frac{f(y)-f(x)}{|y-x|^{d+\alpha}} dy \quad \text{ and} \quad 
Lf(x) := L_\alpha f(x)-\kappa(x)f(x).
\end{equation}
The operator $L$ coincides with the restriction to $C_c^2(D)$ of the generator of the transition semigroup of $Y$ in $C_0(D)$.

\begin{lemma}\label{l:3.1}
Let $0 <p\le q < \alpha$ and  define
$$h_q(x):=\delta_D(x)^q.$$
Then there exist $A_1=A_1(q,d,\alpha,\Lambda,C_2,\eta, R_2)>0$ and $A_2=A_2(q,d,\alpha,\Lambda,C_2,\eta, R_2)\in (0,{(R_2 \wedge 1)}/{4})$ such that the following inequalities hold:

(i) If $q>p$, then
$$
A_1^{-1} \delta_D(x)^{q-\alpha} \le Lh_q(x) \le A_1 \delta_D(x)^{q-\alpha}
$$
for every $x \in D$ with $0<\delta_D(x)<A_2$.

(ii) If $q=p$, then
$$
|Lh_p(x)| \le A_1 ( \delta_D(x)^{p-\eta} + |\log \delta_D(x)|)
$$
for every $x \in D$ with $0<\delta_D(x)<A_2$.
\end{lemma}

\pf 
Without loss of generality, we assume $R_2=1$.
Let $x \in D$ with $\delta_D(x) <1/4$.
Choose a point $z \in \partial D$ such that $\delta_D(x) = |x-z|$. Then there exist a $C^{1,1}$ function $\psi : \R^{d-1} \to \R$ such that $\psi(z)=\nabla \psi(z)=0$ and an orthonormal coordinate system $CS_z$ such that
$$
D \cap B(z,1) = \{ y = (\ty,y_d) \mbox{ in } CS_z: 
y_d > \psi(\ty) \} \cap B(z,1),
$$
and  $z=0$ and $x=(\tx,x_d)=(0,x_d)$ in $CS_z$.
Note that
\begin{align*}
L h_q(x) = L_\alpha  h_q(x) - C(d,\alpha,p)x_d^{q-\alpha} - (\kappa(x) - C(d,\alpha,p)x_d^{-\alpha})h_q(x) =: I - II - III.
\end{align*}
By our assumption, we have
$
|III| \le C_2 x_d^{q-\eta}.
$

For any open subset $U \subset \R^d$, define $\kappa_U(x) = \sA(d,-\alpha)\int_{U^c} \frac{dy}{|y-x|^{d+\alpha}}$.
Recall that $w_q(y) = (y_d)^q$ for $y \in \R^d_+$
and $w_q(y)=0$ otherwise.
Since $h_q(x)=w_q(x) = x_d^q$, by \eqref{e:upperhalfspace} we have 
\begin{align*}
&I  = \sA(d,-\alpha)\; \lim_{\epsilon \downarrow 0} \left[ \int_{\R^d, |y-x|>\epsilon} \frac{h_q(y)-h_q(x)}{|y-x|^{d+\alpha}} dy\right] + \kappa_D(x)h_q(x) \\
& = \sA(d,-\alpha)\; \lim_{\epsilon \downarrow 0} \left[ \int_{\R^d, |y-x|>\epsilon} \frac{h_q(y)-w_q(y)}{|y-x|^{d+\alpha}} dy + \int_{\R^d, |y-x|>\epsilon} \frac{w_q(y)-w_q(x)}{|y-x|^{d+\alpha}} dy\right] +\kappa_D(x)w_q(x)  \\
& = C(d,\alpha,q) x_d^{q-\alpha} + \sA(d,-\alpha)\; \lim_{\epsilon \downarrow 0} \left[ \int_{\R^d, |y-x|>\epsilon} \frac{h_q(y)-w_q(y)}{|y-x|^{d+\alpha}} dy\right] + (\kappa_D(x)-\kappa_{\R^d_+}(x))w_q(x).
\end{align*}

According to \cite[Lemma 5.6]{BBC03}, if $1<\alpha<2$, then there is a constant $c=c(d,\alpha,  \Lambda )$ such that
$
|\kappa_D(x)-\kappa_{\R^d_+}(x)| \le c x_d^{1-\alpha}.
$
By a similar calculation as in \cite[Lemma 5.6]{BBC03}, one can show that for $\alpha \le 1$, 
$
|\kappa_D(x)-\kappa_{\R^d_+}(x)| \le c (|\log x_d|{\bf 1}_{\{\alpha = 1\}}+1).
$
Thus, for any $0<\alpha<2$, we get
\begin{align}
\label{e:new1}
|(\kappa_D(x)-\kappa_{\R^d_+}(x))w_q(x)| \le c x_d^q (x_d^{1-\alpha}+|\log x_d|) \le c.
\end{align}
Therefore, it remains to bound
$I_{\epsilon}:=\int_{\R^d, |y-x|>\epsilon} \frac{h_q(y)-w_q(y)}{|y-x|^{d+\alpha}} dy$.
Since $D$ is a $C^{1,1}$ open set, it satisfies the inner and outer ball conditions. 
Thus we can assume that $B_1= B(\mathbf{e}_d,1) \subset D$ and $B_2=B(-\mathbf{e}_d,1) \subset D^c$ where $\mathbf{e}_d:=(\tilde{0}, 1)$.
We define 
$E:=\{y=(\ty,y_d) : |\ty|<1/4, |y_d|<1/2  \}$, $E_1:=\{ y \in E : y_d>2|\ty|^2\}$ and $E_2:=\{ y \in E : -y_d>2|\ty|^2\}$. Then it is easy to see that $E_1 \subset B_1 \cap E \subset D$ and $E_2 \subset B_2 \cap E \subset D^c$. Thus, since $h_q(y)=w_q(y)=0$ for $y \in E_2$.
\begin{align*}
I_{\epsilon}& = \int_{E^c, |y-x|>\epsilon} \frac{h_q(y)-w_q(y)}{|y-x|^{d+\alpha}} dy +\int_{E_1, |y-x|>\epsilon} \frac{h_q(y)-w_q(y)}{|y-x|^{d+\alpha}} dy\\
& + \int_{E\setminus (E_1 \cup E_2), |y-x|>\epsilon} \frac{h_q(y)-w_q(y)}{|y-x|^{d+\alpha}} dy =: J_{1,\epsilon} +  J_{2,\epsilon} +  J_{3,\epsilon}.
\end{align*}

First,
 since $x=(0, x_d)$ and $x_d = \delta_D(x)<1/4$, we can see that $|y-x| \ge |y|/2$ for all $y \in E^c$. 
Indeed, if $y\in E^c$, then either $|y_d| \ge 1/2$ or $|\wt y| \ge 1/4$. If $|y_d| \ge 1/2>2x_d$, then $|y-x|^2 = |\wt y|^2 + |y_d- x_d|^2 \ge |\wt y|^2 + |y_d|^2/4 \ge |y|^2/4$. If 
$|\wt y| \ge 1/4$, then 
$$
|y-x|^2 - \frac{1}{4}|y|^2 = \frac{3}{4}|\wt y|^2+ \frac{3}{4}|y_d-\frac{4}{3}x_d|^2 - \frac{1}{3}x_d^2 
 \ge\frac{3}{4}|\wt y|^2 - \frac{1}{3}x_d^2  \ge \frac{3}{64}- \frac{1}{48}>0.
$$
Thus, since $|h_q(y)-w_q(y)| \le 2|y|^q$, we get  
$$|J_{1,\epsilon}| \le 2^{1+d+\alpha}\int_{E^c} |y|^{q-d-\alpha}dy \le c \int_{1/4}^{\infty}l^{q-\alpha-1}dl=c.$$

For $y \in E_1$, we have $\delta_D(y) \le \delta_{B_2^c}(y) \le y_d+1-\sqrt{1-|\ty|^2} \le y_d+|\ty|^2 < 2y_d$. Thus, 
\begin{align*}
J_{2,\epsilon} &\le \int_{E_1, |y-x|>\epsilon} \frac{(y_d+|\ty|^2)^q-y_d^q}{|y-x|^{d+\alpha}}dy \le \int_{E_1, |y-x|>\epsilon} \frac{q(1 \vee 2^{q-1})|\ty|^2y_d^{q-1}}{|y-x|^{d+\alpha}}dy\\
& \le cx_d^{q+1-\alpha} \int_{B(0,1/x_d)} \frac{|\tu|^2 u_d^{q-1}}{|u-\mathbf{e}_d|^{d+\alpha}} du.
\end{align*}
We have used the change of the variables $y = x_d u$ in the last inequality above. Note that
\begin{align*}
\int_{B(0,1/x_d)} \frac{|\tu|^2 u_d^{q-1}}{|u-\mathbf{e}_d|^{d+\alpha}} du & = \int_{B(0,2)} \frac{|\tu|^2 u_d^{q-1}}{|u-\mathbf{e}_d|^{d+\alpha}} du + \int_{B(0,1/x_d) \setminus B(0,2)} \frac{|\tu|^2 u_d^{q-1}}{|u-\mathbf{e}_d|^{d+\alpha}} du \\
&\le c \int_{B(0,2)} |u-\mathbf{e}_d|^{2-d-\alpha} du + c \int_{B(0,1/x_d) \setminus B(0,2)} |u|^{q+1-d-\alpha} du \\
&\le c\left( \int_0^2 l^{1-\alpha} dl + \int_2^{x_d^{-1}} l^{q-\alpha} dl \right).
\end{align*}
 It follows that
$
J_{2,\epsilon} \le c(1+|\log x_d|).
$
Besides, for $y \in E_1$, we have
$$\delta_D(y) \ge \delta_{B_1}(y) \ge 1-\sqrt{(1-y_d)^2+|\ty|^2} \ge 1-\sqrt{1-(y_d-|\ty|^2)} \ge \frac{1}{2}(y_d-|\ty|^2)>\frac{1}{4}y_d$$
and
$$ y_d-(1-\sqrt{(1-y_d)^2+|\ty|^2})
=\frac{|\ty|^2}{\sqrt{(1-y_d)^2+|\ty|^2}+1-y_d}
  \le \frac{|\ty|^2}{2(1-y_d)} \le |\ty|^2.$$  Therefore, by the mean value theorem, 
\begin{align*}
J_{2,\epsilon} &\ge \int_{E_1, |y-x|>\epsilon} \frac{( 1-\sqrt{(1-y_d)^2+|\ty|^2})^q-y_d^q}{|y-x|^{d+\alpha}}dy \\
&\ge - q\int_{E_1, |y-x|>\epsilon} \frac{|\ty|^2 \Big( \sup_{\lambda \in [0, |\ty|^2]}( 1-\sqrt{(1-y_d)^2+\lambda})^{q-1}\Big)}{|y-x|^{d+\alpha}}dy \\
&\ge -q(4^{1-q} \vee 1)\int_{E_1, |y-x|>\epsilon} \frac{|\ty|^2y_d^{q-1}}{|y-x|^{d+\alpha}}dy.
\end{align*}
Thus, we have
$
|J_{2,\epsilon}| \le c(1+|\log x_d|).
$

Lastly, let $m_{d-1}(dx)$ be the $(d-1)$-dimensional  Hausdorff  
measure.  Then there exists a constant $c>0$ such that for all $0<l<1$,
$$
m_{d-1}(\{y : |\ty| = l, -2|\ty|^2 \le y_d \le 2|\ty|^2 \}) \le cl^d.
$$
Since $|h_q(y)|,\, |w_q(y)| \le 4^q|\ty|^{2q}$ for $y \in E \setminus (E_1 \cup E_2)$,
\begin{align*}
|J_{3,\epsilon}| \le c \int_0^{1/4} \int_{|\ty|=l, y \in E \setminus (E_1 \cup E_2)} l^{2q-d-\alpha} m_{d-1}(dy)dl \le c \int_0^{1/4} l^{2q-\alpha} dl \le c.
\end{align*}
Combining the above estimates, we conclude that  $|I_{\epsilon}| \le c(1+\log|x_d|)$.

If $q>p$, we note that $C(d,\alpha,q)>C(d,\alpha,p)$ and $q-\alpha< 0 \wedge (q-\alpha+1) \wedge (q-\eta)$. 
\qed

\smallskip

Fix a $q\in (p,\alpha)$ such that $q<p-\eta+\alpha$. Then define $A_3:=A_1(p) \vee A_1(q)$, $A_4:=A_2(p)\wedge A_2(q)$, where $A_1$ and $A_2$ are 
the constants in Lemma \ref{l:3.1}, and 
\begin{equation*}
v_1(x):=h_p(x)+h_q(x). 
\end{equation*}
By Lemma \ref{l:3.1},
for any $x \in D$ with $\delta_D(x) <A_4$, we have
$$
Lv_1(x) \ge A_3^{-1}\delta_D(x)^{q-\alpha}-A_3(\delta_D(x)^{p-\eta} + |\log \delta_D(x)|).
$$
Thus, there exist $A_5 \in (0,A_4)$ and $A_6>0$ 
such that
\begin{equation}\label{e:Lv1}
Lv_1(x) \ge  2  A_6 \delta_D(x)^{q-\alpha} \quad \text{for all} \;\; x \in D \;\; \text{with} \;\; \delta_D(x)<A_5.
\end{equation}
Define $v_2(x):=h_p(x)-\frac{1}{2}h_q(x)$. By the same argument, we can find $A_7 \in (0,A_4)$ and $A_8>0$ such that $Lv_2 (x)\le - 2  A_8 \delta_D(x)^{q-\alpha}$ for all $x\in D$ with $\delta_D(x)<A_7$.\\

Now, we are ready to estimate $\P_x(Y_{\tau_{\U(x,t)}} \in D)$. 
We continue to assume $R_2=1$.
Note that  $D$ is a $\kappa$-fat open set with 
characteristics $(1, \kappa)$.
 Recall that $r_t$ is defined as $r_t=\Phi^{-1}(t)R_1/(3\Phi^{-1}(T))$ in Subsection \ref{s:factorization}. Since in the current setting $\Phi(t)=t^{1/\alpha}$,  we can take $r_t=t^{1/\alpha}/3$ in the definition of $\U(x,t)$. 
 Let $A_9 \in (0, A_5/2]$ be a  constant which will be chosen later. 
Without loss of generality we assume $\kappa < A_7  \wedge A_9$. 

Fix $(t,x) \in (0,1] \times D$. 
If $\delta_D(x) \ge \kappa t^{1/\alpha}/3$, 
then we have $\P_x(Y_{\tau_{\U(x,t)}} \in D) \asymp 1$ in view of Lemma \ref{l:23} and
\eqref{e:c1820}.
Recall that $z_{x,t} \in D$ is a point such that $B(z, 3\kappa r_t) \subset B(x,3r_t) \cap D$.  
 Assume that $\delta_D(x) < \kappa t^{1/\alpha}/3$. In this case, we have $|x-z_{x,t}| \ge \delta_D(z_{x,t})-\delta_D(x) > \kappa t^{1/\alpha} - \kappa t^{1/\alpha}/3 >  \kappa t^{1/\alpha}/2$  and hence we should choose the second definition of $\U(x,t)$ so that $\U(x,t)=B(x,\kappa t^{1/\alpha}/3) \cap D$. Let $w \in \partial D$ be the point such that $|x-w|=\delta_D(x)$. 
 Define $D^{{\rm bdry}}(l):=\{y \in D: |y-w|<l\}$ and 
$D^{{\rm int}}(l):= \{ y \in D^{{\rm bdry}}(2): \delta_D(y)>l\}$.
Note that 
$\U(x,t)  \subset   D^{{\rm bdry}}(A_9)  \subset  D^{{\rm bdry}}(2)$.  Indeed, for every $y \in \U(x,t) \subset D$, we have $|y-w| \le |y-x| + \delta_D(x) \le \kappa/3 + \kappa/3 
< A_9<1/8$.

Let $\varphi \in C_c^\infty(\R^d)$ be a non-negative radial function such that $\varphi(y) = 0$ for $|y| > 1$ and $\int_{\R^d} \varphi(y)dy = 1$. For $k \ge 1$, define $\varphi_k(y):=6^{kd}\varphi(6^k y)$ and 
$f_k:=\varphi_k * (v_1 {\bf 1}_{D^{{\rm int}}(5^{-k})})$.
 Since $6^{-k} < 5^{-k}$, we have $f_k \in C_c^\infty(D)$ and hence $Lf_k$ is well defined everywhere. 
  Pick any $z \in \U(x,t)$ 
 (hence $\delta_D(z)<|z-w|<A_9$ by the calculation above)  
such that 
 $\delta_D(z)>\max\{2^{-k/(q-p)}, 2^{-pk/(d+q)}\}=:a_k$ 
and observe that
\begin{align*}
Lf_k(z) &= L(\varphi_k * v_1)(z) - L(\varphi_k * v_1 - f_k)(z) \\
&= L(\varphi_k * v_1)(z) + \kappa(z) (\varphi_k *v_1 - f_k)(z)  \\
& \quad - \sA(d,-\alpha)\; \lim_{\epsilon \downarrow 0} \int_{D, |y-z|>\epsilon}  
\frac{(\varphi_k * v_1)(y)-f_k(y)-(\varphi_k * v_1)(z)+f_k(z)}{|y-z|^{d+\alpha}}dy\\
&=: M_1(z) + M_2(z) + M_3(z)= M_1 + M_2 + M_3.
\end{align*}

To bound $M_1$, we need some preparation.  
For $|u|<6^{-k}$, let $w_u \in \partial D$ be a point such that $\delta_D(z-u)=|z-u-w_u|$. Note that by the triangle inequality and the assumption that $\delta_D(z) >a_k \ge 2^{-k}$, we have
\begin{equation}\label{e:triangle}
(3^k-1)|u|<(1-3^{-k})\delta_D(z) \le \delta_D(z)-|u| \le \delta_D(z-u) \le \delta_D(z)+|u| \le (1+3^{-k})\delta_D(z).
\end{equation}
Let $\psi_u:\R^{d-1} \to \R$ be a $C^{1,1}$ function and $CS_{w_u}$ an orthonormal coordinate system with origin at $w_u$ such that $\psi_u(\widetilde{0})=0$, $\nabla \psi_u(\widetilde{0})=\widetilde{0}$, $\lVert \nabla \psi_u \rVert_\infty \le \Lambda$,
the coordinate of $z-u$ in $CS_{w_u}$ is $(\wt0, \delta_D(z-u))$ and 
$D \cap B(w_u,1) = \{ y^u = (\ty^u,y^u_d) \mbox{ in } CS_{w_u}:  y^u_d > \psi_u(\ty^u) \} \cap B(w_u,1).$
  Define $D-u:=\{y-u: y \in D\}$ for $u \in \R^d$. Using the coordinate system   
$CS_{w_u}$, we have that 
  for all  $ q_0\in[p,\alpha)$, 
 $\epsilon \in (0,1)$ and $|u|<6^{-k}$, 
\begin{align*}
&\left|\int_{D-u, |y-(z-u)|>\epsilon} \frac{h_{q_0}(y)-h_{q_0}(z-u)}{|y-(z-u)|^{d+\alpha}}dy \right| \\
&\le \left|\int_{B(z-u, \epsilon)^c} \frac{h_{q_0}(y^u)-\delta_D(z-u)^{q_0}}{|y^u-(z-u)|^{d+\alpha}}dy^u \right| + \int_{(D-u)^c} \frac{|h_{q_0}(y^u)-\delta_D(z-u)^{q_0}|}{|y^u-(z-u)|^{d+\alpha}}dy^u \\
&\le \left|\int_{B(z-u, \epsilon)^c} \frac{h_{q_0}(y^u)-(y^u_d \vee 0)^{q_0}}{|y^u-(z-u)|^{d+\alpha}}dy^u \right| + \left|\int_{B(z-u, \epsilon)^c} \frac{(y^u_d \vee 0)^{q_0}-\delta_D(z-u)^{q_0}}{|y^u-(z-u)|^{d+\alpha}}dy^u \right|\\
& \quad  +  \int_{B(z-u, \delta_D(z))^c} \frac{|y^u|^{q_0} + \delta_D(z-u)^{q_0}}{|y^u-(z-u)|^{d+\alpha}}dy^u  \\
&=:N_1(z,u,\epsilon)+N_2(z,u,\epsilon)+N_3(z,u)=N_1+N_2+N_3.
\end{align*}

According to the proof of Lemma \ref{l:3.1} and \eqref{e:triangle}, we can see that for all $|u|<6^{-k}$,
$$
 N_1  \le c_1(1+ 
|\log \delta_D(z-u)|) \le c_2(1+ \log(3/2)+ |\log \delta_D(z)|) 
\quad \text{uniformly  in } \epsilon \in (0,1).
$$ 
Moreover, by \cite[p.120-121]{BBC03} and \eqref{e:triangle}, we obtain 
$ N_2  \le c_3 \delta_D(z-u)^{{q_0}-\alpha} \le c_4\delta_D(z)^{{q_0}-\alpha}$ 
 uniformly  in  $\epsilon \in (0,1)$. Lastly, 
 using the triangle inequality $|y^u|\le |y^u-(z-u)^u|+|(z-u)^u|=|y^u-(z-u)^u|+\delta_D(z-u)$,
 we also have 
 \begin{align*}
  N_3  &\le c_5 \int_{\delta_D(z)}^\infty \big((l + \delta_D(z-u))^{q_0} + \delta_D(z-u)^{q_0} \big)l^{-\alpha-1} dl \\
 & \le c_5 \int_{\delta_D(z)}^\infty \big((l + 2\delta_D(z))^{q_0} + 2^{q_0}\delta_D(z)^{q_0} \big)l^{-\alpha-1} dl \le c_6\delta_D(z)^{{q_0}-\alpha}.
 \end{align*}
Thus, we conclude that for all  $q_0 \in [p,\alpha)$ 
there exists $c_7=c_7(q_0)>0$ such that for all $|u|<6^{-k}$, 
\begin{equation}\label{e:dominant}
\left|\int_{D-u, |y-(z-u)|>\epsilon} \frac{h_{q_0}(y)-h_{q_0}(z-u)}{|y-(z-u)|^{d+\alpha}}dy \right| \le c_7\delta_D(z)^{{q_0}-\alpha} \quad \text{uniformly  in }  \epsilon \in (0,1).
\end{equation}

On the other hand, we also observe that by \eqref{e:triangle}, for 
$|u|<6^{-k}$,
\begin{align*}
&\int_{(D-u) \setminus D} \frac{dy}{|y-(z-u)|^{d+\alpha}} \le 
\int_{((D-u) \setminus D) \cap B(w_u, 1)} \frac{dy^u}{|y^u-(z-u)|^{d+\alpha}} + \int_{B(w_u, 1)^c} \frac{dy^u}{|y^u-(z-u)|^{d+\alpha}} \\
&\le \int_{|\wt y^u|<\delta_D(z)/(2\Lambda)}\int_{\psi_u(\wt y^u)-|u|}^{\psi_u(\wt y^u)} (\delta_D(z-u)-y_d^u)^{-(d+\alpha)}dy_d^ud\wt y^u  \\
& \quad + \int_{\delta_D(z)/(2\Lambda) \le |\wt y^u| \le 1 }\int_{\psi_u(\wt y^u)-|u|}^{\psi_u(\wt y^u)} |\wt y^u|^{-(d+\alpha)}dy_d^ud\wt y^u + 2^{d+\alpha}\int_{B(0, 1)^c} |y|^{-d-\alpha} dy  \\
& \le |u|\Big(\big( \delta_D(z-u)-\delta_D(z)/2\big)^{-(d+\alpha)}\int_{ |\wt y^u|<\delta_D(z)/(2\Lambda)} d\wt y^u + \int_{ |\wt y^u| \ge \delta_D(z)/(2\Lambda)} |\wt y^u|^{-(d+\alpha)}d\wt y^u \Big) + 
c \\
& \le A_{10}(6^{-k}\delta_D(z-u)^{-\alpha-1} + 1),
\end{align*}
for some constant $A_{10}>0$. In the second inequality above, we used the facts that 
the coordinates of $z-u$ in $CS_{w_u}$  are  $(\wt 0, \delta_D(z-u))$,
and for all $y^u \in B(w_u,1)^c$, we have $|y^u-(z-u)| \ge |y^u-w_u| - \delta_D(z-u) \ge 2^{-1} |y^u|$ since $w_u=0$ in $CS_{w_u}$ and $\delta_D(z-u)<2\delta_D(z)\le 2A_9 \le  A_5<A_4<1/2$. Besides, in the third inequality above, we used the fact that  for all $|\wt y^u| < \delta_D(z)/(2\Lambda)$, we have $|\psi_u(\wt y^u)| \le \lVert \nabla \psi_u \rVert_\infty |\wt y^u| \le  \delta_D(z)/2$. 
Since $2^{-k/(q-p)} \vee 2^{-k} \le a_k<\delta_D(z)<A_9$ and $\lim_{k \to \infty} 6^{-k} a_k^{p-q-1} \le \lim_{k \to \infty}6^{-k} (2^{-k/(q-p)})^{p-q} (2^{-k})^{-1} \le
\lim_{k \to \infty}6^{-k}  4^k =0$, 
 by taking $A_9<(A_6/(6\sA(d,-\alpha)A_{10}))^{1/(p+\alpha-q)}$, for all $k$ large enough, we have
 \begin{align}\label{e:Lv1per}
&\sA(d, -\alpha)\lim_{\epsilon \downarrow 0}\int_{D-u, |y-(z-u)|>\epsilon} \frac{v_1(y)-v_1(z-u)}{|y-(z-u)|^{d+\alpha}}dy - \kappa(z-u)v_1(z-u) \nn\\
& = Lv_1(z-u) +  \sA(d,-\alpha) \left( \int_{D \setminus (D-u)} \frac{v_1(z-u)-v_1(y)}{|y-(z-u)|^{d+\alpha}}dy -  \int_{(D-u) \setminus D} \frac{v_1(z-u)}{|y-(z-u)|^{d+\alpha}}dy \right) \nn\\[2pt]
& \ge 2A_6 \delta_D(z-u)^{q-\alpha}
- \sA(d,-\alpha)A_{10}\big(6^{-k}\delta_D(z-u)^{-\alpha-1} + 1\big) v_1(z-u) \nn\\[2pt]
& \ge 2A_6 \delta_D(z-u)^{q-\alpha}  - 2\sA(d,-\alpha)A_{10}\big(6^{-k}\delta_D(z-u)^{p-\alpha-1} + \delta_D(z-u)^p\big)  \nn\\[2pt]
& =2\Big(A_6 - \sA(d,-\alpha)A_{10}\big(6^{-k} \delta_D(z-u)^{p-q-1} + \delta_D(z-u)^{p-q+\alpha}\big)  \Big)\delta_D(z-u)^{q-\alpha} \nn\\
& \ge 2 \Big(A_6 - 2\sA(d,-\alpha)A_{10}\big(6^{-k} a_k^{p-q-1} + A_9^{p-q+\alpha}\big)  \Big)\delta_D(z-u)^{q-\alpha} \ge A_6 \delta_D(z-u)^{q-\alpha}.
\end{align}
The first inequality above is valid since for all $|u|<6^{-k}$ and  $y \in D \setminus(D-u)$, by \eqref{e:triangle}, it holds 
that $\delta_D(y) \le |u| \le \delta_D(z-u)$, implying $v_1(z-u)\ge v_1(y)$.
Moreover, for $k$ large enough, $\delta_D(z-u)\le A_9+6^{-k}<A_5$ so we could use \eqref{e:Lv1}. We used the fact that $v_1(z-u) \le 2\delta_D(z-u)^p$ in the second 
 and \eqref{e:triangle} in the third inequality above.

 Now, 
 since the support of $\varphi_k$ is contained in $B(0, 6^{-k})$,
we have that
 for all sufficiently large $k$,  
\begin{align*}
M_1 & =\lim_{\epsilon \downarrow 0} \sA(d,-\alpha)\int_{D, |y-z|>\epsilon} \int_{\R^d}\varphi_k(u)\frac{v_1(y-u)-v_1(z-u)}{|y-z|^{d+\alpha}}dudy\\
&  \quad  - \int_{\R^d}\kappa(z)\varphi_k(u)v_1(z-u) du\\ 
&=\lim_{\epsilon \downarrow 0} \int_{\R^d}\varphi_k(u)\left(\sA(d,-\alpha)\int_{D-u, |y-(z-u)|>\epsilon} \frac{v_1(y)-v_1(z-u)}{|y-(z-u)|^{d+\alpha}}dy - \kappa(z-u)v_1(z-u) \right)du\\
 &\quad   + \int_{\R^d} (\kappa(z-u) - \kappa(z))\varphi_k(u) v_1(z-u)du \\
 &\ge A_6\int_{\R^d} 
 (\delta_D(z)+|u|)^{q-\alpha}
  \varphi_k(u)du  + C_1\int_{\R^d} (\delta_D(z-u)^{-\alpha}-\delta_D(z)^{-\alpha}) \varphi_k(u)v_1(z-u)du \\
 & \quad - C_2  \int_{\R^d} (\delta_D(z-u)^{-\eta} + \delta_D(z)^{-\eta}) \varphi_k(u)v_1(z-u)du \\[1pt]
 & \ge A_6(1+3^{-k})^{q-\alpha}\delta_D(z)^{q-\alpha}  - C_1(1- (1+3^{-k})^{-\alpha}) \delta_D(z)^{-\alpha}(\varphi_k * v_1)(z)\\[2pt]
 &\quad - C_2 (1 + (1-3^{-k})^{-\eta}) \delta_D(z)^{-\eta}(\varphi_k * v_1)(z).
\end{align*}
 In the second equality above we used Fubini's theorem and the change of variables. In the first inequality we first used the dominated convergence theorem (which is applicable due to  \eqref{e:dominant}) and then  \eqref{e:Lv1per} and \eqref{e:C1}. In the second inequality above, we used \eqref{e:triangle}.

 Since $(\varphi_k * v_1)(z) \le 2(1+3^{-k})^q \delta_D(z)^p $, $q<p-\eta + \alpha$ and  $2^{-k/(q-p)} \le a_k <\delta_D(z)<A_9$, 
 by taking   $A_9< (A_6/(432C_2))^{1/(p-\eta+\alpha-q)}$, 
 for all $k$ large enough, 
we have 
\begin{align*}
 M_1 &\ge \frac{7}{8}A_6 \delta_D(z)^{q-\alpha}  - 3 \big(\alpha C_1 3^{-k} \delta_D(z)^{p-q} + 3C_2 \delta_D(z)^{p-\eta+\alpha-q} \big)  \delta_D(z)^{q-\alpha}  \\
&\ge \left( \frac{7}{8}A_6 - 3\alpha C_1 3^{-k} a_k^{p-q} - 9C_2 A_9^{p-\eta+\alpha-q} \right) \delta_D(z)^{q-\alpha} \ge \frac{5}{6}A_6 \delta_D(z)^{q-\alpha}.
\end{align*}

Note that for every $k \ge 2$, $u \in B(0, 6^{-k})$ and $y \in D$ such that $\delta_D(y)>4^{-k}$ and $|y-w| \le 1$, we have $\delta_D(y-u) \ge 4^{-k} - 6^{-k} >5^{-k}$ and $|y-u-w| \le |y-w| + |u| < 2$ and therefore 
\begin{equation}\label{smallper}
1 - {\bf 1}_{D^{{\rm int}}(5^{-k})}(y-u)  = 0.
\end{equation} 
 In particular, since $\varphi_k$ is supported in $B(0, 6^{-k})$, 
  $\delta_D(z)>2^{-pk/(d+q)}>4^{-k}$  
 and $|z-w| \le |z-x| + |x-w| < 2t^{1/\alpha}/3< 1$, for all $k \ge 2$, we have
\begin{align*}
M_2 = \kappa(z) \int_{\R^d}\big(1-{\bf 1}_
{D^{{\rm int}}(5^{-k})}(z-u) \big) v_1(z-u) \varphi_k(u)   du = 0.
\end{align*}

Finally, using \eqref{smallper},
 by taking $A_9$ sufficiently smaller than $A_6$, for all $k$ large enough, we have
\begin{align*}
&|M_3| \le \sA(d,-\alpha) \lim_{\epsilon \downarrow 0} \int_{D, |y-z|>\epsilon} \int_{\R^d} \varphi_k(u) 
\frac{\big(1-{\bf 1}_{D^{{\rm int}}(5^{-k})}(y-u) \big)v_1(y-u)}{|y-z|^{d+\alpha}}dudy \\
& \le c_1 \left( \int_{D, \delta_D(y) \le 4^{-k}} \int_{\R^d} \varphi_k(u) 
\frac{\delta_D(y-u)^p}{|y-z|^{d+\alpha}}dudy + \int_{D, |y-w|>1} \int_{\R^d}\varphi_k(u)\frac{(\delta_D(y-u)+ 1)^q}{|y-z|^{d+\alpha}} du dy \right) \\
& \le c_1 \left( \int_{D, \delta_D(y) \le 4^{-k}} \int_{\R^d} \varphi_k(u) 
\frac{(\delta_D(y)+|u|)^p}{|y-z|^{d+\alpha}}dudy + \int_{|y-w|>1} \int_{\R^d}\varphi_k(u)\frac{ (|y-w|+ |u|+1)^q}{(3^{-1}|y-w|)^{d+\alpha}} du dy \right)  \\
&  \le c_2 \left( \int_{D, \delta_D(y) \le 4^{-k}}
 \frac{ 4^{-pk} }{|y-z|^{d+\alpha}}   dy \int_{\R^d} \varphi_k(u) du +  \int_{1}^{\infty} l^{q-\alpha-1} dl \int_{\R^d} \varphi_k(u)du \right)  \\
& \le  c_2 \left( \int_{D, \delta_D(y) \le 4^{-k}, |y-z| \le 1}
 \frac{ 4^{-pk} }{(\delta_D(z)-\delta_D(y))^{d+\alpha}}   dy  +  \int_{|y-z| >1}
 \frac{ 4^{-pk} }{|y-z|^{d+\alpha}}   dy + \frac{1}{\alpha-q} \right)  \\[2pt]
 & \le  c_3 \big(4^{-pk}\delta_D(z)^{-d-\alpha} +1 \big) \le c_3 \big(4^{-pk}2^{pk} +A_9^{\alpha-q} \big) \delta_D(z)^{q-\alpha} \le \frac{A_6}{2}\delta_D(z)^{q-\alpha}.
\end{align*}
 In the second inequality above, we have used the facts that  $\delta_D(y)^q \le \delta_D(y)^p$ for $\delta_D(y)<1$ and $\delta_D(y)^q +1 \ge \delta_D(y)^p$ for  all $y$, since $q>p$. 
In the third inequality, we first estimate $|z-w| \le |z-x|+|x-w| <2/3< (2/3)|y-w|$ by using that $|y-w|>1$, which implies that 
$|y-z|\ge |y-w|-|z-w|\ge (1/3)|y-w|$. The estimate $\delta_D(y-u) \le \delta_D(y) + |u| \le |y-w| + |u|$, follows by the choice of $w\in \partial D$. 
In the fourth inequality, we have used  the fact that the support of $\varphi_k$ is contained in $B(0, 6^{-k})$.
 Besides, we have used  the fact that $\int_{\R^d}\varphi_k(u)du=1$ in the fifth inequality and $\delta_D(z) \le |z-w|<A_9$ in the seventh inequality, and the assumption that $\delta_D(z)>a_k\ge 2^{-pk/(d+q)}$ in the sixth and seventh inequalities. 
 
Thus, we conclude that,  for all sufficiently large $k$,
$Lf_k(z) \ge 3^{-1}A_6 \delta_D(z)^{q-\alpha} \ge 0$ for all $z \in \U(x,t)$ such that $\delta_D(z)>a_k$.
Recall that $f_k \in C_c^\infty(D)$ and hence contained in the domain of the generator of $Y$. Thus, by Dynkin's formula, we have that for all sufficiently large $k$,
\begin{align*}
f_k(x) &= \E_x\big[f_k(Y_{\tau_{\U(x,t) \cap D^{{\rm int}}(a_k)}})\big] - 
\E_x \left[\int_0^{\tau_{\U(x,t) \cap D^{{\rm int}}(a_k)}} Lf_k(Y_t)dt \right] \le 
\E_x\big[f_k(Y_{\tau_{\U(x,t) \cap D^{{\rm int}}(a_k)}})\big].
\end{align*}
Since  
$f_k=\varphi_k * (v_1 {\bf 1}_{D^{{\rm int}}(5^{-k})}) 
\to v_1 {\bf 1}_{D^{{\rm bdry}}(2)} \le v_1 $ pointwise 
and $Y_{\tau_{\U(x,t)\cap D^{\rm int, 2}(a_k)}}\to Y_{\tau_{\U(x,t)}}$ (using $\U(x,t) \subset  D^{{\rm bdry}}(2)$),
it follows from the bounded convergence theorem,
\begin{align*}
&\delta_D(x)^p \le v_1(x)  = \lim_{k \to \infty}f_k(x) \le \lim_{k \to \infty} 
\E_x\big[f_k(Y_{\tau_{\U(x,t) \cap D^{{\rm int}}(a_k)}})\big] \\
&= \E_x\big[v_1(Y_{\tau_{\U(x,t)}}) : Y_{\tau_{\U(x,t)}}  \in D^{{\rm bdry}}(2) \big]
\le \E_x\big[v_1(Y_{\tau_{\U(x,t)}})\big].
\end{align*}

Recall that  we have assumed $\kappa  <  A_7\wedge A_9  <A_5 \wedge A_7$.
Set 
$r=r(t):=(A_5 \wedge A_7) t^{1/\alpha}>  \kappa t^{1/\alpha}$.
 Note that for every $n \ge 1$ and $u \in D^{{\rm bdry}}(2^n r)$, we have $v_1(u) \le (\delta_D(x) + 2^n r)^p+ (\delta_D(x) + 2^n r)^q \le 2^{(n+1)p}r^p + 2^{(n+1)q}r^q \le  2^{(n+1)q+1}r^p$. Thus, we have 
\begin{align*}
& \E_x\big[v_1(Y_{\tau_{\U(x,t)}})\big] \\
\le &
\E_x\big[v_1(Y_{\tau_{\U(x,t)}}) : Y_{\tau_{\U(x,t)}} \in D^{{\rm bdry}}(r)\big]+
\sum_{n=0}^{\infty} \E_x\big[v_1(Y_{\tau_{\U(x,t)}}) : Y_{\tau_{\U(x,t)}} \in 
D^{{\rm bdry}}(2^{n+1}r) \setminus D^{{\rm bdry}}(2^n r)\big] \\
\le& c_0 r^p\P_x\big(Y_{\tau_{\U(x,t)}} \in D^{{\rm bdry}}(r)\big)
+c_0\sum_{n=0}^{\infty} 2^{(n+1)q+1}r^p \P_x\big(Y_{\tau_{\U(x,t)}} \in  
D^{{\rm bdry}}(2^{n+1}r) \setminus D^{{\rm bdry}}(2^n r)\big)
\end{align*}
and that for every $n \ge 0$,
\begin{align*}
&\P_x\big(Y_{\tau_{\U(x,t)}} \in D^{{\rm bdry}}(2^{n+1}r) \setminus D^{{\rm bdry}}(2^n r)\big)\\
&\quad \le c_1 \E_x \int_0^{\tau_{\U(x,t)}} 
\int_{ D^{{\rm bdry}}(2^{n+1}r) \setminus D^{{\rm bdry}}(2^n r))} |Y_s-z|^{-d-\alpha} dz ds \\
&\quad \le c_2 (2^{n+1}r)^d (2^n r)^{-d-\alpha} \E_x\big[\tau_{\U(x,t)}\big]=c_3 2^{-n\alpha} r^{-\alpha} \E_x\big[\tau_{\U(x,t)}\big].
\end{align*}
Since 
$$
\P_x\big(Y_{\tau_{\U(x,t)}} \in D^{{\rm bdry}}(r)\big) 
\ge c_4 \E_x \left[\int_0^{\tau_{\U(x,t)}} 
\int_{ D^{{\rm bdry}}(r)} |Y_s-z|^{-d-\alpha} dz ds \right] 
\ge c_5 r^{-\alpha}\E_x[\tau_{\U(x,t)}],$$
we deduce that
\begin{align*}
\delta_D(x)^p & \le c_0r^p\P_x\big(Y_{\tau_{\U(x,t)}} \in D^{{\rm bdry}}(r)\big)+
c_6 \sum_{n=0}^{\infty} 2^{(n+1)q-n\alpha}r^p \P_x\big(Y_{\tau_{\U(x,t)}} \in 
D^{{\rm bdry}}(r)
\big)\\
&\le c_7 r^p\P_x\big(Y_{\tau_{\U(x,t)}} \in D^{{\rm bdry}}(r)\big)
\le c_7 r^p\P_x\big(Y_{\tau_{\U(x,t)}} \in D\big),
\end{align*}
where in the second inequality we used the fact that $q < \alpha$.

By applying the similar argument to the function $g_k:=\varphi_k * (v_2 {\bf 1}_{D^{{\rm int}}(5^{-k})})$, we also have that
\begin{align*}
\delta_D(x)^p \ge v_2(x) &= \lim_{k \to \infty} g_k(x) \ge \lim_{k \to \infty} \E_x\big[g_k(Y_{\tau_{\U(x,t) \cap D^{{\rm int}}(a_k)}})\big] \\
&= \E_x [ (v_2 {\bf 1}_{D^{{\rm bdry}}(2)}) (Y_{\tau_{\U(x,t)}})] \ge \frac{1}{2}r^p \P_x(Y_{\tau_{\U(x,t)}} \in \V(x,t)) .
\end{align*}
The last inequality holds since $\V(x,t) \subset D^{{\rm bdry}}(2)$.

Therefore, in view of Lemma \ref{l:23}, we get 
$\P_x\big(\zeta>t\big) \asymp  \left(\frac{\delta_D(x)}{r}\right)^p$.
Finally, from \eqref{e:qDf} we conclude that 
\begin{theorem}\label{t:DHKE1} 
Suppose that $D$ is a $C^{1,1}$ open set in $\R^d$, $d\ge 2$, with characteristics $(R_2$, $\Lambda)$. 
For all $T>0$,  $p \in [\alpha-1, \alpha)\cap (0, \alpha)$ 
and $\eta \in [0, \alpha)$, 
there exists a constant $c=c(C_1, C_2, p, \alpha, d, \eta, T, R_2, \Lambda) \ge 1$ such that for all $\kappa \in \sH_\alpha (p)$, 
  the transition density $q^D(t,x,y)$ of the 
 Hunt process $Y$ on $D$ corresponding to the Feynman-Kac semigroup of $X^D$  
via the multiplicative functional $e^{-\int_0^t \kappa (X^D_s)ds}$ 
 satisfies that 
\begin{align*}
&c^{-1} \left(1 \wedge \frac{\delta_D(x)}{t^{1/\alpha}}\right)^p \left(1 \wedge \frac{\delta_D(y)}{t^{1/\alpha}}\right)^p \left[t^{-d/\alpha} \wedge \frac{t}{|x-y|^{d+\alpha}}\right] \\
\le& 
q^D(t,x,y) \le c \left(1 \wedge \frac{\delta_D(x)}{t^{1/\alpha}}\right)^p \left(1 \wedge \frac{\delta_D(y)}{t^{1/\alpha}}\right)^p \left[t^{-d/\alpha} \wedge \frac{t}{|x-y|^{d+\alpha}}\right]
\end{align*}
for $(t,x,y) \in (0,T] \times D \times D$. 
\end{theorem}

In the case $D=\R^d_+$ and $\kappa(x)=C(d,\alpha,p)x^{-\alpha}_d$, 
one can use the scaling property to
get that the two-sided heat estimates in Theorem \ref{t:DHKE1}  is valid for all $t>0$.

\smallskip
\begin{remark}
Theorem \ref{t:DHKE1}  also holds in $d=1$. In fact, let $D \subset \R$ be a union of open intervals with a localization radius $r_0$ and $C(1,\alpha,p) = \sA(1,-\alpha)\gamma(\alpha,p)$. 
The first difference of the proof appears in the bound of $|III|$ in Lemma \ref{l:3.1}. 
We use the following calculation instead of \cite[Lemma 5.6]{BBC03}:
\begin{align*}
|\kappa_D(x)-\kappa_{\R_+}(x)| &\le |\kappa_{(-r_0,0)^c}(x)-\kappa_{(0,r_0)}(x)| = \sA(1,-\alpha)\left[\int_{-\infty}^{-r_0}+\int_{r_0}^{\infty} \frac{dy}{|y-x|^{1+\alpha}}\right]\\
&=\frac{\sA(1,-\alpha)}{\alpha}((r_0-x)^{-\alpha}+(r_0+x)^{-\alpha}) \le c,
\end{align*}
provided $\delta_D(x)<r_0/2$. 
Moreover, the bound for $|I_{\epsilon}|$ is easy in Lemma \ref{l:3.1}: Since $h_q(y)=w_q(y)$ for $y \in (-\infty, r_0)$,
$
I_{\epsilon} \le c \int_{r_0}^{\infty} y^{q-1-\alpha} dy = c.
$
\end{remark}

\begin{remark}\label{r:new}
It follows from \cite[pp.94--95]{BBC03} that $Z^D$ can  be obtained from 
$X^D$
via a Feynman-Kac perturbation of the form $e^{-\int^t_0\kappa_D(X^D_s)ds}$.
In view of \eqref{e:new1},
$\kappa_D$ satisfies condition \eqref{e:C1}
with $C_1= \frac{\sA(d,-\alpha)}{\alpha} \frac{\omega_{d-1}}{2} \beta(\frac{\alpha+1}{2}, \frac{d-1}{2})$. 
By direct calculation, we can see that $\gamma(\alpha,\alpha/2)=1/\alpha$. 
This means that $C_1=C(d,\alpha,\alpha/2)$. 
Thus Theorem 
\ref{t:DHKE1} recovers \eqref{e:est_p_Dc}.
 When $\alpha\in (1,2)$, 
$C_1=0=C(d,\alpha,\alpha-1) $ is allowed. 
 Thus, by taking $\kappa=0$, Theorem 
\ref{t:DHKE1} recovers  \eqref{e:est_p_Dc1} as well.

We also remark here that  Theorem 
\ref{t:DHKE1} provides examples of processes studied in \cite{CKS15} (see \eqref{e:C1} and \cite[Proposition 4.1(ii)]{CKS15}).
\end{remark}

\subsection{Non-local perturbation in bounded $C^{1,1}$ open set}\label{ss-nl}

Recall that
$ {\mathcal A}(d, \alpha)= \alpha2^{\alpha-1}\pi^{-d/2}
\Gamma(\frac{d+\alpha}2) \Gamma(1-\frac{\alpha}2)^{-1}$. 
We also recall that   we write $y=(\wt{y},y_d)$ for $y\in \R^d$. 
  For $u:\R^d_+\to [0,\infty)$, $\lambda \in (0, \infty)$ and $\beta\in (-\infty, 2)$, we define
$$
L_{d, \lambda}^{\beta} u(x) 
:= \lim_{\eps \downarrow
0}\int_{\{y\in \bR_+^d: \, \eps <|y-x|<\lambda \}} (u(y)-u(x))
\frac{dy}{|x-y|^{d+\beta}} \, , \qquad x\in \R^d_+\, .
$$
In the remainder of this subsection, we will assume that $\beta\in (-\infty, 2)$ .
For any real number $p$ and $y \in \R^d_+$, let $g(y):=y_d^p=\delta_{\R^d_+}(y)^p$. 

\begin{lemma}\label{l:4.1}
For all positive $p, \lambda$  and $\beta\in (-\infty, 2)$, there exist $c_1=c_1(p,d,\beta, \lambda)>0$ 
and $c_2=c_2(p,d,\beta, \lambda)\in (0,\frac{1}{4})$ such that, for every $x \in \R^d_+$ with $0<x_d<c_2$, the following inequalities hold:
$$
| L_{d, \lambda}^{\beta} g(x)| \le c_1
\begin{cases}
1 & \text{ if } p>\beta;\\
 |\log x_d|& \text{ if } p=\beta;\\
 x_d^{p-\beta}& \text{ if } p<\beta.
\end{cases}
$$
\end{lemma}

\pf
When $\beta \le 0$, then clearly for $x \in \R^d_+$,
\begin{align*}
&  \int_{\R^d_+}\frac{|y_d^p-x_d^p|}{|y-x|^{d+\beta}}{\bf 1}_{\{|y-x|<\lambda \}} \, dy
 \le  c\int_{B(0, \lambda)}|z|^{-d-\beta+p}  dz \le c \int_0^{\lambda} s^{-1-\beta+p} ds =c \lambda^{p-\beta}.
\end{align*}

We now assume $\beta> 0$. 
For simplicity, take $x=(\wt{0}, x_d)$ and denote $\mathbf{e}_d=(\tilde{0}, 1)$.
Then by the change of variables $z=y/x_d$, we have
\begin{align*}
L_{d, \lambda}^{\beta} g(x)=
& \textrm{p.v.} \int_{\R^d_+}\frac{y_d^p-x_d^p}{|y-x|^{d+\beta}}{\bf 1}_{\{|y-x|<\lambda \}} \, dy \nn\\
=& x_d^{p-\beta} \textrm{p.v.} \int_{\R^{d-1}}\int_0^{\infty}\frac{z_d^p-1}{|z-\mathbf{e}_d|^{d+\beta}}{\bf 1}_{\{|z-\mathbf{e}_d|<\lambda/x_d \}}   dz_d d \wt z 
=:x_d^{p-\beta} I_1 \, \nn .
\end{align*}
Using the change of variables $\wt{z}=|z_d-1|\wt{u}$, we get
\begin{align*}
I_1
=&\int_{\R^{d-1}}\frac{1}{(|\wt{u}|^2+1)^{(d+\beta)/2}} \left(\textrm{p.v.} \int_0^{\infty}\frac{z_d^p-1}{|z_d-1|^{1+\beta}} {\bf 1}_{\{|z_d-1|<\lambda(|\wt{u}|^2+1)^{-1/2}/x_d \}}  dz_d   \right) d\wt{u}\\
=:& \int_{\R^{d-1}}\frac{1}{(|\wt{u}|^2+1)^{(d+\beta)/2}} I_2 d\wt{u}\, .
\end{align*}
Fix $\wt{u}$ and let $M:=(|\wt{u}|^2+1)^{1/2}$. Then 
\begin{align}
I_2=\lim_{\epsilon\to 0}\left(\int_{(1-\frac{\lambda}{Mx_d})_+}^{1-\epsilon}\frac{z_d^p-1}{|z_d-1|^{1+\beta}} dz_d+\int_{1+\epsilon}^{\frac{Mx_d+\lambda}{Mx_d}}\frac{z_d^p-1}{|z_d-1|^{1+\beta}} dz_d \right). \label{e:MM}
\end{align}
By using the change of variables $w=1/z_d$, we get that, for $\epsilon<\lambda/(Mx_d)$, the second integral in \eqref{e:MM} is equal to
\begin{eqnarray*}
\int_{\frac{Mx_d}{Mx_d+\lambda}}^{\frac{1}{1+\epsilon}}\frac{w^{\beta-1-p}-w^{\beta-1}}{(1-w)^{1+\beta}} dw
=\int_{\frac{Mx_d}{Mx_d+\lambda}}^{1-\epsilon}\frac{w^{\beta-1-p}-w^{\beta-1}}{(1-w)^{1+\beta}}\, dw+\int_{1-\epsilon}^{\frac{1}{1+\epsilon}}\frac{w^{\beta-1-p}-w^{\beta-1}}{(1-w)^{1+\beta}}\, dw\, .
\end{eqnarray*}
Note that from \cite[p.121]{BBC03}, we see that 
$$
\left|\int_{1-\epsilon}^{\frac{1}{1+\epsilon}}\frac{w^{\beta-1-p}-w^{\beta-1}}{(1-w)^{1+\beta}}\, dw\right|\le c \epsilon^{2-\beta}\, .
$$
 By writing the first integral in \eqref{e:MM} as
$$
\int_{\frac{Mx_d}{Mx_d+\lambda}}^{1-\epsilon}\frac{w^p-1}{(1-w)^{1+\beta}}dw - \int_{(1-\frac{\lambda}{Mx_d})_+}^\frac{Mx_d}{Mx_d+\lambda} \frac{1-w^p}{(1-w)^{1+\beta}} dw,
$$
and by  using 
\begin{equation}\label{e:wwww}
 (w^p-1)+(w^{\beta-1-p}-w^{\beta-1})=(1-w^p)(1-w^{p-(\beta-1)})w^{\beta-1-p},
 \end{equation}
we have 
\begin{align}\label{e:I2}
I_2=\lim_{\epsilon\to 0}\int_{\frac{Mx_d}{Mx_d+\lambda}}^{1-\epsilon} \frac{(1-w^p)(1-w^{p-(\beta-1)})}{(1-w)^{1+\beta}}w^{\beta-1-p}\, dw -\int^{\frac{Mx_d}{Mx_d+\lambda}}_{(1-\frac{\lambda}{Mx_d})_+}\frac{1-w^p}{(1-w)^{1+\beta}} dw =: I_{21}-I_{22}\, .
\end{align}
First,  it is easy to see that
$$
0<I_{22}\le \int^{\frac{Mx_d}{Mx_d+\lambda}}_{0}\frac{1-w^p}{(1-w)^{1+\beta}} dw\le c  
\begin{cases}
1 & \text{if } \beta \in (0,1);\\
\log (1+ Mx_d/\lambda)& \text{if } \beta =1;\\
(1+ Mx_d/\lambda)^{\beta-1} & \text{if } \beta \in (1,2).
\end{cases}
$$
Next, since $\beta<2$, the fraction in $I_{21}$  is integrable near 1. Thus, 
$$
I_{21}=\int_{\frac{Mx_d}{Mx_d+\lambda}}^1  \frac{(1-w^p)(1-w^{p-(\beta-1)})}{(1-w)^{1+\beta}}w^{\beta-1-p}\, dw.
$$
Note that, if $\frac{Mx_d}{Mx_d+\lambda} \ge 1/4$, then clearly,
$
I_{21}  \le c <\infty.
$
If $\frac{Mx_d}{Mx_d+\lambda} < 1/4$, then 
\begin{align*}
I_{21} \le& c+
\int_{\frac{Mx_d}{Mx_d+\lambda}}^{1/2}  \frac{(1-w^p)(1-w^{p-(\beta-1)})}{(1-w)^{1+\beta}}w^{\beta-1-p}\, dw\le c+
c\int_{\frac{Mx_d}{Mx_d+\lambda}}^1w^{\beta-1-p}\, dw.
\end{align*}
Thus
$$
I_{21} \le c \begin{cases}
(1+ \lambda/(Mx_d))^{p-\beta} & \text{if } p>\beta ;\\
  \log (1+ \lambda/(Mx_d))& \text{if } p=\beta.
\end{cases}
$$
Therefore, if $ p>\beta$,  then for small $x_d$,
\begin{align*}
&|x_d^{p-\beta} I_1|\\
\le&c x_d^{p-\beta}
\int_{\R^{d-1}}\frac{1}{(|\wt{u}|^2+1)^{(d+\beta)/2}}\\
&
\times\begin{cases}
(1+  \lambda/((|\wt{u}|^2+1)^{1/2}x_d))^{p-\beta}  d\wt{u} & \text{if } 
\beta \in (0,1);\\
((1+  \lambda/((|\wt{u}|^2+1)^{1/2}x_d))^{p-\beta}  + \log(1+ (|\wt{u}|^2+1)^{1/2}x_d/\lambda))d\wt{u}& \text{if } \beta =1;\\
((1+  \lambda/((|\wt{u}|^2+1)^{1/2}x_d))^{p-\beta} +(1+ (|\wt{u}|^2+1)^{1/2}x_d/\lambda)^{\beta-1})d\wt{u} & \text{if } \beta \in (1,2)
\end{cases}\\
\le&c 
\int_{\R^{d-1}}\frac{1}{(|\wt{u}|^2+1)^{(d+\beta)/2}}\\
&
\times\begin{cases}
(x_d+  \lambda/(|\wt{u}|^2+1)^{1/2})^{p-\beta}  d\wt{u} & \text{if } 
\beta \in (0,1);\\
((x_d+  \lambda/(|\wt{u}|^2+1)^{1/2})^{p-\beta}  +\log(1+ (|\wt{u}|^2+1)^{1/2}/\lambda))d\wt{u}& \text{if } \beta =1;\\
((x_d+  \lambda/(|\wt{u}|^2+1)^{1/2})^{p-\beta} +(1+ (|\wt{u}|^2+1)^{1/2}/\lambda)^{\beta-1})d\wt{u} & \text{if } \beta \in (1,2)
\end{cases}\\
\le&c (\lambda)
\begin{cases}
\int_{\R^{d-1}}{(|\wt{u}|^2+1)^{-(d+\beta)/2}} d\wt{u} & \text{if } 
\beta \in (0,1);\\
\int_{\R^{d-1}}\frac{\log(1+ (|\wt{u}|^2+1)^{1/2})}{(|\wt{u}|^2+1)^{(d+\beta)/2}}d\wt{u}& \text{if } \beta =1;\\
\int_{\R^{d-1}}{(|\wt{u}|^2+1)^{-(d+1)/2}}{d\wt{u}} & \text{if } \beta \in (1,2)
\end{cases}=c(\lambda, \beta) < \infty.
\end{align*}
If $ p=\beta>0$, then for small $x_d$,
\begin{align*}
&| I_1|
\le c 
\int_{\R^{d-1}}\frac{1}{(|\wt{u}|^2+1)^{(d+\beta)/2}}
\times\\
&
\begin{cases}
\log(1+ \lambda/((|\wt{u}|^2+1)^{1/2}x_d)) d\wt{u} & \text{if } \beta \in (0,1);\\
[\log(1+ \lambda/((|\wt{u}|^2+1)^{1/2}x_d)) + \log(1+ (|\wt{u}|^2+1)^{1/2}x_d/\lambda)]d\wt{u}& \text{if } \beta =1;\\
[\log(1+ \lambda/(((|\wt{u}|^2+1)^{1/2}x_d))+(1+ (|\wt{u}|^2+1)^{1/2}x_d/\lambda)^{\beta-1})]d\wt{u} & \text{if } \beta \in (1,2)
\end{cases}\\
\le&c 
\int_{\R^{d-1}}\frac{1}{(|\wt{u}|^2+1)^{(d+\beta)/2}}
\begin{cases}
\log(1+ \lambda/x_d)   d\wt{u} & \text{if } \beta \in (0,1);\\
[\log(1+ \lambda/x_d)   +\log(1+ (|\wt{u}|^2+1)^{1/2}/\lambda)]d\wt{u}& \text{if } \beta =1;\\
[\log(1+ \lambda/x_d)   +(1+ (|\wt{u}|^2+1)^{1/2}/\lambda)^{\beta-1})]d\wt{u} & \text{if } \beta \in (1,2)
\end{cases}\\
\le& c  \log(1+ \lambda/x_d).
\end{align*}

We now assume that $0<p<\beta$. Note that  by \eqref{e:I2}, \eqref{e:wwww} and simple algebra, 
\begin{eqnarray*}
I_2=\int_{(1-\frac{\lambda}{Mx_d})_+}^1  \frac{(1-w^p)(1-w^{p-(\beta-1)})}{(1-w)^{1+\beta}}w^{\beta-1-p}\, dw-\int_{(1-\frac{\lambda}{Mx_d})_+}^{\frac{Mx_d}{Mx_d+\lambda}}\frac{w^{\beta-1-p}-w^{\beta-1}}{(1-w)^{1+\beta}}\, dw.
\end{eqnarray*}
Since $w\mapsto w^{\beta-1-p}$ is integrable near 0, 
\begin{eqnarray*}
x_d^{p-\beta} I_1 \le 
x_d^{p-\beta}\int_{\R^{d-1}}\frac{1}{(|\wt{u}|^2+1)^{(d+\beta)/2}}\left( \int_{0}^1  \frac{(1-w^p)(1-w^{p-(\beta-1)})}{(1-w)^{1+\beta}}w^{\beta-1-p}\, dw\right)d\wt{u} .
\end{eqnarray*}
On the other hand, $-x_d^{p-\beta} I_1 \le c(d) x_d^{p-\beta} I_{1,2}$, where 
\begin{align*}
I_{1,2}:=\int_{0}^{\infty} \frac{1}{(u^2+1)^{(d+\beta)/2}}
 \int_{(1-\frac{\lambda}{(u^2+1)^{1/2}x_d})_+}^{\frac{(u^2+1)^{1/2}x_d}{(u^2+1)^{1/2}x_d+\lambda}}\frac{w^{\beta-1-p}(1-w^{p})}{(1-w)^{1+\beta}}\, dw u^{d-2}du.
\end{align*}
Note that 
\begin{align*}
&\sup_{v \ge 2\lambda} \int_{1-\lambda/v}^{v/(v+\lambda)} \frac{w^{\beta-1-p}(1-w^{p})}{(1-w)^{1+\beta}}\, dw
\le c\sup_{v \ge 2\lambda} \int_{1-\lambda/v}^{v/(v+\lambda)}(1-w)^{-\beta}dw\\
&= c\sup_{v \ge 2\lambda} \int_{\lambda/(v+\lambda)}^{\lambda/v}t^{-\beta}dt\
\le c\sup_{v \ge 2\lambda} (v+\lambda)^{\beta} (\frac{1}{v}-\frac{1}{v+\lambda} )\le c \sup_{v \ge 2\lambda} v^{\beta-2} <\infty,
\end{align*}
and, for $x_d < \lambda$,  
\begin{align*}
&\sup_{x_d \le v < 2\lambda} \int_{0}^{v/(v+\lambda)} \frac{w^{\beta-1-p}(1-w^{p})}{(1-w)^{1+\beta}}\, dw
\le c\sup_{x_d \le v < 2\lambda} \int_{0}^{2/3}w^{\beta-1-p}dw <\infty.
\end{align*}
Thus for $x_d < \lambda$, 
\begin{align*}
0<I_{1,2} \le& 
\int_{0}^{\infty} \frac{{\bf 1}_{(u^2+1)^{1/2}x_d< 2\lambda}}{(u^2+1)^{(d+\beta)/2}}
 \int_{0}^{\frac{(u^2+1)^{1/2}x_d}{(u^2+1)^{1/2}x_d+\lambda}}\frac{w^{\beta-1-p}(1-w^{p})}{(1-w)^{1+\beta}}\, dw u^{d-2}du\\
 &+\int_{0}^{\infty} \frac{{\bf 1}_{(u^2+1)^{1/2}x_d \ge 2\lambda}}{(u^2+1)^{(d+\beta)/2}}
 \int_{(1-\frac{\lambda}{(u^2+1)^{1/2}x_d})_+}^{\frac{(u^2+1)^{1/2}x_d}{(u^2+1)^{1/2}x_d+\lambda}}\frac{w^{\beta-1-p}(1-w^{p})}{(1-w)^{1+\beta}}\, dw u^{d-2}du\\
 \le &
\int_{0}^{\infty} \frac{u^{d-2}}{(u^2+1)^{(d+\beta)/2}}
du <\infty.
\end{align*}
\qed

Throughout the remainder of this subsection we assume that $D$  is a bounded 
$C^{1,1}$ open subset of $\R^d$ and 
$\alpha \in (0 \vee \beta, 2)$.
We also assume that 
$b(x,y)$ is a symmetric Borel function on $D \times D$ such that
$C_{b,1}:=\sup_{x,y \in D} |b(x,y)|<\infty$ and
the function 
$$B(x,y) :={\cal A} (d, \, -\alpha) +|x-y|^{\alpha-\beta}b(x,y), \quad x, y\in D, $$
is bounded below by a positive constant, that is,
$C_{b, 2}\le  B(x,y)$ for some $C_{b, 2} \in (0, \infty)$. Clearly, $B(x,y)$ is bounded above by ${\cal A} (d, \, -\alpha) +($diam$(D))^{\alpha-\beta}C_{b,1}$.
 
We further assume that the first partials of $B(x, y)$ are bounded on $D\times D$.
Note that, $\beta$ and $b$ can be negative, as long as the condition above is satisfied.
Let $({\EE}^{(B)}, \overline{\FF})$ be the
Dirichlet form on $L^2(D, dx)$ defined by
\begin{eqnarray*}
{\EE}^{(B)} (u, v):=  \frac{1}{2}
\int_D \int_D{(u(x)-u(y))(v(x)-v(y))}
\frac{B(x,y)}{|x-y|^{d+\alpha}} dxdy, \quad   u, v \in \overline{\FF}.
\end{eqnarray*}
By \cite{CK03}, $({\EE}^{(B)},\, \overline{\FF}\,)$ is a regular
Dirichlet form on $\overline D$ and
its associated Hunt process $X^{(B)}$ is conservative and lives on $\overline D$. 
Moreover, since $B(x,y)$ is bounded on $D\times D$ between two strictly positive constants, the form $({\EE}^{(B)},\, \overline{\FF}\,)$ satisfies the assumptions of \cite[Remark 2.4]{BBC03}, so we can freely use results of \cite[Section 2]{BBC03}.
Further, $X^{(B)}$  
admits a strictly positive and
jointly continuous transition density $p(t,x,y)$ with respect to the Lebesgue measure $dx$ such that
\begin{align}\label{e:est_p_Dr2}
C_0^{-1} \left[t^{-d/\alpha} \wedge \frac{t}{|x-y|^{d+\alpha}}\right] \le p(t,x,y) \le C_0 \left[t^{-d/\alpha} \wedge \frac{t}{|x-y|^{d+\alpha}}\right],
\end{align}
for $(t,x,y) \in (0,1)\times D\times D$.
 
Let $L^{(B)}$ be the generator of $X^{(B)}$ in the $L^2$ sense. Similar to \cite[Section 4]{SV}, cf.~also \cite{Kom88}, we can show that $C^2_c(D)$ is contained in the domain of $L^{(B)}$ and give an explicit expression for $L^{(B)}f$ when $f\in C^2_c(D)$. Using these, one can check that the 
process $X^{(B)}$ satisfies Assumptions {\bf A} and {\bf U}.

If $m>0$, by taking $\beta=\alpha-2$ and
$b(x,y)={\cal A} (d, \, -\alpha)(\varphi(m^{1/\alpha}|x-y|)-1) |x-y|^{-2}$ with
$$
\varphi (r):= 2^{-(d+\alpha)} \, \Gamma \left(
\frac{d+\alpha}{2} \right)^{-1}\, \int_0^\infty s^{\frac{d+\alpha}{
2}-1} e^{-\frac{s}{ 4} -\frac{r^2}{ s} } \, ds,
$$
we cover the reflected relativistic $\alpha$-stable process $X^m$ with weight $m>0$ in $\overline{D}$. When $\alpha\in (1, 2)$, the killed process $X^{m, D}$ is the censored relativistic $\alpha$-stable process in $D$. When $\alpha\in (0, 1]$, it follows from \cite[Section 2]{BBC03} that, starting from inside $D$, the process $X^m$ neither hits nor approaches $\partial D$ at any finite time. 
Thus, the killed process $X^{m, D}$ is simply $X^m$ restricted to $D$. 

Recall that 
for $u:D\to [0,\infty)$, 
$$
L_{\beta} u(x) = \sA(d,-\beta)\lim_{\eps \downarrow
0}\int_{\{y\in D: \, \eps <|y-x| \}} (u(y)-u(x))
\frac{dy}{|x-y|^{d+\beta}} \, , \qquad x\in D\, .
$$

Let 
$$
 L_{\beta, b} u(x) := 
 \lim_{\eps \downarrow
0}\int_{\{y\in D: \, \eps <|y-x|\}} (u(y)-u(x))
\frac{b(x,y)}{|x-y|^{d+\beta}}dy \, , \qquad x\in D\, .
$$
Let $p\in [\alpha-1, \alpha)\cap (0, \alpha)$, $\kappa\in \sH_\alpha(p)$. If $\beta \ge p$,
then  we always assume that, there exist $C_{b,3}>0$ and 
$\beta_1>\beta-p$ such that 
\begin{align}
\label{e:bholder}
|b(x,y)-b(x,x)| \le  C_{b,3} |x-y|^{\beta_1}, \quad x, y \in D.
\end{align}
Note that, under \eqref{e:bholder}, for any bounded Borel function $u$ satisfying $|u(x)-u(y)| \le c |x-y|^p$ on $D$, 
\begin{align}
\label{e:bhL}
| L_{\beta, b} u(x)|& \le
\int_{\{y\in D: \, \eps <|y-x| \}} |u(y)-u(x)|
\frac{|b(x,y)-b(x,x)|}{|x-y|^{d+\beta}}dy+\frac{|b(x,x)|}{\sA(d,-\beta)}
 |L_{\beta} u(x)|\nn\\
& \le  c_1+ c_2  |L_{\beta} u(x)|.
\end{align}

Recall that for an open set $D$ and  $q\ge 0$, $h_q(x)=\delta_D(x)^q$.

\begin{lemma}\label{l:4.2}
Let $D$ be a bounded $C^{1,1}$ open set with characteristics $(R_2, \Lambda)$.
For any $q \ge p$, there exist constants $c_1>0$ and $c_2\in (0,(R_2 \wedge 1)/{4})$
depending only on $p,q,d,\beta, R_2, \Lambda$, $\text{diam}(D)$, $ C_{b,1}, C_{b,2}, C_{b,3}, \beta_1$ 
such that for every $x \in D$ with $0<\delta_D(x)<c_2$, the following inequalities hold:
$$
 |  L_{\beta, b} h_q(x)| \le c_1
\begin{cases}
1 & \text{ if } q>\beta;\\
 |\log\delta_D(x)|& \text{ if } q=\beta;\\
\delta_D(x)^{q-\beta}& \text{ if } q<\beta.
\end{cases}
$$
\end{lemma}

\pf
Without loss of generality we assume  diam$(D) \le 1$ and  
let $x \in D$ with $\delta_D(x)<R_2/4$. 
Choose a point $z \in \partial D$ such that $\delta_D(x) = |x-z|$. Then, there exists a $C^{1,1}$ function $\Gamma : \R^{d-1} \mapsto \R$ such that $\Gamma(z)=\nabla \Gamma(z)=0$ and an orthonormal coordinate system $CS_z$ with origin at $z$ such that
$$
D \cap B(z,R_2) = \{ y = (\ty,y_d) \mbox{ in } CS_z: 
y_d > \Gamma(\ty) \} \cap B(z,R_2),
$$
and  $z=0$ and $x=(\tx,x_d)=(0,x_d)$ in $CS_z$.
Define $w_q(y) := (y_d)^q$ for $y \in \R^d_+$ and $w_q(y):=0$, otherwise.
For any open subset $U \subset \R^d$, 
define $\wh \kappa_U(x) := \sA(d,-\beta)\int_{U^c\cap B(x,1)} {|y-x|^{-d-\beta}}{dy}. $
 Since $h_q(x)=w_q(x) = x_d^q$, using \eqref{e:bhL} we have 
\begin{align*} 
&L_{\beta} h_q(x)  = \sA(d,-\beta)\; \lim_{\epsilon \downarrow 0} \left[ \int_{1 >|y-x|>\epsilon} \frac{h_q(y)-h_q(x)}{|y-x|^{d+\beta}} dy + \wh \kappa_D(x)h_q(x)\right] \\
& = \sA(d,-\beta)\; \lim_{\epsilon \downarrow 0} \left[ \int_{1 >|y-x|>\epsilon} \frac{h_q(y)-w_q(y)}{|y-x|^{d+\beta}} dy + \int_{1> |y-x|>\epsilon} \frac{w_q(y)-w_q(x)}{|y-x|^{d+\beta}} dy +\wh \kappa_D(x)w_q(x)\right]  \\
& =  L_{d, 1}^{\beta}w_q(x) 
+ \sA(d,-\beta)\; \lim_{\epsilon \downarrow 0} \left[ \int_{1 > |y-x|>\epsilon} \frac{h_q(y)-w_q(y)}{|y-x|^{d+\beta}} dy + (\wh \kappa_D(x)-\wh \kappa_{\R^d_+}(x))w_q(x)\right].
\end{align*}
By a similar calculation as in \cite[Lemma 5.6]{BBC03}, for any $0<\beta<2$, we get
$$
|(\wh \kappa_D(x)-\wh \kappa_{\R^d_+}(x))w_q(x)| \le c x_d^q (x_d^{1-\beta}+|\log x_d|) \le c.
$$
Therefore, it remains to bound $I_{\epsilon}:=\lim_{\epsilon \downarrow 0} \int_{1>|y-x|>\epsilon} \frac{h_q(y)-w_q(y)}{|y-x|^{d+\beta}} dy$.
When $q <\beta$, by the proof of Lemma \ref{l:3.1},
$$\sup_{\eps<1/2}
 |I_{\epsilon}|  \le
\sup_{\eps<1/2}
\left|\int_{\R^d, |y-x|>\epsilon} \frac{h_q(y)-w_q(y)}{|y-x|^{d+\beta}} dy \right| \le c(1+\log|x_d|).
$$
When $q \ge \beta$, by \cite[(3.13)]{CKSV12}, we get 
$\sup_{\eps<1/2}
 |I_{\epsilon}|  \le c$.
The lemma now follows from these bounds, Lemma \ref{l:4.1} and \eqref{e:bhL}.
\qed

Let $p\in [\alpha-1, \alpha)\cap (0, \alpha)$, $\kappa\in \sH_\alpha(p)$, 
 and define 
\begin{equation*}
\wt Lf(x) := L_\alpha f(x)+L_{\beta, b} f(x)-\kappa(x)f(x)=Lf(x)+L_{\beta,b} f(x).
\end{equation*}

Combining Lemmas \ref{l:3.1} and \ref{l:4.2}, we get the following lemma.
\begin{lemma}\label{l:4.3}
Let $0<p\le q < \alpha$ and $ \beta< \alpha$ and  define
$$h_q(x):=\delta_D(x)^q.$$
Then there exist 
$c_1>0$ and  $c_2\in (0,(R_2\wedge 1)/{4})$ depending only on $R_2, p,q,d$, $\alpha, \beta, \Lambda$, $C_2,\eta, C_{b,1}$, $C_{b,2}, C_{b,3}, \beta_1$  such that the following inequalities hold:

(i) If $q>p$,
$$
c_1^{-1} \delta_D(x)^{q-\alpha} \le 
\wt Lh_q(x) \le c_1 \delta_D(x)^{q-\alpha}
$$
for every $x \in D$ with $0<\delta_D(x)<c_2$.

(ii) If $q=p$,
$$
|\wt Lh_p(x)| \le c_1 
( \delta_D(x)^{p-(\beta \vee \eta)}  + |\log \delta_D(x)|)
$$
for every $x \in D$ with $0<\delta_D(x)<c_2$.
\end{lemma}

Recall that  $X^{(B), D}$ denote the process $X^{(B)}$ killed upon exiting $D$. 
Note that
the operator $\wt L$ coincides with the restriction to $C_c^2(D)$ of the generator of 
of the Feynman-Kac semigroup of 
  $X^{(B), D}$ 
  via the multiplicative functional 
  $e^{-\int_0^t \kappa (X^{(B), D}_s)ds}$ in $C_0(D)$.
We now follow the argument of the previous subsection 
(choosing  $q\in (p,(p-(\eta \vee \beta)+\alpha) \wedge \alpha)$) 
and can conclude the following.

\begin{theorem}\label{t:DHKE2} 
Suppose that $D$ is a bounded $C^{1,1}$ open set in $\R^d$, $d\ge 2$, with characteristics $(R_2$, $\Lambda)$. For all $T>0$,  $p \in [\alpha-1, \alpha)\cap (0, \alpha)$, $\beta<\alpha$ and $\eta \in [0, \alpha)$, 
there exists $c=c(C_1, C_2, p, \alpha, \beta, d, \eta, b, \text{diam}(D), T, C_{b,1}, C_{b,2}, C_{b,3}, \beta_1) \ge 1$ such that for all $\kappa \in \sH_\alpha (p)$, 
  the transition density $q^D(t,x,y)$ of the
  Hunt process $Y$ on $D$ corresponding to the Feynman-Kac semigroup of 
  $X^{(B), D}$ 
  via the multiplicative functional 
  $e^{-\int_0^t \kappa (X^{(B), D}_s)ds}$ satisfies that 
\begin{align*}
&c^{-1} \left(1 \wedge \frac{\delta_D(x)}{t^{1/\alpha}}\right)^p \left(1 \wedge \frac{\delta_D(y)}{t^{1/\alpha}}\right)^p \left[t^{-d/\alpha} \wedge \frac{t}{|x-y|^{d+\alpha}}\right] \\
\le& 
q^D(t,x,y) \le c \left(1 \wedge \frac{\delta_D(x)}{t^{1/\alpha}}\right)^p \left(1 \wedge \frac{\delta_D(y)}{t^{1/\alpha}}\right)^p \left[t^{-d/\alpha} \wedge \frac{t}{|x-y|^{d+\alpha}}\right]
\end{align*}
for $(t,x,y) \in (0,T] \times D \times D$. 
\end{theorem}

We remark here that Theorem \ref{t:DHKE2} recovers \cite[Theorem 4.8]{CKS15}.
Let $\kappa_D^m$ be the killing function of the killed relativistic
$\alpha$-stable process $Z^{m, D}$ in $D$.
It follows from \cite[pp.94--95]{BBC03} that the killed relativistic
$\alpha$-stable process $Z^{m, D}$ can  be obtained from 
$X^{m, D}$
via a Feynman-Kac perturbation of the form $e^{-\int^t_0\kappa_D^m(X^{m, D}_s)ds}$.
  It follows \cite[p. 278]{CS-JFA} that
$0\le \kappa_D(x)-\kappa_D^m(x)\le c\delta_D(x)^{2-\alpha}$ for all $x\in D$.
Combining this with \eqref{e:new1}, we get
$$
|(\kappa_D^m(x)-\kappa_{\R^d_+}(x))w_q(x)| \le c x_d^q (x_d^{1-\alpha}+|\log x_d|) \le c.
$$
Now by the same argument as in Remark \ref{r:new}, we see that Theorem \ref{t:DHKE2}
recovers the main result of \cite{CKS-AOP} for bounded $C^{1, 1}$ open set $D$.

\subsection{$\R^d\setminus\{0\}$}\label{ss-eorigin}
In this subsection we assume that $\X=\R^d$, $d\ge 2$,  $X$ is an isotropic $\alpha$-stable process
on $\R^d$ and $D=\R^d\setminus\{0\}$.
Obviously, $D$ is a $(1/2)$-fat open set with characteristics $(\infty, 1/2)$ and  $X$ satisfies Assumptions {\bf A} and
{\bf U}. Since $X$ does not hit $\{0\}$, 
the killed process $X^D$ is simply the restriction of $X$ to $D$.

Recall that
$\sA(d,-\alpha)=\alpha 2^{\alpha-1} \pi^{-d/2} \Gamma((d+\alpha)/2)\Gamma(1-\alpha/2)^{-1}$.
Let $p\in (0, \alpha)$ and define
$$
H(s)=2\pi\frac{\pi^{\frac{d-3}{2}}} { \Gamma ( \frac{d-1}{2} )}\int_{0}^{\pi}\sin^{d-2}\!\theta \, \frac{ (  \sqrt{s^2-\sin^2\theta}+\cos\theta)^{1+{\alpha}}}{\sqrt{s^2-\sin^2 \theta}} \, d\theta, \quad s \ge 1, \quad 
$$
and
$$\wt C(\alpha, d, p):=\sA(d,-\alpha) \, \int_1^{+\infty}
(s^p-1)(1-s^{-d+\alpha-p})
s( s^2-1
  )^{-1-{\alpha}} H(s)ds.$$
  Note that $p \to \wt{C}(\alpha, d, p)$ is strictly increasing on $(0, \alpha)$.
The function  $H(s)$ is positive and continuous on $[1, +\infty)$ with  
$H(s)\asymp s^{{\alpha}}$ for large $s$ and 
$$
s( s^2-1)^{-1-{\alpha}} H(s)  \asymp (s-1)^{-1-{\alpha}},  \quad s \ge 1,
$$
(see the paragraph after \cite[Theorem 1.1]{FV12}).
Thus
\begin{align}
\label{e:palphab}
\lim_{p \downarrow 0}\wt C(\alpha, d, p)=0 \quad \text{and}\quad \lim_{p \uparrow \alpha}\wt C(\alpha, d, p)=\infty.
\end{align}
Applying  \cite[Theorem 1.1]{FV12} to $u_p:=|x|^p$, we get that
 \begin{equation}\begin{split}\label{rapresth}
-(-\Delta)^{\alpha/2}u_p(x)= \wt C(\alpha, d, p) \, |x|^{p-\alpha}, \quad |x|>0, x \in \mathbb{R}^d.
\end{split}
\end{equation}

Let $\sG_\alpha$ be the collection of  non-negative functions on 
$D$ such that for each $\kappa  \in \sG_\alpha$
there exist constants $C_1>0$, $C_2 \ge 0 $ 
and $\eta \in [0, \alpha)$ 
such that $\kappa(x)\le C_2$ for all $x$ with $|x|\ge 1$ and 
\begin{equation}\label{e:kappainf}
\big|\kappa(x)-C_1 |x|^{-\alpha}\big| \le C_2 |x|^{-\eta},
\end{equation} 
for all $x \in D$ with $|x|<1$.
By \eqref{e:palphab}
we can find a unique $p \in (0,  \alpha)$ such that $C_1 = \wt C(\alpha,d, p)$. 
Define
\begin{align}
\label{e:palphape}
\sG_\alpha (p):=\{\kappa \in \sG_\alpha:  \text{ the constant } C_1 \text{ in \eqref{e:kappainf} is } \wt C(\alpha,d, p)\}.
\end{align}
Note that $\sG_\alpha=\cup_{0< p < \alpha} \sG_\alpha (p)$.
We fix a $\kappa \in \sG_\alpha (p)$ and let $Y$  be a Hunt process on $D$ corresponding to the Feynman-Kac semigroup of 
$X^D$ via the multiplicative functional $e^{-\int_0^t \kappa (X^D_s)ds}$,
that is, 
\begin{align*}
\E_x\left[  f(Y_t)\right] = \E_x\left[ e^{-\int_0^t \kappa (X^D_s)ds} f(X^D_t)\right], \quad t\ge 0 , x\in D.
\end{align*}
Since, by Example \ref{e:1}, 
$\kappa(x)dx \in \mathbf{K}_1(D)$,
it follows from Theorem \ref{t:f1} that
 $Y$ has a transition density $q^D(t,x,y)$ with the following estimate
\begin{align}\label{e:qDfe}
q^D(t,x,y) \asymp \P_x(\zeta>t) \P_y(\zeta>t)\left[t^{-d/\alpha} \wedge \frac{t}{|x-y|^{d+\alpha}}\right],
\end{align} for $(t,x,y) \in (0,1) \times D \times D$,
where  $\zeta$ is the lifetime of $Y$.
Moreover, when $C_2=0$, 
$\kappa(x)dx \in \mathbf{K}_{\infty}(D)$ by Example \ref{e:2}.
Thus, by Theorem \ref{t:f2},
\eqref{e:qDfe} holds for all $t>0$.

Define 
\begin{equation*}
Lf(x) := -(-\Delta)^{\alpha/2} f(x)-\kappa(x)f(x).
\end{equation*}

Fix a $q\in (p,\alpha)$ such that $q<p-\eta+\alpha$ and let
$A=\wt C(\alpha, d, q)-\wt C( \alpha,d,  p)>0$. Define
\begin{equation*}
v_1(x):=u_p(x)+u_q(x), \quad v_2(x):=u_p(x)-\frac{1}{2}u_q(x).
\end{equation*}
Since, for $|x|<C_2^{-1}$,  in view of \eqref{rapresth} and \eqref{e:kappainf}, 
\begin{align*}
Lv_1(x) &\ge
A|x|^{q-\alpha}-|x|^p\big|\wt C(d, \alpha, p)|x|^{-\alpha}-\kappa(x)\big|-|x|^q\big|\wt C(d, \alpha, p)|x|^{-\alpha}-\kappa(x)\big|\\
&\ge
A|x|^{q-\alpha}-2C_2(|x|^{p-\eta}+|x|^{q-\eta})
\end{align*}
and
\begin{align*}
Lv_2(x) &\le -2^{-1}A|x|^{q-\alpha}+|x|^p\big|\wt C(d, \alpha, p)|x|^{-\alpha}-\kappa(x)\big|+2^{-1
}|x|^q\big|\wt C(d, \alpha, p)|x|^{-\alpha}-\kappa(x)\big|\\
&\le
-2^{-1} A|x|^{q-\alpha}+(3/2) C_2(|x|^{p-\eta}+|x|^{q-\eta}),
\end{align*}
 there exists $c_1 >0$ such that $Lv_1(x) \ge 0$ and $Lv_2(x)  \le 0$ whenever  $0<|x|<c_1$.
 Pick any $(t,x) \in (0,1) \times D$ and set $r=r(t)=c_1 t^{1/\alpha}$ for $t<1$.
Now we can follow the argument before the statement of Theorem \ref{t:DHKE1} and 
 get 
$\P_x(\zeta>t) \asymp \left(1 \wedge {|x|}/{r}\right)^p$ for $t<1$.

Moreover, if $\kappa(x)=\wt C(\alpha, d, p) |x|^{-\alpha}$, we can simply take $v_1(x)=v_2(x)=u_p(x)$ and $r(t)=t^{1/\alpha}$ for all $t>0$
 and get 
$\P_x(\zeta>t) \asymp \left(1\wedge {|x|}/{r(t)}\right)^p$ for all $t>0$.
 
Therefore, we conclude that 
\begin{theorem}\label{t:DHKE3} 
For all positive $T>0$,   $p \in (0, \alpha)$ and $\eta \in [0, \alpha)$, 
there exists a $c=c(C_1, C_2, p, \alpha, d, \eta, T) \ge 1$ such that for all $\kappa \in \sG_\alpha (p)$,  
the transition density $q(t, x, y)$ of $Y$, 
the Hunt process
on $\R^d\setminus\{0\}$ associated with the Feynman-Kac semigroup of the isotropic $\alpha$-stable process $Z$
via the multiplicative functional
$e^{-\int^t_0\kappa(Z_s)ds}$, 
satisfies that
\begin{align*}
&c^{-1} \left(1 \wedge \frac{|x|}{t^{1/\alpha}}\right)^p \left(1 \wedge \frac{|y|}{t^{1/\alpha}}\right)^p \left[t^{-d/\alpha} \wedge \frac{t}{|x-y|^{d+\alpha}}\right]\\
&\le q(t,x,y) \le c \left(1 \wedge \frac{|x|}{t^{1/\alpha}}\right)^p \left(1 \wedge \frac{|y|}{t^{1/\alpha}}\right)^p \left[t^{-d/\alpha} \wedge \frac{t}{|x-y|^{d+\alpha}}\right], \end{align*}
for $(t,x,y) \in (0, T) \times (\R^d \setminus \{0\}) \times (\R^d \setminus \{0\})$. 
Moreover, if $\kappa(x)=\wt C(\alpha, d, p) |x|^{-\alpha}$, then the above estimates holds for all $t>0$. 
\end{theorem}

The last claim in Theorem \ref{t:DHKE3} can be proved using the scaling property and the finite time estimates in Theorem \ref{t:DHKE3}. 
This was proved independently in \cite{JW} using a different method.

Let $\alpha \in (1, 2)$ and $g$ be an $\R^d$-valued $C^{1}$ function with $\|g\|_\infty+\|\nabla g\|_\infty  <\infty$. 
Let $\wt{X}^g$ be an $\alpha$-stable process with drift $g$, that is, a non-symmetric 
Hunt process with generator $ -(-\Delta)^{\alpha/2} f(x)+ g \cdot \nabla f(x)$, 
see \cite{BJ}.
Let $X^g$ be the Hunt process obtained from $\wt{X}^g$ by killing with rate $\|{\rm div\,} g\|$.
The generator of $X^g$ is 
$ -(-\Delta)^{\alpha/2} f(x)+ g \cdot \nabla f(x)-\|{\rm div} g\|_\infty f(x)$.
By \cite{BJ}, the transition density $p(t,x,y)$ of $X^g$ satisfies
\begin{align}
\label{e:nerw1}
p(t,x,y) \asymp t^{-d/\alpha} \wedge \frac{t}{|x-y|^{d+\alpha}}, \quad (t,x,y) \in (0,1]\times \R^d \times \R^d.
\end{align}

The dual of $ -(-\Delta)^{\alpha/2} f(x)+ g \cdot \nabla f(x)-\|{\rm div} g\|_\infty f(x)$ is
$ -(-\Delta)^{\alpha/2} f(x)- g \cdot \nabla f(x)-{\rm div\,}g(x)f(x) -\|{\rm div} g\|_\infty f(x)$, 
which is the generator of a Hunt process $\wh{X}^g$ 
which can be obtained from an $\alpha$-stable process
with drift via the killing potential $-{\rm div\,}g(x) -\|{\rm div} g\|_\infty$.
It is easy to check that $X^g$ and $\wh{X}^g$ are strong duals 
of each other with respect to 
the Lebesgue measure. 
It is also easy to check that $X^g$ and $\wh{X}^g$ satisfy the sector condition, thus, by
\cite[Theorem 4.17]{Fitz},
all semipolar sets are polar. 
Moreover, since $\alpha \in (1,2)$, 
Assumption $\mathbf{U}$ holds true.

Fix a $\kappa \in \sG_\alpha (p)$ and let $Y^g$  be a Hunt process on $D$ corresponding to the Feynman-Kac semigroup of 
$X^{g, D}$ via the multiplicative functional $e^{-\int_0^t \kappa (X^{g, D}_s)ds}$,
that is, 
\begin{align*}
\E_x\left[  f(Y^g_t)\right] = \E_x\left[ e^{-\int_0^t \kappa (X^{g, D}_s)ds} f(X^{g, D}_t)\right], \quad t\ge 0 , x\in D.
\end{align*}
Since $\kappa(x) dx \in \mathbf{K}_1(D)$ 
by \eqref{e:nerw1}, it follows from 
\eqref{e:nerw1} and Theorem \ref{t:f1} that
 $Y^g$ has a transition density $q^g(t,x,y)$ with the following estimate
\begin{align}\label{e:qDfe2}
q^g(t,x,y) \asymp \P_x(\zeta>t) \P_y(\zeta>t)\left[t^{-d/\alpha} \wedge \frac{t}{|x-y|^{d+\alpha}}\right],
\end{align} for $(t,x,y) \in (0,1] \times D \times D$,
where  $\zeta$ is the lifetime of $Y^g$.
Note that with $u_p=|x|^p$, we get that
 \begin{equation}\begin{split}\label{rapresth1}
| g \cdot \nabla u_p(x)|+\|{\rm div} g\|_\infty |u_p(x)|  \le  \wt C(\alpha, d, p) \, |x|^{p-1}, 
\quad 0<|x|<1.
\end{split}
\end{equation}
From \eqref{rapresth}, \eqref{rapresth1} and the assumption $\alpha \in (1, 2)$, we see that terms
 $ g \cdot \nabla f(x)-\|{\rm div} g\|_\infty f(x)$ and 
$- g \cdot \nabla f(x)-{\rm div\,}g(x)f(x) -\|{\rm div} g\|_\infty f(x)$ can be treated as lower order terms. 
Thus, using \eqref{rapresth1} and the assumption $\alpha \in (1, 2)$, by repeating the argument of the first part of this subsection, we can easily get the following result.

\begin{theorem}\label{t:DHKE-ns} 
Suppose that  $\alpha \in (1, 2)$. 
For all positive $T>0$,   $p \in (0, \alpha)$ and $\eta \in [0, \alpha)$, 
there exists a 
$c=c(C_1, C_2, p, \|g\|_\infty, \alpha, d, \eta, T, \| \nabla g\|_\infty) \ge 1$ 
such that for all $\kappa \in \sG_\alpha (p)$,  
the transition semigroup $q^{g}(t, x, y)$ of $Y^g$ satisfies that
\begin{align*}
&c^{-1} \left(1 \wedge \frac{|x|}{t^{1/\alpha}}\right)^p \left(1 \wedge \frac{|y|}{t^{1/\alpha}}\right)^p \left[t^{-d/\alpha} \wedge \frac{t}{|x-y|^{d+\alpha}}\right]\\
&\le q^g(t,x,y) \le c \left(1 \wedge \frac{|x|}{t^{1/\alpha}}\right)^p \left(1 \wedge \frac{|y|}{t^{1/\alpha}}\right)^p \left[t^{-d/\alpha} \wedge \frac{t}{|x-y|^{d+\alpha}}\right], \end{align*}
for $(t,x,y) \in (0, T) \times (\R^d \setminus \{0\}) \times (\R^d \setminus \{0\})$. 
\end{theorem}

Note that, Theorem \ref{t:DHKE-ns} 
also holds for the fundamental solution to $\partial_t = -(-\Delta)^{\alpha/2} + g \cdot~\nabla  -\kappa(x)$.

\section{Appendix: Continuous additive functionals for killed non-symmetric processes}

We keep the assumptions  and the notations in Sections \ref{s:setup}--\ref{s:3p}.
In this section, $D$ is an open subset of $\X$ and $U$ is a relatively compact subset of $D$.

\begin{lemma}\label{l1}
If $h\in \DD( \wh \LL)$ is non-negative, bounded and has compact support contained in $U$, then for any $t\ge 0$,
$$
\limsup_{\epsilon\to 0}\frac1{\epsilon}\int_U \wh P^U_th(x)\P_x(
\tau^X_U
\le \epsilon)m(dx)<\infty.
$$
\end{lemma}

\pf Noticing $h(\wh X_{\tau^X_U})=0$, we get
$$
\wh P^U_th(x)=\wh \E_x[h(\wh X_t)1_{t<\wh 
\tau^X_U}]=\wh \E_xh(\wh X_{t\wedge
\wh \tau^X_U})=
h(x)+\wh \E_x\int^{t\wedge \wh
\tau^X_U}_0\wh\LL h( \wh X_s)ds.
$$
Using this and the duality,
we have
\begin{align*}
&\int_U\wh P^U_th(x)\P_x(
\tau^X_U\le \epsilon)m(dx) =\int_U \wh P^U_th(x)(1-P^U_\epsilon1(x))m(dx) 
\\&=\int_U(\wh P^U_th(x)- \wh P^U_{t+\epsilon}h(x))m(dx) =-\int_U \wh \E_x\int^{t+\epsilon}_t \wh \LL h(\wh X_s) 1_{s<\wh 
\tau^X_U}dsm(dx) \\
&\le \epsilon \left(\sup_{x\in U}|\wh \LL h(x)|  \right)m(U),
\end{align*}
from which the conclusion follows immediately.
\qed

\

\begin{lemma}\label{l2}
Let $\mu\in  { \bf K}_T(D) $ for some $T>0$. 
If $A$ is the continuous additive functional of 
$X^D$ associated with $\mu$, $h\in \DD(\wh \LL)$ is non-negative, bounded and has compact support contained in $U$, then for any bounded  Borel function $f$
on $U$ and $t\ge 0$,
\begin{align*}
\lim_{\epsilon\to 0}\frac1{\epsilon}\int_U \wh P^U_th(x)\left(\E_x\int^\epsilon_0 {\bf 1}_{\tau^X_U\le s}f(X^D_s)dA_s \right) m(dx)=0.
\end{align*}
\end{lemma}

\pf
Since by the strong Markov property $$
\E_x\int^\epsilon_01_{
\tau^X_U\le s}f(X^D_s)dA_s=\E_x\left[{\bf 1}_{
\tau^X_U<\epsilon}
\E_{X^D_{\tau^X_U}}\int^{\epsilon-\tau^X_U}_0 f(X^D_s)dA_s
\right],
$$
we have 
\begin{align*}
\left|\E_x\int^\epsilon_01_{
\tau^X_U\le s}f(X^D_s)dA_s\right|\le \left(\sup_{y\in 
\X}\E_{y}\int^\epsilon_0|f(X^D_s)|dA_s \right) \P_x(
\tau^X_U\le \epsilon).
\end{align*}
The assertion now follows from Lemma \ref{l1} and 
condition (2) in Definition \ref{d:KT}.
\qed

\begin{proposition}
Let $\mu\in  { \bf K}_T(D)$ for some $T>0$. 
If $A$ is the continuous additive functional of $X$ associated with $\mu$, 
then $(A_{t\wedge\tau^X_U})$
is the continuous additive functional of $X^U$ associated with $\mu_U$.
\end{proposition}

\pf 
Let $A^U_t:=A_{t\wedge
\tau^X_U}$.  Then $A^U$ is a continuous additive functional of $X^U$.  
Let $h\in \DD(\LL)$ be non-negative, bounded and have compact support contained in $U$, and let $f$ be a bounded Borel function supported in $U$. Define
$$
g_t:=\int_Uh(x)\E_x\int^t_0f(X^U_s)dA^U_sm(dx).
$$
Since 
\begin{align*}
&g_{t+\epsilon}-g_t =\int_Uh(x)\E_x\int^{t+\epsilon}_t f(X^U_s)dA^U_sm(dx)
=\int_Uh(x)\E_x\int^{\epsilon}_0 f(X^U_{s+t})dA^U_{s+t}m(dx)\\
&=\int_U h(x) P^U_t(\E_\cdot \int^\epsilon_0 f(X^U_s)dA_s)(x)m(dx)
=\int_U \wh P^U_th(x)\E_x\int^\epsilon_0 {\bf 1}_{\tau^X_U> s} f(X^U_s)dA_sm(dx),
\end{align*}
it follows from Lemma \ref{l2} that
$$
\lim_{\epsilon\to 0}\frac{g_{t+\epsilon}-g_t}{\epsilon}
=\lim_{\epsilon\to 0}\frac1\epsilon\int_U \wh P^U_th(x)
\E_x\int^\epsilon_0f(X^D_s)dA_sm(dx)
=\int_U \wh P^U_th(x)f(x)\mu(dx)
$$
which implies
$$
\int_Uh(x)\E_x\int^t_0f(X^U_s)dA^U_sm(dx)
=\int^t_0\int_U \wh P^U_sh(x)f(x)\mu(dx)ds=\int_U  h(x)\int^t_0 P_s^Uf(x)ds\mu(dx).
$$
Using   the dominated convergence theorem and the monotone convergence theorem, one can show that the equality above is valid for all bounded non-negative Borel functions $h$ and $f$ supported in $U$. 
Therefore, 
$$\int_U f(x) \mu (dx) = \lim_{t\downarrow 0}
\E_{m_U} \left[ \frac1t \int_0^t f(X^U_s) dA^U_s \right].$$ 
\qed

\bigskip
\noindent
 
{\bf Acknowledgments:}  We thank the referee for the very helpful comments.

\vskip 0.1truein

\parindent=0em

{\bf Soobin Cho}

Department of Mathematical Sciences,

Seoul National University, Building 27, 1 Gwanak-ro, Gwanak-gu Seoul 08826, Republic of Korea

E-mail: \texttt{soobin15@snu.ac.kr}

\bigskip

{\bf Panki Kim}

Department of Mathematical Sciences and Research Institute of Mathematics,

Seoul National University, Building 27, 1 Gwanak-ro, Gwanak-gu Seoul 08826, Republic of Korea

E-mail: \texttt{pkim@snu.ac.kr}

\bigskip

{\bf Renming Song}

Department of Mathematics, University of Illinois, Urbana, IL 61801,
USA

E-mail: \texttt{rsong@math.uiuc.edu}

\bigskip

{\bf Zoran Vondra\v{c}ek}

Department of Mathematics, University of Zagreb, Zagreb, Croatia,

Email: \texttt{vondra@math.hr}

\end{document}